\documentclass[aps,twocolumn,prd,superscriptaddress,preprintnumbers,floatfix,showpacs,showkeys,letter]{revtex4-2}

\usepackage{amscd}
\usepackage{amsmath}
\usepackage{amssymb}
\usepackage{graphicx}
\usepackage{scalerel}

\usepackage{dsfont}
\usepackage[ruled]{algorithm2e}
\usepackage{needspace}
\usepackage{newtxtext}
\usepackage{booktabs}
\usepackage{bm}

\usepackage[colorlinks=true, pdfstartview=FitV, linkcolor=blue,
            citecolor=blue, urlcolor=blue]{hyperref} 
\usepackage{subcaption}
\usepackage{sidecap}
\usepackage{grffile}
\usepackage{xcolor}

\graphicspath{{./Figures/},{../Figures/}}

\providecommand{\U}[1]{\protect\rule{.1in}{.1in}}
\newenvironment{proof}[1][Proof]{\noindent\textbf{#1.} }{\ \rule{0.5em}{0.5em}}

\def\W{{\boldsymbol{W}_t}}
\def\bfeta{{\boldsymbol{\eta}_t}}

\newcommand{\norm}[1]{\left\lVert#1\right\rVert}
\def\m{\mathfrak{m}}
\def\Re{\mathrm{Re}} 
\def\Im{\mathrm{Im}} 

\def\bea{\begin{equation} \begin{aligned}}
\def\eea{\end{aligned} \end{equation}}
\def\beas{\begin{equation*} \begin{aligned}}
\def\eeas{\end{aligned} \end{equation*}}
\def\bes{\begin{equation*}}
\def\ees{\end{equation*}}

\def\d{\, \mathrm{d}}

\def\be{\begin{equation}}
\def\ee{\end{equation}}
\def\adots{
  \mathinner{\mkern1mu\raise1pt\hbox{.}\mkern2mu\raise4pt\hbox{.}
  \mkern2mu\raise7pt\vbox{\kern7pt\hbox{.}}\mkern1mu}}

\def\bftau{{\boldsymbol{\tau}}}

\def\c{\mathfrak{c}}
\def\s{\mathfrak{f}}

\def\L{\mathcal{L}}
\def\B{\mathcal{B}}

\def\P{\mathcal{P}}

\newcommand{\x}{\bm{x}}
\newcommand{\y}{\bm{y}}

\newtheorem{thm}{Theorem}[section]
\newtheorem{lem}{Lemma}[section]
\newtheorem{defi}{Definition}[section]

\newtheorem{rem}{Remark}[section]
\newtheorem{cor}{Corollary}[section]
\def\bt{\begin{thm}}
\def\et{\end{thm}}
\def\bl{\begin{lem}}
\def\el{\end{lem}}
\def\bd{\begin{defi}}
\def\ed{\end{defi}}
\def\bc{\begin{cor}}
\def\ec{\end{cor}}
\def\bp{\begin{proof}}
\def\ep{\end{proof}}
\def\br{\begin{rem}}
\def\er{\end{rem}}
\def\bi{\begin{itemize}}
\def\ei{\end{itemize}}

\begin{document}

\bibliographystyle{apsrev}

\title[]{Non-Markovian Reduced Models to Unravel Transitions in Non-equilibrium Systems}

\author{Micka{\"e}l D. Chekroun}
\email{mchekroun@atmos.ucla.edu} 
\affiliation{Department of Earth and Planetary Sciences, Weizmann Institute, Rehovot 76100, Israel} 
\affiliation{Department of Atmospheric and Oceanic Sciences, University of California, Los Angeles, CA 90095-1565, USA}

\author{Honghu Liu}
\affiliation{Department of Mathematics, Virginia Tech, Blacksburg, VA 24061, USA}

\author{James C. McWilliams}
\affiliation{Department of Atmospheric and Oceanic Sciences and Institute of Geophysics and Planetary Physics, University of California, Los Angeles, CA 90095-1565, USA}

\date{\today}

\begin{abstract}
This work proposes a general framework for analyzing noise-driven transitions in  spatially extended non-equilibrium systems and explains the emergence of coherent patterns beyond the instability onset. The framework relies on stochastic parameterization formulas to reduce the complexity of the original equations while preserving the essential dynamical effects of unresolved scales. The approach is flexible and operates for both Gaussian noise and non-Gaussian noise with jumps.

Our stochastic parameterization formulas offer two key advantages. First, they can approximate stochastic invariant manifolds when these manifolds exist. Second, even when such manifolds break down, our formulas can be adapted through a simple optimization of its constitutive parameters. This allows us to handle scenarios with weak time-scale separation where the system has undergone multiple transitions, resulting in large-amplitude solutions not captured by invariant manifolds or other time-scale separation methods.

The optimized stochastic parameterizations capture then how small-scale noise impacts larger scales through the system's nonlinear interactions. This effect is achieved by the very fabric of our parameterizations incorporating non-Markovian (memory-dependent) coefficients into the reduced equation. These coefficients account for the noise's past influence, not just its current value, using a finite memory length that is selected for optimal performance. The specific ``memory" function, which determines how this past influence is weighted, depends on both the strength of the noise and how it interacts with the system's nonlinearities.

Remarkably, training our theory-guided reduced models on a single noise path effectively learns the optimal memory length for out-of-sample predictions. This approach retains indeed good accuracy in predicting noise-induced transitions, including rare events, when tested against a large ensemble of different noise paths. This success stems from our ``hybrid" approach, which combines analytical understanding with data-driven learning. This combination avoids a key limitation of purely data-driven methods: their struggle to generalize to unseen scenarios, also known as the ``extrapolation problem."

\end{abstract}




\maketitle

\tableofcontents

\bibliographystyle{aipnum4-1}

\section{Introduction}
Non-equilibrium systems are  irreversible systems characterized by a continuous flow of energy that is driven by external forces or internal fluctuations. There are many different types of non-equilibrium systems, and they can be found in a wide variety of fields, including physics, chemistry, biology, and engineering.
Non-equilibrium systems can exhibit complex behavior, including self-organization, pattern formation, and chaos.
These complex behaviors arise from the interplay of nonlinear dynamics and statistical fluctuations.
The emergence of coherent, complex patterns is ubiquitous in many spatially extended non-equilibrium systems; see e.g.~\cite{hohenberg1977theory,pearson1993complex,cross1993pattern,berloff1999large,murray2003mathematical,baurmann2007instabilities,sagues2007spatiotemporal,zelnik2015gradual,meron2018patterns,mcwilliams2019survey}. Mechanisms to explain the emergence of these patterns include the development of instability saturated by nonlinear effects whose calculations can be conducted typically near the onset of linear instability; see e.g.~\cite{crawford1991introduction,ma2019phase}. 
Historically, hydrodynamic systems have been commonly used as
prototypes to study instabilities out of equilibrium in
spatially extended systems, from both an experimental
and a theoretical point of view \cite{crawford1991symmetry,cross1993pattern}.

To describe how physical instability develops in such systems with infinitely many degrees of freedom, we often focus on the amplitudes' temporal evolution of specific normal modes. The latter correspond typically to those that are mildly unstable and that are only slightly damped in linear theory. When the number of these nearly marginal modes is finite, their amplitudes are governed by ordinary differential equations (ODEs) in which the growth rates of the linear theory have been renormalized by nonlinear terms \citep{hohenberg1977theory,crawford1991introduction,cross1993pattern,ma2019phase}.
Intuitively, the reason for this reduction is a simple separation of time scales. Modes that have just crossed the imaginary axis have a small real part and are evolving slowly 
on long time scales, all the other fast modes rapidly adapting themselves to these slow modes.

However, for non-equilibrium systems away from the onset of linear instability, such as when the Reynolds number for a fluid flow increases far beyond a laminar regime, the emergence of coherent patterns does not fit within this instability/nonlinear saturation theory as the reduction principle of the fast modes onto the slow ones breaks down \cite{crawford1991introduction} and calls for new reduction techniques to explain emergence. 
Furthermore, in the presence of random fluctuations, reduced equations under the form of deterministic ODEs are inherently incapable to capture phenomena like noise-induced transitions.  Noise can have unexpected outcomes such as altering the sequence of transitions, possibly suppressing instabilities \cite{Chekroun_al2016}, the onset of turbulence \cite{brand1985external} or, to the opposite,  exciting transitions \cite{kai1989structure,sagues2007spatiotemporal}. Examples in which noise-induced transitions 
have an important role include the generation of convective rolls  \cite{ahlers1989deterministic,garcia1993effects}, electroconvection in nematic liquid crystals \cite{kai1989structure,rehberg1991thermally}, certain oceanic flows  \cite{sura2002sensitivity,simonnet2021multistability}, and climate phenomena \cite{saravanan1998advective,roulston2000response,eisenman2005westerly,Margazoglou2021}. 

This work aims to gain deeper insights into the efficient derivation of reduced models for such non-equilibrium systems subject to random fluctuations \cite{cross1993pattern,sagues2007spatiotemporal}. In that respect, we seek reduced models able to capture the emergence of noise-driven spatiotemporal patterns and predict their transitions triggered by the subtle coupling between noise and the nonlinear dynamics. 
Our overall goal is to develop a unified framework for identifying the key variables and their interactions (parameterization) in complex systems undergoing noise-driven transitions. We specifically focus on systems that have experienced multiple branching points (bifurcations), leading to high-multiplicity regimes with numerous co-existing metastable states whose amplitude is relatively large (order one or more) resulting from a combination of strong noise, inherent nonlinear dynamics, or both. The purpose is thus to address limitations of existing approaches like amplitude equations and center manifold techniques in dealing with such regimes. 

To address these limitations, the framework of optimal parameterizing manifolds (OPMs) introduced in \citep{CLM17_Lorenz9D,CLM20_closure,chekroun2023optimal} for forced-dissipative systems offers a powerful solution. This framework allows for the efficient derivation of parameterizations away from the onset of instability, achieved through continuous deformations of near-onset parameterizations that are optimized by data-driven minimization. Typically, the loss function involves a natural  discrepancy metric measuring the defect of parameterization of the unresolved variables by the resolved ones.
The parameterizations that are optimized are derived from the governing equations
by means of backward-forward systems providing useful approximations of the unstable growth saturated by the nonlinear terms, even when the spectral gap between the stable and unstable modes is small \cite{CLM20_closure,CLW15_vol2,chekroun2023optimal}. These approximations, once optimized, do not suffer indeed  the restrictions of invariant manifolds (spectral gap condition \cite[Eq.~(2.18)]{CLM20_closure}) allowing to handle situations with e.g.~weak time-scale separation between the resolved and unresolved modes.

In this article, we carry over this variational approach to the case of spatially extended non-equilibrium systems within the framework of stochastic partial differential equations (SPDEs). 
In particular, we extend the parameterization formulas obtained in \citep{CLW15_vol1,CLW15_vol2,Chekroun_al2023} for SPDEs driven by multiplicative noise (parameter noise \cite{garcia1993effects,sagues2007spatiotemporal}), to SPDEs driven, beyond this onset, by more realistic, spatio-temporal noise either of Gaussian or of non-Gaussian nature with jumps. For this class of SPDEs, our framework allows for dealing with  the important case of cutoff scales larger than the scales forced stochastically.

 For this forcing scenario, the stochastic parameterizations  derived in this work, give rise to reduced systems taking the form of stochastic differential equations (SDEs) with non-Markovian, path-dependent coefficients depending on the noise's past; see e.g.~Section \ref{Sec_above_forcingscaleACE} below.
We mention that many time series records from real-world datasets, are known to exhibit long-term memory. This is the case of long-range memory features assumed to represent the internal variability of the climate on time scales from years to centuries \cite{schulz2002redfit,moberg2005highly,fraedrich2009continuum,rypdal2014long}. 
There, the surface temperature is considered as a superposition of internal variability (the response to stochastic forcing) and a forced signal which is the linear response to external forcing. The background variability may be modeled using a stochastic process with memory, or a different process that incorporates non-Gaussianity if this is considered more appropriate \cite{rypdal2016late}.
As will become apparent at the end of this paper, our reduction approach is adaptable to PDEs under the influence of stochastic forcing with long-range dependence; see Section \ref{Sec_discussion}.

 In parallel and since Hasselmann's seminal work \cite{hasselmann1976stochastic},  several climate models have been incorporating stochastic forcing to represent unresolved processes \cite{lin2000influence,palmer2001nonlinear,majda2001mathematical,Franzke.ea.2015,berner2017stochastic,LC2023}. Our approach shows that reduced models derived from these complex systems are expected to exhibit finite-range memory effects depending on the past history of the stochastic forcing, even if the latter is white in time.
These exogenous memory effects, stemming from the forcing noise's history, are distinct from endogenous memory effects encountered in the reduction of nonlinear, unforced systems \cite{Chorin2002,GKS04,LC2023}. The latter, predicted by the Mori-Zwanzig (MZ) theory, are functionals of the past of the resolved state variables. Endogenous memory effects often appear when the validity of the conditional expectation in the MZ expansion, breaks down  \cite{Stinis06,CLM17_Lorenz9D,Chekroun2021c,LC2023}.
The exogenous memory effects dealt with in this work, have a different origin, arising as soon as the cutoff scale is larger than the scales forced stochastically.

It is worth mentioning that reduced models with non-Markovian coefficients responsible for such exogenous memory effects, have been encountered in the reduction of stochastic systems driven by white noise, albeit near a change of stability. These are obtained from reduction approaches benefiting from timescale separation such as stochastic normal forms \cite{NL91,Coullet_al_Hopf,Arnold98},  stochastic invariant manifold approximations \citep{FS09,CLW15_vol1,CLW15_vol2,Chekroun_al2023}, or multiscale methods \cite{BH05,Blom07,BHP07,klepel2014amplitude}. Our reduction approach is not limited to timescale separation.

 In that respect, our reduction framework is tested against a stochastic Allen-Cahn equation (sACE), a powerful tool for modeling non-equilibrium phase transitions in various scientific fields \cite{hohenberg1977theory,cross1993pattern,bray2002theory,tauber2014critical}.  The latter model is set in a parameter regime with weak timescale separation and far from the instability onset, in which the system exhibits multiple coexisting metastable states connected to the basic state through rare or typical stochastic transition paths.

We show that our resulting OPM reduced systems, when trained over a single path, are not only able to retain a remarkable predictive power to  emulate ensemble statistics (correlations, power spectra, etc.) obtained as average over a large ensemble of noise realizations,  but also anticipate what are the system's typical and rare metastable states and their statistics of occurrence. The OPM reduced system's ability to reproduce accurately such transitions is rooted in its very structure.
Its coefficients are nonlinear functionals of the aforementioned non-Markovian terms (see  Eq.~\eqref{Eq_closure} below) allowing for an accurate representation of the genuine nonlinear interactions between noise and nonlinear terms in the original sACE, which drive fluctuations in the large-mode amplitudes.

We emphasize, that because the cutoff scale, defining the retained large-scale dynamics, is chosen here to be larger than the forcing scale, any deterministic reduction (e.g.~Galerkin, invariant manifold, etc.) would filter out the underlying noise's effects. This leaves them blind to the subtle fluctuations driving the system's behavior and leading eventually to stochastic transitions. In contrast, our stochastic parameterization approach sees right through this filter. It tracks how the noise, acting in the ``unseen" part of the system, interacts with the resolved modes through our reduction method's non-Markovian coefficients.

To demonstrate the broad applicability of our reduction approach, we apply it in Section \ref{Sec_Jump} to another significant class of spatially extended non-equilibrium systems: jump-driven SPDEs. 
The unforced, nonlinear dynamics is here chosen to exhibit an S-shaped solution curve. Many systems sharing this attribute are indicative of multistability, tipping points, and hysteresis. Such behaviors are observed in various fields, including combustion theory  \cite{bebernes2013mathematical,frank2015diffusion}, plasma physics  \cite{grad1958hydromagnetic,Shafranov1958}, ecology \cite{meron2018patterns}, neuroscience \cite{izhikevich2007dynamical}, climate science \cite{ghil1976climate,north1981energy}, and oceanography \cite{thual1992catastrophe,weijer2019stability}. By understanding the complex interactions between noise and nonlinearity into these phenomena, we can gain insights into critical transitions and tipping points, which are increasingly relevant to addressing global challenges like climate change \cite{lenton2008tipping,boers2021,boers2021critical,wunderling2021interacting,LC2023}.

Jump processes offer a powerful tool for modeling complex systems with non-smooth dynamics.
They offer valuable insights into multistable behaviors \cite{serdukova2017metastability,zheng2020maximum} and have been applied to various real-world phenomena, such as paleoclimate events  \cite{ditlevsen1999observation,rypdal2016late}, chaotic transport \cite{solomon1993observation}, atmospheric dynamics \cite{khouider2003coarse,stechmann2011stochastic,penland2012alternative,thual2016simple,chen2017simple,chen2024stochastic}, and cloud physics \cite{horenko2011nonstationarity,Chekroun_al2022_SciAdv}.

SPDEs driven by jump processes are gaining increasing attention in applications \cite{debussche2013dynamics}. 
For instance, jump processes can replace (non-smooth) "if-then" conditions in such models, enabling more efficient simulations \cite{chen2024stochastic}. Nevertheless, efficient reduction techniques to disentangle the jump interactions with other  smoother nonlinear components of the model, are still under development \cite{yuan2019slow,liu2021random,yuan2022modulation}. Our stochastic parameterization framework provides a promising solution in this direction, as demonstrated in Section \ref{Sec_Jump}. By capturing the intricate interplay between jump noise and nonlinear dynamics within the reduced equations, our approach can significantly simplify the analysis of such complex systems.


\section{Invariance Equation and Approximations}\label{Sec_Invariance_Eq}

\subsection{Spatially extended non-equilibrium systems}\label{Sec_SPDEs}
This article is concerned with the efficient reduction of spatially extended non-equilibrium systems. 
To do so, we work within the framework of Stochastic Equations in Infinite Dimensions \cite{DPZ08} and its Ergodic Theory  \cite{DZ96}. Formally, these equations take the following  form        
\be \label{Eq_SPDE}
\d u= (A u + G(u))\d t + \d \bfeta,
\ee
 in which $\bfeta$ is a stochastic process, either composed of Brownian motions or jump processes and whose exact structure is specified below. This formalism provides a convenient way to analyze stochastic partial differential equations (SPDEs) by means of semigroup theory \cite{Hen81,Pazy83}, in which the unknown $u$ evolves typically in a Hilbert space $H$.

The operator $A$ represents a linear differential operator while  $G$ is a nonlinear operator that accounts for the nonlinear terms. Both of these operators may involve loss of spatial regularity when applied to a function $u$ in $H$. To have a consistent existence theory of solutions and their stochastic invariant manifolds for SPDEs requires to take into account such loss of regularity effects \cite{DPZ08,DZ96}.  

General assumptions encountered in applications  are made on the linear operator $A$ following \cite{Hen81}. More specifically, we assume
\be\label{Eq_A}
A=-\L + \B, 
\ee
where  $\L$ is sectorial with domain $D(\L)\subset H$ which is compactly and densely embedded in $H$. We assume also that $-\L$ is stable, while $\B$ is a low-order perturbation of $\L$, i.e.~$\B$ is a  bounded linear operator such that $\B: D(\L^{\alpha}) \rightarrow H$  for some $\alpha$ in $[0,1)$. We refer to \cite[Sec.~1.4]{Hen81} for an intuitive presentation of fractional power of an operator and characterization of the so-called operator domain $D(\L^{\alpha})$.
In practice, the choice of $\alpha$ should match the loss of regularity effects caused by the nonlinear terms so that 
\bes
G \colon  D(\L^{\alpha})  \rightarrow H,
\ees
is a well-defined $C^p$-smooth mapping that satisfies $G(0)=0$, $DG(0)=0$ (tangency condition), and  
\be\label{F_Taylor}
G(u)= G_k(u, \cdots, u) + O(\|u\|^{k+1}_\alpha), 
\ee
with $p > k \ge 2$, and $G_k$ denoting the leading-order operator (on $D(\L^{\alpha})$) in the Taylor expansion of $G$, near the origin.  A broad class of spatially extended stochastic equations from physics can be recasted into this framework (see \cite{MW05,DPZ08,CLW15_vol1,Lord2014introduction,Chekroun_al2023} for examples), as well as time-delay systems subject to stochastic disturbances \cite{Chekroun_al2022_SciAdv,chekroun2024effective}.

Throughout this article, we assume that the driving noise in Eq.~\eqref{Eq_SPDE} is an   $H$-valued stochastic process that takes the form 
\be \label{Eq_eta}
\bfeta(\omega) = \sum_{j=1}^{N} \sigma_j  \eta_t^j(\omega)\boldsymbol{e}_j, \; t \in \mathbb{R}, \; \sigma_j\geq 0,\;  \omega \in \Omega,
\ee
where the $\bm{e}_j$ are eigenmodes of $A$, the $\eta_t^j$ denote either a finite  family of mutually independent jump processes, or mutually independent Brownian motions, $W_t^j$, over their relevant probability space $(\Omega,\mathbb{P},\mathcal{F})$ endowed with its canonical probability measure $\mathbb{P}$ and filtration $\mathcal{F}$; see  e.g.~\cite[Appendix A.3]{Arnold98}.

Within this framework, the reduction of spatially extended non-equilibrium systems is organized in terms of resolved and unresolved spatial scales.  The resolved scales are typically spanned by large wavenumbers and the unresolved by smaller ones. In that respect, the eigenmodes of the operator $A$ plays a central role throughout this paper to rank the spatial scales.  Typically, we assume that the state space $H_\c$ of resolved variables is spanned by the following modes
\be \label{Hc}
H_{\c}=\mbox{span}\{\boldsymbol{e}_1,\cdots,\boldsymbol{e}_{m_c}\},
\ee
where the $\boldsymbol{e}_j$ correspond to large-scale modes up to a cutoff-scale associated with some index $m_c$ (see Remark \ref{rmk_wave_vectors1}).
The subspace of unresolved modes is then the orthogonal complement of $H_{\c}$ in $H$, namely
\begin{equation} \label{H_decomposition}
\begin{aligned}
H_{\c}\oplus H_{\s} =H.
\end{aligned}
\end{equation}
To these subspaces, we associate their respective canonical projectors denoted by
\bea \label{Pc Ps}
\Pi_{\c}: H \rightarrow H_{\c}, \textrm{ and }  \Pi_{\s}: H \rightarrow H_{\s}.
\eea 
Throughout Section \ref{Sec_OPMwhite} below and the remaining of this Section, we focus on the case of driving Brownian motions; the case of  driving jump processes is dealt with in Section \ref{Sec_Jump}.

Our goal is to provide a general approach to derive reduced models that preserve the essential features of the large-scale dynamics without resolving the small scales.

 We present hereafter and in Section \ref{Sec_OPMwhite} below, the  formulas of the underlying small-scale parameterizations   in the Gaussian noise case. The formalism is flexible and easily adaptable to the case of non-Gaussian  noises with jumps such as discussed in Sections \ref{Sec_Jump} and \ref{Sec_Levy}, below.

\br \label{rmk_wave_vectors1}
Within our working assumptions, the spectrum $\sigma(A)$ consists only of isolated eigenvalues with finite multiplicities. This combined with the sectorial property of $A$ implies that  there are at most finitely many eigenvalues with a given real part. The sectorial property of $A$ also implies that the real part of the spectrum, $\Re( \sigma(A))$, is bounded above (see also \cite[Thm.~II.4.18]{EN00}).  These two properties of $\Re (\sigma(A))$  allow us in turn to label elements in $\sigma(A)$ according to the lexicographical order which we adopt throughout this article:
\bea  \label{eq:ordering-1}
\sigma(A) = \{\lambda_n \:|\:  n \in \mathbb{N}^\ast\}, 
\eea
such that for any $1\le n < n'$ we have either 
\bea
\Re (\lambda_{n}) > \Re( \lambda_{n'}), 
\eea
or
\bea  \label{eq:ordering-3}
\Re (\lambda_{n})= \Re (\lambda_{n'}), \quad \text{ and } \quad \Im (\lambda_{n}) \geq \Im (\lambda_{n'}). 
\eea
This way, we can rely on a simple labeling of the eigenvalues/eigenmodes by positive integers to organize the resolved and unresolved scales and the corresponding parameterization formulas derived hereafter.  In practice, when the spatial domain is 2D or 3D, it is usually  physically more intuitive to label the eigenelements by wave vectors. The parameterization formulas presented below can be easily recast within this convention; see e.g.~\cite{chekroun2022transitions}. 
\er

\subsection{Invariance equation and backward-forward systems} \label{sec_inv_eq}
As mentioned in Introduction, this work extends the parameterization formulas in \cite{CLM20_closure,chekroun2023optimal} (designed for constant or time-dependent forcing) to the stochastic setting. In both deterministic and stochastic contexts, deriving the relevant parameterizations and reduced systems relies on solving Backward-Forward (BF) systems arising in the theory of invariant manifolds. Similar to \cite{CLM20_closure,chekroun2023optimal}, our BF framework allows us to overcome limitations of traditional invariant manifolds, particularly the spectral gap condition (\cite[Eq.~(2.18)]{CLM20_closure}). This restrictive condition, which requires large gaps in the spectrum of the operator $A$, is bypassed through data-driven optimization (detailed in Section \ref{Sec_Optim}) of the backward integration time over which the BF systems are integrated.

Before diving into the parameterization formulas retained to address these limitations (Section \ref{Sec_fully_coupled_BF_gen}), we recall the basics of reducing an SPDE to a Stochastic Invariant Manifold (SIM). In particular, this detour enables us to better appreciate the emergence of BF systems as a natural framework  for parameterizations.
For that purpose, adapting the approach from \citep{CLW15_vol1}, we transform the SPDE reduction problem into the more tractable problem of reducing a PDE with random coefficients. This simplification makes the problem significantly more amenable to analysis than its original SPDE form.

 To do so, consider the stationary solution $z(t,\omega)$ to  the Langevin equation 
\be
\d z = A z \d t+  \d \W,
\ee
(Ornstein-Uhlenbeck (OU) process), where $\W$ is defined in Eq.~\eqref{Eq_eta} with the corresponding mutually independent Brownian motions, $W_t^j$, in place of the $\eta_t^j$, i.e.:
\be \label{Eq_W}
\W(\omega) = \sum_{j=1}^{N} \sigma_j  W_t^j(\omega)\boldsymbol{e}_j, \; t \in \mathbb{R}, \;  \omega \in \Omega.
\ee

Then, for each noise's path $\omega$, the change of variable 
\be 
v = u -  z(t, \omega),
\ee
transforms the SPDE,
\bes
\d u= (A u + G(u))\d t + \d \boldsymbol{W}_t, 
\ees
into the following PDE with random coefficients in the $v$-variable: 
\be \label{RPDE}
\frac{\d v}{\d t} = A v + G(v  + z(t, \omega)). 
\ee

Under standard conditions involving the spectral gap between the resolved and unresolved modes \cite{BF95}, the  path-dependent PDE \eqref{RPDE}  admits a SIM, and the underlying stochastic SIM mapping, $\phi (t,X)$ ($X\in H_\c$), satisfies the stochastic invariance equation:   
 \begin{widetext}
 \begin{equation} \label{RIE}
 \partial_t \phi +\mathcal{L}_A[\phi] (X) = - D \phi(t,X,\omega) \Big(\Pi_{\c} G(X + \phi + z(t, \omega))\Big) + \Pi_{\s} G(X + \phi + z(t, \omega)),
\end{equation}
\end{widetext}
where $\mathcal{L}_A$ denotes the operator acting on differentiable mappings $\psi$ from $H_\c$ into $H_\s$, as follows:
\be \label{h1_eqn}
\mathcal{L}_A[\psi] (X)=D \psi(X) A_{\c} X  - A_{\s} \psi(X), \; X\in H_\c,
\ee
with $A_{\c} = \Pi_{\c} A$ and $A_{\s} = \Pi_{\s} A$.

To simplify, assume that $H_\c$ and $H_\s$  are chosen such that $\Pi_\s \W=0$. Then in particular $z_{\s}=\Pi_{\s} z=0$. Consider now the Lyapunov-Perron integral
\be \label{Eq_phi_stationary}
\mathfrak{J}(t,X,\omega) = \int_{-\infty}^t e^{(t-s)A_{\s}} \Pi_{\s} G_k(e^{(s-t)A_{\c}}X + z_{\c}(s,\omega)) \d s,
\ee
 where $z_{\c}=\Pi_{\c} z$.
 In what follows, we denote  by $\langle \cdot, \cdot \rangle$ the Hermitian  inner product on $H$ defined as  $\langle f, g \rangle =\int f(x) \overline{g}(x) \d x$,  while $G_k$ denotes the leading-order term of order $k$ in the Taylor expansion of $G(X)$ around $X=0$, and $
\bm{e}_n^{\ast}$ denotes the $n$-th eigenmode of the adjoint operator of $A$ adopting the same labelling convention as for $A$; see Remark.~\ref{rmk_wave_vectors1}. 
To the order $k$ we associate  the set $K$ of indices $\{1,\cdots,k\}$, and denote by $K_z$ any subset of indices made of  (possibly empty) of disjoint elements of $K$.

We have then the following result. 
\bt \label{Lemma_BF}
Under the previous assumptions, assume furthermore that the following  non-resonance condition holds, for any $(j_1,\cdots, j_k)$ in $(1,\cdots,m_c)^k$, $n \ge m_c+1$, and any subset $K_z$ such as defined above: 
\bea \label{Eq_NR}
&\left( G_{j_1 \cdots j_k}^n \neq 0 \text{ and } \Pi_{q \in K_z} \sigma_q \neq 0 \right)  \implies\\
&\hspace{3cm}\left(\Re\Big(\lambda_n-\sum_{p \in K \setminus K_z} \lambda_{j_p}\Big) < 0\right),
\eea
with $G_{j_1 \cdots j_k}^n=\langle G_k(\boldsymbol{e}_{j_1},\cdots,\boldsymbol{e}_{j_{k}}),\boldsymbol{e}_n^{\ast}\rangle$.
 
 Then, the Lyapunov-Perron integral $\mathfrak{J}$ (Eq.~\eqref{Eq_phi_stationary}) is well-defined almost surely,  and is a solution to the following homological equation  with random coefficients:
 \be\label{Eq_homoligical}
\partial_t \psi +\mathcal{L}_A[\psi] (X) = \Pi_{\s}G_k(X+z_{\c}(t,\omega)).
\ee

 Moreover, we observe that the integral $\mathfrak{J}$  is (formally) obtained as the limit, when $\tau$ goes to infinity, of the $q$-solution to the following backward-forward (BF) auxiliary system: 
\begin{subequations} \label{Eq_BF_RPDE}
\begin{align}
& \frac{\d p}{\d t} =  A_\c p, \hspace{3.5cm} s \in [t-\tau, t],    \label{BF1_RPDE} \\
& \frac{\d q}{\d t} =  A_\s q + \Pi_{\s} G_k(p+ z_{\c}(s,\omega)),\hspace{.2cm} s \in [t-\tau, t], \label{BF2_RPDE}\\
& \mbox{with } p(s)\vert_{s=t} = X \in H_{\c}, \mbox{ and } q(s)\vert_{s=t-\tau}=0. \label{BF3_RPDE}
\end{align}
\end{subequations}
That is
\be \label{Eq_PB_limit}
 \lim_{\tau \rightarrow \infty} \| \mathfrak{J}(t,X,\omega) - q_{\tau}(t,X, \omega)\|_{H_\s} = 0,
\ee
where $q_{\tau}(t,X, \omega)$ denotes the solution to Eq.~\eqref{BF2_RPDE} at time $s = t$, when initialized with $q=0$ at $s=t-\tau$. 

\et
For a proof of this Theorem, see Appendix \ref{Sec_Proof}.  This theorem extends to the stochastic context Theorem III.1 of \cite{chekroun2023optimal}.  To simplify the notations, we omit below the $\omega$-dependence in certain notations unless specified otherwise.

Theorem \ref{Lemma_BF} teaches us that solving the BF system \eqref{Eq_BF_RPDE} gives  the solution to the homological equation \eqref{Eq_phi_stationary} which is an approximation of the full invariance equation Eq.~\eqref{RIE}. In the deterministic setting,  solutions to the  homological equation \eqref{Eq_phi_stationary} are known to provide actual approximations of the invariant manifolds with rigorous estimates; see \cite[Theorem 1]{CLM20_closure} and \cite[Theorem III.1]{chekroun2023optimal}.  Such results extend to the case of SPDEs driven by multiplicative (parameter) noise; see \cite[Theorem 2.1]{Chekroun_al2023}  and \cite[Theorem 6.1 and Corollary 7.1]{CLW15_vol1}. It is not the scope of this article to deal with rigorous error estimates regarding the approximation problem of stochastic invariant manifolds for SPDEs driven by additive noise, but the rationale stays  the same:  solutions to the  homological equation \eqref{Eq_phi_stationary} or equivalently to the BF system   \eqref{Eq_BF_RPDE} provide actual approximations of the underlying stochastic invariant manifolds, in the additive noise case as well.  

Going back to the SPDE variable $u$, the BF system \eqref{Eq_BF_RPDE} becomes 
\begin{subequations} \label{Eq_BF_SPDE_lowb}
\begin{align}
& \d \widehat{p} =  A_\c \widehat{p} \d s+\Pi_\c \d \bm{W}_s, \hspace{1.75cm} s \in [t-\tau, t],    \label{BF1_RPDEb} \\
& \d \widehat{q} =  \left( A_\s \widehat{q} + \Pi_{\s} G_k(\widehat{p})\right) \d s,  \hspace{1.3cm} s \in [t-\tau, t], \label{BF2_RPDEb}\\
& \mbox{with } \widehat{p}(s)\vert_{s=t} = X \in H_{\c}, \mbox{ and } \widehat{q}(s)\vert_{s=t-\tau}=0. \label{BF3_RPDEb}
\end{align}
\end{subequations}

Note that in Eq.~\eqref{Eq_BF_SPDE_lowb}, only the low-mode variable $\widehat{p}$ is forced stochastically, and thus from what precedes, $\lim_{\tau\rightarrow\infty} \widehat{q}_{\tau}(t,X)$, provides a legitimate approximation of the SIM at time $t$ when in the original SPDE, only the low modes are forced stochastically, i.e.~$N=m_c$ in \eqref{Eq_W}.  
In the next section, we consider the general case when the low and high modes are stochastically forced.

\subsection{Approximations of fully coupled backward-forward systems}\label{Sec_fully_coupled_BF_gen}
Backward-forward systems have been actually proposed in the literature for the construction of SIMs \cite{DPD96}, through a different route than what presented above, i.e.~without exploiting the invariance equation.
The idea pursued in    \cite{DPD96} is to envision SIMs as a   fixed point of  an integral form of the following fully coupled BF system   
\begin{subequations} \label{Eq_BF_full_SPDE}
\begin{align}
& \d p =  \big[ A_\c p  + \Pi_{\c} G(p+q) \big] \mathrm{d} s +  \Pi_{\c} \mathrm{d} \boldsymbol{W}_{s}, \; \; s \in I_{t,\tau},    \label{BF1_full_SPDE} \\
& \d q =  \big[ A_\s q  + \Pi_{\s} G(p+q) \big] \mathrm{d} s +  \Pi_{\s} \mathrm{d} \boldsymbol{W}_{s}, \;\; s \in I_{t,\tau}, \label{BF2_full_SPDE}\\
& \mbox{with } p(s)\vert_{s=t} = X \in H_{\c}, \mbox{ and } q(s)\vert_{s=t-\tau}=0. \label{BF3_full_SPDE}
\end{align}
\end{subequations}
where $I_{t,\tau}=[t-\tau,t]$. Note that in Eq.~\eqref{Eq_BF_full_SPDE}, we do not assume here the noise term $\bm{W}$ to be a finite sum of independent Brownian motions. In the case of an infinite sum  one can thus covers the case of space-time white noise due to      \cite{jetschke1986equivalence}.

In any case, it is known that this type of {\it fully coupled nonlinear problems} involving backward stochastic equations does not have always solutions in general and we refer to \cite[Proposition 3.1]{DPD96} for conditions on $A$ and $G$ ensuring existence of solutions and thus SIM. This (nonlinear) BF approach to SIM is also subject to a spectral gap condition that requires gap to be large enough as in  \cite{BF95}.

Denoting by $q_{\tau}(t,X)$ the solution to Eq.~\eqref{BF2_full_SPDE} at $s = t$,  Proposition 3.4 of \cite{DPD96} ensures then that $\lim_{\tau \rightarrow \infty} q_{\tau}(t,X)$ exists in $L^2(\Omega)$ and that this limit gives the sought SIM, i.e.
\be \label{Eq_SIM}
\Psi^{\text{sim}}(X,t) = \lim_{\tau \rightarrow \infty} q_\tau  (t,X), \quad X \in H_{\c}.
\ee
The SIM is thus obtained, almost surely, as the asymptotic graph of  the mapping $X\mapsto  q_\tau (t, X)$, when $\tau$ is sent to infinite.

Instead of relying on an integral form of  Eq.~\eqref{Eq_BF_full_SPDE}, one can address though the existence and construction of a SIM via a more direct approach, exploiting iteration schemes built directly from  Eq.~\eqref{Eq_BF_full_SPDE}.
It is not the purpose of this article to analyze the convergence of such iterative schemes (subject also to a spectral gap condition as in  \cite{DPD96}) but rather to illustrate how informative such schemes can be in designing stochastic parameterizations in practice.

 To solve \eqref{Eq_BF_full_SPDE}, given an initial guess, $(p^{(0)},q^{(0)})$, we propose the following iterative scheme:
\begin{subequations} \label{Eq_BF_fully_coupled_additive2}
\begin{align}
& \mathrm{d} p^{(\ell)} =  \Big[ A_\c p^{(\ell)} + \Pi_{\c} G\big( p^{(\ell-1)} +  q^{(\ell-1)} \big) \Big] \mathrm{d} s+  \Pi_{\c} \mathrm{d} \boldsymbol{W}_{s},   \label{Eq_BF2_DPDa} \\
& \mathrm{d} q^{(\ell)} = \Big[ A_{\s} q^{(\ell)}  +  \Pi_{\s} G\big( p^{(\ell-1)} +  q^{(\ell-1)} \big) \Big] \mathrm{d} s +  \Pi_{\s} \mathrm{d} \boldsymbol{W}_{s}, \label{Eq_BF2_DPDb}
\end{align}
\end{subequations}
where $\ell \geq 1$ and with $p^{(\ell)}(s)\vert_{s=t} = X$ in $H_{\c}$ and $q^{(\ell)}(s)\vert_{s=t-\tau}=0$. Here again, the first equation (Eq.~\eqref{Eq_BF2_DPDa}) is integrated backward over $[t-\tau,t]$, followed by a forward integration of the second equation (Eq.~\eqref{Eq_BF2_DPDb}) over the same interval. To simplify the notations, we will often omit  to point out the interval $[t-\tau,t]$ in the BF systems below. 

To help interpret the parameterization produced by such an iterative scheme, we restrict momentarily ourselves to the case of a nonlinearity $G$ that is quadratic (denoted by $B$)  and to noise terms  that are scaled by a parameter $\epsilon$.  This case covers the important case of the Kardar--Parisi--Zhang equation \cite{corwin2012kardar}. The inclusion of this scaling factor $\epsilon$ allows us to group the terms constituting the stochastic parameterization according to powers of $\epsilon$  providing in particular useful small-noise expansions; see Eq.~\eqref{Eq_q2_explicit} below.

Under this working framework, if we start with $(p^{(0)},q^{(0)}) \equiv (0,0)$, we get then 
\bea \label{Eq_1st_iteration}
&p^{(1)}(s) = e^{(s-t) A_{\c}}X - \epsilon \int_{s}^t   e^{(s-s') A_{\c}}  \Pi_{\c} \d \boldsymbol{W}_{s'}, \\
&q^{(1)}(s) =  \epsilon \int_{t-\tau}^s  e^{(s-s') A_{\s}}  \Pi_{\s} \d \boldsymbol{W}_{s'},
\eea 
which leads to
\beas 
p^{(2)}(s) &= e^{(s-t) A_{\c}}X  \\
& \quad -  \int_{s}^t  e^{(s-s') A_{\c}} \Pi_{\c} B \big( p^{(1)}(s') +  q^{(1)}(s') \big) \d s' \\
& \qquad - \epsilon \int_{s}^t  e^{(s-s') A_{\c}}  \Pi_{\c} \d \boldsymbol{W}_{s'}, 
\eeas
and
\bea\label{Eq_2nd_iteration}
&q^{(2)}(s) = \int_{t-\tau}^s  e^{(s-s') A_{\s}} \Pi_{\s} B \big( p^{(1)}(s') +  q^{(1)}(s') \big) \d s'  \\
&\hspace{2cm}+  \epsilon \int_{t-\tau}^s  e^{(s-s') A_{\s}}  \Pi_{\s} \d \boldsymbol{W}_{s'}.
\eea
To further make explicit $q^{(2)}(s)$ that provides the stochastic parameterization we are seeking (after 2 iterations), we introduce 
\bea\label{Eq_fg}
f_i(s) &= -\int_{s}^t e^{(s-s') \lambda_i}  \d W^i_{s'}(\omega),  \quad i = 1, \ldots, m_c,\\
g_n(s) &= \int_{t-\tau}^s  e^{(s-s') \lambda_n} \d W^n_{s'}(\omega), \quad n \ge m_c+1.
\eea
We can then rewrite $(p^{(1)}, q^{(1)})$ given in \eqref{Eq_1st_iteration} as follows: 
\bea \label{Eq_1st_iteration_v2}
&p^{(1)}(s) = \sum_{i = 1}^{m_c} \big(e^{(s-t) \lambda_i } X_i  + \epsilon f_i(s) \big) \bm{e}_i, \\
&q^{(1)}(s) =  \epsilon \sum_{n \ge m_c+1} g_n(s) \bm{e}_n.
\eea 

By introducing additionally
\bes
B_{ij}^n = \langle B(\boldsymbol{e}_{i}, \boldsymbol{e}_{j}), \boldsymbol{e}^\ast_{n} \rangle,
\ees
we get for any $n \ge m_c+1$ the following explicit expressions 
{\small
\beas
& \Pi_{n}B(p^{(1)}(s), p^{(1)}(s)) \\
&= \epsilon^2 \sum_{i,j=1}^{m_c} B^{n}_{ij} f_{i}(s) f_{j}(s) + \epsilon \sum_{i,j=1}^{m_c} B^{n}_{ij} f_{i}(s) e^{(s-t) \lambda_{j}} X_{j} \\
&\;\; + \epsilon \sum_{i,j=1}^{m_c}\! B^{n}_{ij} f_{j}(s) e^{(s-t) \lambda_{i} } X_{i} + \sum_{i,j=1}^{m_c} \! B^{n}_{ij} e^{(s-t) (\lambda_{i} + \lambda_{j})} X_{i} X_{j}, 
\eeas}
{\small 
\beas
 \Pi_{n} B(p^{(1)}(s),q^{(1)}(s)) &= \epsilon^2 \sum_{i = 1}^{m_c} \sum_{n' =m_c+1}^{\infty} B^{n}_{i, n'} f_{i}(s) g_{n'}(s) \\
& \; + \epsilon \sum_{i=1}^{m_c} \sum_{n'=m_c+1}^{\infty} B^{n}_{i n'} g_{n'}(s) e^{(s-t)\lambda_{i}} X_{i},
\eeas}
and
\bes
\Pi_{n} B(q^{(1)}(s),q^{(1)}(s)) = \epsilon^2 \sum_{n',n'' \ge m_c+1} B^{n}_{n' n''} g_{n'}(s) g_{n''}(s).
\ees
Let us denote by $\Pi_n$ the projector onto the mode $\boldsymbol{e}_n$.
Using the above identities in \eqref{Eq_2nd_iteration}, and setting $s = t$, we get for $n\geq m_c+1$,
\bea \label{Eq_q2_explicit}
\Pi_n q^{(2)}(t) & = \epsilon^2 a_n(\tau) + \epsilon \sum_{i = 1}^{m_c} b_{in}(\tau) X_i \\
& \quad + \sum_{i,j = 1}^{m_c} c_{ij}^n(\tau) B_{ij}^n X_{i} X_{j} \\
& \quad  + \epsilon \int_{t-\tau}^t  e^{(t-s') \lambda_n} \d W^n_{s'}(\omega),
\eea 
where 
\begin{widetext}
\bea \label{Eq_an}
a_n(\tau) & =  \underbrace{\int_{t-\tau}^t e^{(t-s')\lambda_n} \sum_{i,j=1}^{m_c} B^{n}_{ij} f_{i}(s') f_{j}(s') \d s'}_{(\alpha)}     +  \underbrace{\int_{t-\tau}^t e^{(t-s')\lambda_n} \bigg(\sum_{i = 1}^{m_c} \sum_{n' =m_c+1}^{\infty} (B^{n}_{in'} + B^{n}_{n' i}) f_{i}(s') g_{n'}(s') \bigg)\d s'}_{(\beta)}\\
&\hspace{3cm} + \underbrace{\int_{t-\tau}^t e^{(t-s')\lambda_n} \bigg( \sum_{n',n'' \ge m_c+1} B^{n}_{n' n''} g_{n'}(s') g_{n''}(s') \bigg) \d s'}_{(\gamma)},
\eea
\end{widetext}
{\small 
\bea\label{Eq_bn}
b_{in}(\tau) &=  \int_{t-\tau}^t  \hspace{-1ex}  e^{(t-s')(\lambda_n - \lambda_{i})}  \sum_{j=1}^{m_c} (B^{n}_{ij} + B^{n}_{ji}) f_{j}(s') \d s' \\
&  \hspace{-1ex}  + \int_{t-\tau}^t  \hspace{-1ex}  e^{(t-s')(\lambda_n - \lambda_{i})}  \hspace{-2ex} \sum_{n'=m_c+1}^{\infty}  \hspace{-1ex}  (B^{n}_{in'} + B^{n}_{n'i}) g_{n'}(s') \d s',
\eea
}
and 
\beas
c_{ij}^n(\tau) &=  \int_{t-\tau}^t e^{(t-s')(\lambda_n - (\lambda_{i} + \lambda_{j}))} \d s'\\
                 &= \begin{cases}
\frac{1 - \exp(-\delta_{ij}^n\tau)}{\delta_{ij}^n}, & \text{if $\delta_{ij}^n\neq 0$}, \\
\tau, & \text{otherwise},
\end{cases}
\eeas
with $\delta_{ij}^n = \lambda_{i} + \lambda_{j} - \lambda_{n}$. 

In the classical approximation theory of SIM, one is interested in conditions ensuring convergence of the integrals involved in the random coefficients $a_n$ and  $b_{in}$ as $\tau\rightarrow \infty$, since $q^{(2)}(t)$ (and its higher-order analogues $q^{(\ell)}(t)$ with $\ell > 2$) aims to approximate the SIM defined in Eq.~\eqref{Eq_SIM}.
In our approach to stochastic parameterization, we do not restrict ourselves to such a limiting case but rather seek optimal backward integration time $\tau$ that minimizes a parameterization defect as explained below in Section \ref{Sec_Optim}.

However, computing this parameterization defect involves the computation of the random coefficients in the course of time (see Eq.~\eqref{Eq_minQnHn} below). The challenge is that the structure of  $a_n$ and $b_{in}$ involves repeated stochastic convolutions in time, and as such one wants to avoid a direct computation by quadrature. 
We propose below an alternative and efficient way to compute such random coefficients, for a simpler, more brutal approximation than $q^{(2)}$.  As shown in Sections \ref{Sec_ACE} and \ref{Sec_Jump}, this other class of stochastic parameterizations turns out to be highly performant for the important case where the stochastically forced scales are exclusively part of the neglected scales. 

This approximation consists of setting $q^{(1)}=0$ in Eq.~\eqref{Eq_2nd_iteration}, namely to deal with the stochastic parameterization
 \bea \label{Eq_us1_for_LIA}
\Phi(t,X,\omega) & =  \int_{t-\tau}^t  e^{(t-s) A_{\s}} \Pi_{\s} G \big( e^{(s-t) A_{\c}}X\big) \d s \\
& \qquad +  \epsilon \int_{t-\tau}^t  e^{(t-s) A_{\s}}  \Pi_{\s} \d \boldsymbol{W}_{s}(\omega),
\eea
after restoring $G$ as a more general nonlinearity than $B$. 

Note also that this parameterization $\Phi(t,X,\omega)$ is exactly the solution $q(s,X,\omega)$ at $s=t$ of the following BF system
\begin{subequations} \label{Eq_BF_SPDE_highmode_forced}
\begin{align}
& \!\d p =  A_\c p \d s, \hspace{3.9cm} s \in I_{t,\tau}, \\
& \!\d q =  \left( A_\s q + \Pi_{\s} G(p)\right) \d s + \epsilon \Pi_\s \d \bm{W}_s,  \;\;  s \in  I_{t,\tau},\\
& \!\mbox{with } p(s)\vert_{s=t} = X \in H_{\c}, \mbox{ and } q(s)\vert_{s=t-\tau}=0,
\end{align}
\end{subequations}
with $ I_{t,\tau}=[t-\tau,t]$.

Section \ref{Sec_Zn_simul}  details an efficient computation of the single stochastic convolution in equation \eqref{Eq_us1_for_LIA} using auxiliary ODEs with random coefficients. It is worth noting that previous SPDE reduction approaches have, in certain circumstances or under specific assumptions (e.g., $\Pi_c B(q,q)=0$ as seen in \cite[Eq.~(2.2)]{BHP07} and \cite[Eq.~(2.10)]{majda2001mathematical}), deliberately opted to avoid directly computing these convolutions via quadrature \cite{roberts2003step,roberts2006resolving,wang2013macroscopic}.

The stochastic convolution in Eq.~\eqref{Eq_us1_for_LIA} accounts for finite-range memory effects stemming from the stochastic forcing to the original equation. Such terms have been encountered  in the approximation of low-mode amplitudes, albeit in their asymptotic format when $\tau \rightarrow \infty$ see e.g.~\cite{BHP07,BM13}. As shown in our examples below (Sections \ref{Sec_ACE} and \ref{Sec_Jump}), stochastic convolution terms become crucial for efficient reduction when the cutoff scale, defining the retained large-scale dynamics, is larger than the forcing scale. 
There, we show furthermore that the finite-range memory content (measured by $\tau$) is a key factor for a  skillful reduction when optimized properly.
Also, as exemplified in applications, optimizing the nonlinear terms in Eq.~\eqref{Eq_us1_for_LIA} may turn to be of utmost importance to reproduce accurately the average motion of the SPDE dynamics; see Section \ref{Sec_weak_timescale}  below.

The efficient computation of repeated convolutions involved in the coefficients $a_n$ and $b_{in}$ is  however more challenging and will be communicated elsewhere.   This is not only a technical aspect though, as these repeated convolutions in $a_n$ given by Eq.~\eqref{Eq_an} characterize important triad interactions reflecting how the noise interacts through the nonlinear terms into three groups: the low-low interactions in ($\alpha$), the low-high interactions in ($\beta$) and high-high interactions in ($\gamma$).  
By using the simpler parameterization $\Phi$ defined in Eq.~\eqref{Eq_us1_for_LIA} we do not keep these interactions at the parameterization level, however as $\Phi$  approximates the high-mode amplitude it still allows us to account for triad interactions into the corresponding reduced models; see e.g.~Eq.~\eqref{Eq_closure} below.


\section{Non-Markovian Parameterizations: Formulas and Optimization}\label{Sec_OPMwhite}

The previous discussion leads us to consider, for each $n\geq m_c+1$,  the following scale-aware BF systems
\begin{subequations} \label{Eq_BF_SPDE}
\begin{align}
 & \hspace{-0.2cm}\mathrm{d} p =  A_\c  p \d s, \hspace{3.7cm} s \in [t-\tau, t],    \label{BF1_SPDE} \\
&  \hspace{-0.2cm}\mathrm{d} q  = \Big( \lambda_n q +  \Pi_{n} G\big(p \big)  \Big) \d s + \sigma_n \d W_{s}^n, s \in [t-\tau, t], \label{BF2_SPDE}\\
 & \hspace{-0.2cm} \mbox{with } p(s)\vert_{s=t} = X \in H_{\c},  \mbox{ and } q(s)\vert_{s=t-\tau}=Y  \in \mathbb{R}, \label{BF3_SPDE}
\end{align}
\end{subequations}
where $\Pi_n$ denotes the projector onto the mode $\boldsymbol{e}_n$.
Here, $\tau$ is a free parameter to be adjusted, no longer condemned to approach $\infty$ as in the approximation theory of SIM discussed above. Note that compared to Eq.~\eqref{Eq_BF_SPDE_highmode_forced}, we break down the forcing mode by mode for each high mode  to allow for adjusting the free backward parameter, $\tau$,  per scale to parameterize (scale-awareness).     
This strategy allows for a greater degree of freedom to calibrate useful parameterizations, as will be apparent in applications dealt with in Sections \ref{Sec_ACE} and \ref{Sec_Jump}.

Also,  compared to Eq.~\eqref{Eq_BF_SPDE_highmode_forced}, the initial condition for the forward integration in  Eq.~\eqref{Eq_BF_SPDE} is a scale-aware parameter  $Y$ in Eq.~\eqref{BF3_SPDE}. 
It is aimed at resolving the time-mean 
of the $n$-th high-mode amplitude, $\langle u(t),\boldsymbol{e}_n^\ast\rangle$ with $u$ solving Eq.~\eqref{Eq_SPDE}. In many applications, it is enough to set $Y=0$ though. 

In the following, we make explicit the stochastic parameterization obtained by integration of the BF system \eqref{Eq_BF_SPDE} for SPDEs with cubic nonlinear terms,  for which the stochastically forced scales are part of the neglected scales.

\subsection{Stochastic parameterizations for systems with cubic interactions} \label{Sec_LIA_cubic_case}
We consider
\be \label{Eq_SPDE_cubic}
\d u =(A u + G_2(u,u) +G_3(u,u,u)) \d t + \d \W.
\ee
Here, $A$  is a partial differential operator as defined in \eqref{Eq_A} for a suitable ambient Hilbert space $H$; see also \cite[Chap.~2]{CLW15_vol2}.
The operators $G_2$ and $G_3$ are quadratic and cubic, respectively. Here again, we assume the eigenmodes of $A$ and its adjoint $A^*$ to form a bi-orthonormal basis of $H$; namely, the set of eigenmodes $\{\boldsymbol{e}_j \}$ of $A$ and $\{\boldsymbol{e}^*_j \}$ of $A^*$ each forms a Hilbert basis of $H$, and they satisfy the bi-orthogonality condition $\langle \boldsymbol{e}_j, \boldsymbol{e}^*_k \rangle_{H} = \delta_{j,k}$, where $\delta_{j,k}$ denotes the Kronecker delta function. The noise term $\W$ takes the form defined in Eq.~\eqref{Eq_W}. 
Stochastic equations such as Eq.~\eqref{Eq_SPDE_cubic} arise in many areas of physics such as  in the description of phase separation \cite{cahn1958free,da1996stochastic,bray2002theory}, pattern formations \cite{cross1993pattern,sagues2007spatiotemporal} or transitions in geophysical fluid models \cite{sura2002sensitivity,pierini2010coherence,SD13,galfi2021applications,Margazoglou2021,simonnet2021multistability,Chekroun_al2023}.

We assume that Eq.~\eqref{Eq_SPDE_cubic} is well-posed in the sense of possessing for any $u_0$ in $H$ a unique mild solution,  i.e.~there exists a stopping time $T>0$ such that the integral equation,  
\bea
&u(t)=e^{A t}u_0\hspace{-.45ex}+\hspace{-.25ex}\int_{0}^t e^{A (t-s)} \Big(\hspace{-.2ex}G_2(u(s)) +G_3(u(s)) \Big)\hspace{-.5ex}\d s \\
& \hspace{2.5cm}+ \int_{0}^t e^{A (t-s)} \d \bm{W}_s, \;\; t\in[0,T],
\eea
possesses, almost surely, a unique solution in $C([0,T],H)$ (continuous path in $H$) for all stopping time $t \leq T$.

Here the stochastic integral can be represented as a series of one-dimensional It\^o integrals
\be \label{Eq_proj_stochconv}
 \int_{0}^t e^{A (t-s)} \d \bm{W}_s=\sum_k  \sigma_k \int_{0}^t e^{\lambda_k(t-s)} \d W_s^k \bm{e}_k,
\ee
where the $W_s^k$ are the independent Brownian motions in  \eqref{Eq_W}.  
It is known that Eq.~\eqref{Eq_SPDE_cubic} admits a unique mild solution in $H$ with the right conditions on the coefficients of $G_2$ and $G_3$ (confining  potential)  \cite{DPZ08};  see also \cite[Prop.~3.4]{hairer2011theory}.
We refer to \cite{DZ96,DPZ08} for background on the stochastic convolution term $\int_0^t e^{(t-s)A} \d \boldsymbol{W}_s$.

Throughout this section we assume that  only a subset of  unresolved modes  in  the decomposition \eqref{H_decomposition} of $H$ are stochastically forced according to  Eq.~\eqref{Eq_W}, namely that the following condition holds:
\bi
\item[({\bf H})] $\sigma_j=0$ for $1\leq j\leq m_c$ in \eqref{Eq_W}, and that there exists at least one index $k$ in $(m_c+1,N)$ such that $\sigma_k\neq 0$.
\ei
In this case, the BF system \eqref{Eq_BF_SPDE} becomes
\bea\label{Eq_BF_SPDEcubic}
& \hspace{-2ex}\mathrm{d} p =  A_\c p \d s,  \qquad \qquad\qquad \qquad \; \; \; s \in [t-\tau, t], \\
& \hspace{-2ex}\mathrm{d} q_n  = \Big( \lambda_n q_n  +  \Pi_{n} G_2\big(p\big) +  \Pi_{n} G_3\big(p \big) \Big) \d s\\
& \hspace{2.6cm} + \sigma_n \d W_{s}^n, \; \;\; \; s \in [t-\tau, t],\\
&\hspace{-2ex}  \mbox{with } p(s)\vert_{s=t} = X \in H_{\c}, \mbox{ and } q_{n}(s)\vert_{s=t-\tau}=Y. 
\eea

The solution $q_{n}$ to Eq.~\eqref{Eq_BF_SPDEcubic} at $s = t$ provides a parameterization
$\Phi_{n}(\tau,X,t)$ aimed at approximating the $n$th high-mode amplitude, $u_n(t)=\Pi_n u(t)=\langle u(t),\boldsymbol{e}_n\rangle$, of the SPDE solution $u$ when $X$ is equal to the low-mode amplitude $u_\c(t)=\Pi_\c u(t)$.  The BF system \eqref{Eq_BF_SPDEcubic} can be solved analytically. Its solution provides the stochastic parameterization, $\Phi_n$, expressed, for $n\geq m_c+1$, as the following  stochastic nonlinear mapping of $X$  in $H_\c$:
\bea\label{eq:h1tau_stoch}
\Phi_{n} &(\tau, X, t; \omega) \\
 &=  q_n(t,X, \omega) \\
 & = e^{\lambda_n \tau} Y +\sigma_n \big(W^n_{t}(\omega) -  e^{\tau \lambda_n} W^n_{t-\tau}(\omega) \big)   \\
 & \;\;  + Z^{n}_{\tau}(t; \omega)   + \sum_{i,j = 1}^{m_c} \Bigl(D^{n}_{i j}(\tau) B_{i j}^n X_{i} X_{j} \Bigr) \\
 & \;\; +\sum_{i,j,k = 1}^{m_c} \Big( E_{ijk}^n(\tau)  C_{ijk}^n  X_{i} X_{j}X_{k} \Big),
\eea
with $X=\sum_{\ell=1}^{m_c} X_{\ell} \boldsymbol{e}_{\ell}.$

The stochastic term $Z^{n}_{\tau}$ is given by
\be  \label{Eq_Zn_decomp} 
Z^{n}_{\tau}(t; \omega) =  \sigma_n \lambda_n e^{\lambda_n t} \int_{t-\tau}^t e^{-\lambda_n s} W_s^{n}(\omega) \d s,
\ee
and as such, is dependent on the noise's path and the ``past" of the Wiener process forcing the $n$-th scale: it conveys exogenous memory effects, i.e.~it is non-Markovian in the sense of \cite{HO07}.

The coefficients $B_{ij}^n$ and $C_{ijk}^n$ in \eqref{eq:h1tau_stoch} are given respectively by
\bea
B_{ij}^n = \langle G_2(\boldsymbol{e}_{i}, \boldsymbol{e}_{j}), \boldsymbol{e}^\ast_{n} \rangle, \; C_{ijk}^n = \langle G_3(\boldsymbol{e}_{i}, \boldsymbol{e}_{j}, \boldsymbol{e}_{k}), \boldsymbol{e}_n^\ast \rangle,
\eea
while the coefficients $D^{n}_{i j}(\tau)$ and $E_{ijk}^n(\tau)$ are given by 
\bea \label{Eq_Dhat}
D_{ij}^n(\tau)=  \begin{cases}
\frac{1 - \exp\big(-\delta_{ij}^n\tau\big)}{\delta_{ij}^n}, & \text{if $\delta_{ij}^n\neq 0$}, \\
\tau, & \text{otherwise},
\end{cases}
\eea
and
\bea \label{Eq_Ehat}
E_{ijk}^n(\tau)=  \begin{cases}
\frac{1 - \exp\big(-\delta_{ijk}^n\tau\big)}{\delta_{ijk}^n}, & \text{if $\delta_{ijk}^n\neq 0$}, \\
\tau, & \text{otherwise},
\end{cases}
\eea
with $\delta_{ij}^n = \lambda_{i} + \lambda_{j} - \lambda_{n}$ and $\delta_{ijk}^n = \lambda_{i} + \lambda_{j} + \lambda_{k}- \lambda_{n}$.  From Eqns.~\eqref{Eq_Dhat} and \eqref{Eq_Ehat}, one observes that the parameter $\tau$ (depending on $n$) allows, in principle, for balancing the small denominators due to small spectral gaps, i.e.~when the $\delta_{ij}^n$ or the $\delta_{ijk}^n$ are small. This attribute is shared with the parameterization formulas obtained by the BF approach in the deterministic context; see \cite[Remark III.1]{chekroun2023optimal}.
It plays an important role in the ability of our parameterizations to handle physically relevant situations with a weak time-scale separation; see Section \ref{Sec_weak_timescale} below.

\begin{rem} \label{Rmk_OU}
In connection with the stochastic convolution $ \xi_{t,\tau}=\int_{t-\tau}^t  e^{(t-s) A_{\s}}  \Pi_{\s} \d \boldsymbol{W}_{s}$ in Eq.~\eqref{Eq_us1_for_LIA}, we note that the stochastic terms   in Eq.~\eqref{eq:h1tau_stoch}, involving the processes $W^n_{t}$, $W^n_{t-\tau}$, and $Z^{n}_{\tau}$, 
 can be rewritten as a scalar stochastic integral, more precisely:
 \be\label{Eq_OU_corresp}
 \int_{t-\tau}^t  e^{\lambda_n(t-s)}  \d W^n_{t}= W^n_{t}-  e^{\tau \lambda_n} W^n_{t-\tau} + \sigma_n^{-1}Z^{n}_{\tau}(t).
 \ee
 Note that the RHS in Eq.~\eqref{Eq_OU_corresp} is 
 simply the projection of $\xi_{t,\tau}$ onto mode $\bm{e}_n$ (Eq.~\eqref{Eq_proj_stochconv}). 
 
Equation \eqref{Eq_OU_corresp}  is a direct consequence  of It\^o's formula applied to the product $e^{-\lambda s} W_s^n$ followed by integration over $[t-\tau, t]$. 
Indeed, we have  (dropping the $n$-dependence),
\beas\label{Eq_toto}
\d \Big(e^{-\lambda s} W_s\Big)&= -\lambda e^{-\lambda s} W_s \d s  +e^{-\lambda s} \d W_s + \d  (e^{-\lambda s}) \d W_s\\
&=-\lambda e^{-\lambda s}W_s  \d s  + e^{-\lambda s} \d W_s,
\eeas
by treating  $\d s \d W_s$ as zero due to It\^o calculus \cite{oksendal2013stochastic}.

By introducing now $ \delta W= W_t -e^{\lambda \tau} W_{t-\tau}$ and denoting by  $J$ the stochastic convolution  $\int_{t-\tau}^t  e^{\lambda (t-s)}  \d W_{s}$, we arrive due to Eq.~\eqref{Eq_Zn_decomp} at  
\be
e^{-\lambda t} \delta W+\lambda e^{-\lambda t}Z_\tau=e^{-\lambda t} J,
\ee
which gives the relation $J=\delta W+\lambda Z_\tau$, where we have set  $\sigma_n=1$ to simplify. 

In the expression of $\Phi_n$ given by Eq.~\eqref{eq:h1tau_stoch}, we opt for the expression involving the Riemann–Stieltjes  integral $Z^{n}_{\tau}$ for practical purpose.  The latter can be indeed readily generalized to the case of jump noise (Eqns.~\eqref{eq:h1tau_jumpnoise}--\eqref{Eq_J_inital} below) without relying on stochastic calculus involving a jump measure (Section \ref{Sec_Levy}).
\end{rem}

 In what follows we drop the dependence of $\Phi_{n}$ on $\omega$ and sometimes on $t$, unless specified.
The parameterization of the neglected scales writes then as
\be\label{Eq_LIA}
\Phi_{\bm{\tau}}(X,t)=\sum_{n\geq m_c+1}\Phi_n(\tau_n,X,t) \bm{e}_n,\\ 
\ee 
with $\bftau= (\tau_n)_{n\geq m_c+1}$, and $\Phi_n$ defined in Eq.~\eqref{eq:h1tau_stoch}.

We observe that letting $\tau$ approach infinity in Eq.~\eqref{Eq_LIA} recovers  Eq.~\eqref{Eq_us1_for_LIA} when $G = G_2 + G_3$. This implies that the manifold constructed from Eq.~\eqref{Eq_LIA} can be viewed as a homotopic deformation of the SIM approximation defined in Eq.~\eqref{Eq_us1_for_LIA}, controlled by the parameters $\tau_n$.

This connection motivates utilizing the variational framework presented in  \cite{CLM20_closure,chekroun2023optimal} to identify the ``optimal" stochastic manifold within the family defined by Eq.~\eqref{Eq_LIA}.  In Section \ref{Sec_Optim}, we demonstrate that a simple data-driven minimization of a least-squares parameterization error, applied to a single training path, yields trained parameterizations (in terms of $\tau$) with remarkable predictive power. For instance, these parameterizations can be used to infer ensemble statistics of SPDEs for a large set of unseen paths during training.

However, unlike the deterministic case, the stochastic framework necessitates addressing the efficient simulation of the coefficients involved in Eq.~\eqref{Eq_LIA} for optimization purposes. This challenge is tackled in the following subsection.

\subsection{Non-Markovian path-dependent coefficients: Efficient simulation}  \label{Sec_Zn_simul}
Our stochastic parameterization given by \eqref{Eq_LIA} contains random coefficients $Z^{n}_{\tau}(t; \omega)$, each of them involving an integral of the history of the Brownian motion making them non-Markovian in the sense of \cite{HO07}, i.e.~depending on the noise path history. We present below an efficient mean to compute these random coefficients by solving auxiliary ODEs with path-dependent coefficients, avoiding this way the cumbersome computation of integrals over noise paths that would need to be updated at each time $t$. This approach is used for the following purpose: 
\begin{itemize}
\item[(i)] To find the optimal vector $\bftau^\ast$ made of optimal backward-forward integration times and form in turn  the  stochastic optimal parameterization $\Phi_{\bftau^\ast}$, and
\item[(ii)] To efficiently simulate the corresponding optimal reduced system built from $\Phi_{\bftau^\ast}$ along with its path-dependent coefficients.
\end{itemize}
According to \eqref{Eq_Zn_decomp}, the computation of $Z^{n}_{\tau}$  
 boils down to the computation of the random integral 
 \bes
 e^{\lambda_n t}  \int_{-\tau + t}^t e^{-\lambda_n s} W_{s}^{n}(\omega) \d s,
 \ees 
 which is of the following form: 
\be \label{Eq_def_I}
I_{\tau}^{j}(t, \omega; \kappa) = e^{\kappa t} \int_{-\tau + t}^t e^{-\kappa s} W_s^{j}(\omega) \d s,
\ee
where  $\tau > 0$, and $j\geq m_c+1$ due to assumption ({\bf H}). We only need to consider the case that $\kappa < 0$, since we  consider here only cases in which the unstable modes are included within the space $H_\c$ of resolved scales.

By taking time derivative on both sides of \eqref{Eq_def_I}, we obtain that $I_{\tau}^j(t, \omega; \kappa)$ satisfies the following ODE with path-dependent coefficients also called a random differential equation (RDE):
\be \label{Eq_for_I}  
\frac{\d I}{\d t} = \kappa I + W_t^{j}(\omega) - e^{\kappa \tau} W_{t-\tau}^{j}(\omega).
\ee
Since $\kappa<0$, the linear part $\kappa I$ in \eqref{Eq_for_I} brings a stable contribution to \eqref{Eq_for_I}. As a result, $I_{\tau}^j(t, \omega; \kappa)$ can be  computed by integrating  \eqref{Eq_for_I} forward in time starting at $t=0$ with initial datum given by: 
\be \label{IC_for_I_case1}
I(0,\omega) = I_{\tau}^j(0, \omega; \kappa),
\ee
where $I_{\tau}^j(0, \omega; \kappa)$ is computed using \eqref{Eq_def_I}. 

Now, to determine $Z^{n}_{\tau}(t; \omega)$ in \eqref{Eq_Zn_decomp} it is sufficient to observe that 
\be \label{Eq_Zn_decomp2}
Z^{n}_{\tau}(t; \omega) =   \sigma_n \lambda_n  I_{\tau}^n(t, \omega; \lambda_n).
\ee
The random coefficient $Z^{n}_{\tau}(t; \omega)$ is thus computed by using \eqref{Eq_Zn_decomp2} after $I_{\tau}^n(t, \omega; \lambda_n)$  is computed by solving, forward in time, the corresponding RDE of the form \eqref{Eq_for_I} with initial datum  
\be\label{Eq_init}
I(0,\omega) = I_{\tau}^n(0, \omega; \lambda_n) =  \int_{-\tau}^0 e^{-\lambda_n s} W_s^{n}(\omega) \d s. 
\ee
This way, we compute only the integral \eqref{Eq_init} and then propagate it through  Eq.~\eqref{Eq_for_I} to evaluate $Z^{n}_{\tau}(t; \omega)$, instead of computing for each $t$ the integral in the definition  \eqref{Eq_Zn_decomp}. The resulting procedure is thus much more efficient to handle numerically than a direct computation of $Z^{n}_{\tau}(t; \omega)$ based on the integral formulation \eqref{Eq_Zn_decomp} which would involve a careful and computationally more expensive quadrature at each time-step. 
Note that a similar treatment has been adopted in \cite[Chap.~5.3]{CLW15_vol2} for  the computation of the time-dependent coefficients arising in stochastic parameterizations of SPDEs driven by linear multiplicative noise, and in \cite[Sec.~III.C.2]{chekroun2023optimal} for the case of non-autonomous forcing.

\subsection{Data-informed optimization and training path}\label{Sec_Optim}
 Thanks to Section \ref{Sec_Zn_simul}, we are now in position to 
propose an efficient variational approach to optimize the stochastic parameterizations of Section \ref{Sec_LIA_cubic_case} in view of handling parameter regimes located away from the instability onset and thus a greater wealth of possible stochastic transitions.

\begin{figure} 
\centering
\includegraphics[width=1\linewidth, height=0.4\linewidth]{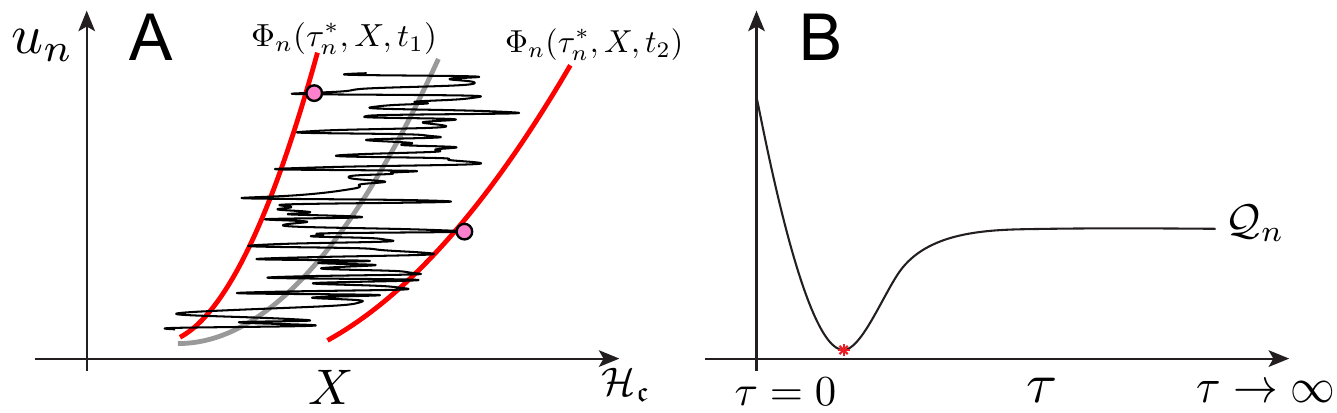} 
\caption{{\bf Panel A}: This panel illustrates a training solution path (black curve) for the stochastic optimal parameterization method. This path is  shown to be transverse to the underlying stochastic manifolds  (red curves) as occurring typically for the reduction of stochastic systems using such objects (see \cite[Fig.~2]{Chekroun_al2023}) for a concrete simple example).  
The optimal parameterization (denoted by $\Phi_n(\tau_n^\ast,X,t)$) aims to approximate the actual state (pink dots) of the system's solution (represented by the black curve) at any given time (here labeled as $t_1$ and $t_2$).  A key strength of this parameterization is its ability to capture the solution's average behavior by its expected value, $\mathbb{E}(\Phi_n)$ (gray curve).
{\bf Panel B}: This panel schematically shows how the ``goodness-of-fit" (parameterization defect, $\mathcal{Q}_n$, defined by Eq.~\eqref{Eq_minQnHn}) changes with the parameter $\tau$. The red asterisk marks the value $\tau_n^\ast$ that minimizes this defect.}
\label{Fig_intro}
\end{figure}

Given a solution path $u(t)$ to Eq.~\eqref{Eq_SPDE_cubic} available over an interval $I_T$ of length $T$, we aim for determining the optimal parameterization $\Phi_n$ (given by \eqref{eq:h1tau_stoch}) that minimizes---in the $\tau$-variable---the following {\it parameterization defect}
\be\label{Eq_minQnHn}
\mathcal{Q}_n(\tau,T)= \overline{\big|u_n(t) -\Phi_{n}(\tau,u_\c (t), t) \big|^2},
\ee
for each $ n \geq m_c+1 $. Here $\overline{(\cdot )}$ denotes the time-mean over $I_T$ while $u_n(t)$ and $u_\c(t)$ denote the projections of $u(t)$ onto the high-mode $\boldsymbol{e}_n$ and the reduced state space $H_\c$, respectively. 
Figure \ref{Fig_intro} provides a schematic of the stochastic optimal parameterizing manifold (OPM) found this way. 

 The objective is obviously not to optimize the parameterizations $\Phi_n$ for every solution path, but rather to optimize $\Phi_n$ on a single solution path---called the {\it training path} and then use the optimized parameterization for predicting dynamical behaviors for any other noise path, or at least in some statistical sense. To do so, the optimized $\tau_n$-values (denoted by $\tau_n^\ast$ below) is used to build the resulting stochastic  OPM  given by \eqref{Eq_LIA} whose optimized  non-Markovian  coefficients are aimed at  encoding efficiently the interactions between the nonlinear terms and the noise, in the course of time.

Sections~\ref{Sec_ACE} and \ref{Sec_Jump} below illustrate how this single-path training strategy can be used efficiently to  predict from the corresponding reduced systems  the statistical behavior as time and/or the noise's path is varied. 
The next section delves into the justification of this variational approach, and the theoretical characterization of the notion of stochastic OPM by relying on ergodic arguments. The reader interested in applications can jump to Sections \ref{Sec_ACE} and \ref{Sec_Jump}.

\section{Non-Markovian Optimal Parameterizations and Invariant Measures}\label{Sec_OPM_IM_main}

 The problem of ergodicity and mixing of dissipative properties of PDEs subject to a stochastic external force has been an active research topic over the last two decades. It is rather 
well understood in the case when all deterministic modes are forced; see e.g.~\cite{goldys2005exponential,kuksin2012mathematics}. The situation in which only finitely many modes are forced such as considered in this study is much more subtle to handle in order to prove unique ergodicity.  The works  \cite{hairer2006ergodicity,hairer2011theory} provide answers in such a degenerate situation based on a theory of hypoellipticity for nonlinear SPDEs. In parallel, the work \cite{hairer2011asymptotic} generalizing the asymptotic coupling method introduced in \cite{hairer2002exponential}, allows for covering a broad class of nonlinear stochastic evolution equations including stochastic delay equations.

Building upon these results, we proceed under the assumption of a unique ergodic invariant measure, denoted by $\mu$. Given a reduced state space $H_\c$, we  demonstrate in this section the following:
\bi
    \item[(i)] The path-dependence of the non-Markovian optimal reduced model arises from the random measure denoted by $\rho_\omega$ that is obtained through the disintegration of $\mu$ over the underlying probability space $\Omega$.
     \item[(ii)] The non-Markovian optimal reduced model provides an averaged system in $H_\c$ that is still path-dependent. For a given noise path $\omega$, it provides the reduced system in $H_\c$ that averages out the unresolved variables with respect to the disintegration of $\rho_\omega$ over $H_\c$. A detailed explanation is provided in Theorem \ref{Thm_OPMclosure2}.
\ei

\subsection{Theoretical insights}\label{Sec_OPM_IM}
Adopting the framework of random dynamical systems (RDSs) \cite{Arnold98}, recall that an RDS, $S(t,\omega)$, is said to be {\it white noise}, if the associated ``past" and ``future" $\sigma$-algebras $\mathcal{F}^-$ (see \eqref{F_}) and $\mathcal{F}^+$ are independent.
Given a white noise RDS, $S(t,\omega)$, the relations between random invariant measures for $S$ and invariant
measures for the associated Markov semigroup $P_t$ is well known \cite{crauel2008measure}. We recall below these relationships and enrich the discussion with Lemma \ref{Lemma_random_from_mu} below from \cite{flandoli2017synchronization}.

Assume that Eq.~\eqref{Eq_SPDE} generates a white noise RDS, $S(t,\omega)$, associated with a Markov semigroup $P_t$ having $\mu$ as an invariant measure, then the following limit taken in the weak$^\ast$ topology:
\be\label{def1}
\rho_{\omega}=\lim_{t\rightarrow \infty} S(t,\theta_{-t}\omega)  \mu,
\ee
exists $\mathbb{P}$-a.s., and provides a random probability measure that is $S$-invariant in the sense that 
\be\label{Eq_Sinv}
S(t,\omega)  \rho_{\omega}=\rho_{\theta_t \omega}, \;\;\; \mathbb{P}\mbox{-a.s.}, \;\;  t\geq 0.
\ee
A random probability measure $\rho_\omega$ that satisfies \eqref{Eq_Sinv} is called below a (random) statistical equilibrium.  
Recall that $\theta_t$ here denotes the standard ergodic transformation on the set of Wiener path $\Omega$ defined through  the helix identity  \cite[Def.~2.3.6]{Arnold98}, $W_s(\theta_t\omega)=W_{t+s}(\omega)-W_t(\omega)$, $\omega$ in $\Omega$.

A statistical equilibrium  such as $\rho_\omega$ defined above is furthermore Markovian \cite{crauel1991markov}, in the sense it is measurable with respect to the past $\sigma$-algebra
\be\label{F_}
\mathcal{F}^-=\sigma\{\omega\mapsto S(\tau,\theta_{-t}\omega)\;:\; 0\leq \tau \leq t\}; 
\ee
see \cite[Prop.~4.2]{crauel2008measure}; see also \cite[Theorems 1.7.2 and 2.3.45]{Arnold98}.

Note that the limit \eqref{def1} exists in the sense that for every bounded measurable function $f :H \rightarrow \mathbb{R}$, the real-valued stochastic process 
\be
(t,\omega)\mapsto \int_{H} f(x) \d  \big(S(t,\theta_{-t}\omega)\mu\big)=\int_{H} f( S(t,\theta_{-t}\omega)x) \d \mu,
\ee
is a bounded martingale, and therefore converges $\mathbb{P}$-a.s by the Doob's first martingale convergence theorem; see e.g.~\cite[Thm.~C.5 p.~302]{oksendal2013stochastic}. This implies, in
particular, the $\mathbb{P}$-almost sure convergence of the measure-valued random variable $\omega\mapsto S(t,\theta_{-t}\omega)\mu$ in the topology of weak convergence of $Pr_{\Omega}(H)$, which in turn implies the convergence  in the narrow topology of $Pr_{\Omega}(H)$; see \cite[Chap.~3]{crauel2003random}.

Reciprocally, given a Markovian $S$-invariant random measure $\rho_{\omega}$ of an RDS for which the $\sigma$-algebras $\mathcal{F}^-$ and $\mathcal{F}^+$ are independent, the probability measure, $\mu$, defined by 
\be\label{cond1}
\mu=\mathbb{E} (\rho_{\bullet})=\int_\Omega \rho_{\omega} \d \mathbb{P}(\omega),
\ee
is an invariant measure of the Markov semigroup $P_t$; see \cite[Thm.~1.10.1]{chueshov2002monotone}.

In case of uniqueness of $\mu$ (and thus ergodicity) satisfying $P_t \mu=\mu,$ then this one-to-one correspondence property implies that any other random probability measure $\nu_{\omega}$ different from the random measure $\rho_{\omega}$ obtained from  \eqref{def1}, is non-Markovian which means in particular that $\int_\Omega \nu_{\omega} \d \mathbb{P}\neq \mu$ otherwise we would have $\nu_\omega=\rho_\omega$, $\mathbb{P}$-a.s.

In case a random attractor $\mathcal{A}(\omega)$ exists, since the latter is a forward invariant compact random set, the Markovian random measure $\rho_{\omega}$ must be supported by $\mathcal{A}(\omega)$; see \cite[Thm.~6.17]{crauel2003random} (see also \cite[Prop.~4.5]{CF94}). 
We conclude thus \textemdash\, in the case of a unique ergodic measure for $P_t$ \textemdash\, to the existence of {\it a unique} Markov  $S$-invariant random measure $\rho_\omega$ supported by $\mathcal{A}$, and that all other $S$-invariant random measures (also necessarily supported by $\mathcal{A}$) are non-Markovian. 

RDSs where the future and past $\sigma$-algebras are independent, are generated for many stochastic systems in practice. This is the case for instance of a broad class of stochastic differential equations (SDEs); see \cite[Sect.~2.3]{Arnold98}. The problem of generation of white noise RDSs in infinite dimension is much less clarified and it is beyond the scope of this article to address this question. 

Instead, we point out below that there exists other ways to associate to an RDS which is not necessarily of white-noise type,  a meaningful random probability measure $\rho_\omega$ that still for instance satisfies \eqref{cond1}. 
This is indeed the case if the RDS satisfies some weak mixing property with respect to a (non-random) probability $\mu$ on $H$ as the following Lemma from \cite{flandoli2017synchronization}, shows:
\bl\label{Lemma_random_from_mu}
Assume there exists $\mu$ in Pr($H$) satisfying the following weak mixing property
\be\label{weak_mix}
\frac{1}{t}\int_0^t \d s \int_{H}\mathbb{E} f(S(s,\theta_{-s}\bullet)u) \d \mu    \underset{t\rightarrow \infty}\longrightarrow \int_{H} f \d \mu,
\ee
for each bounded, continuous $f: H\rightarrow \mathbb{R}$.  

Then there exists an $\mathcal{F}_0$-measurable random probability measure $\rho_\omega$
given by
\be\label{Eq_rho_def}
\rho_{\theta_t\omega}=\frac{1}{t}\int_0^t S(s,\theta_{-s}\omega)\mu \d s,
\ee
that satisfies 
\be\label{cond2}
\mathbb{E} (\rho_{\bullet})=\int_\Omega \rho_{\omega} \d \mathbb{P} (\omega)=\mu,
\ee
and is $S$-invariant in the sense of \eqref{Eq_Sinv}, i.e.~that is a statistical equilibrium. 
\el
For a proof of Lemma \ref{Lemma_random_from_mu} see that of \cite[Lemma 2.5]{flandoli2017synchronization}.

Thus, either built from an invariant measure of the Markov semigroup, in the case of a white-noise RDS, or 
from a weakly mixing measure in the sense of \eqref{weak_mix}, a statistical equilibrium satisfying \eqref{cond1}  can be naturally associated with Eq.~\eqref{Eq_SPDE_cubic} as long as the appropriate assumptions are satisfied. 
At this level of generality, we do not enter into addressing the important question dealing with sufficient conditions on the linear part $A$ and and nonlinear terms that ensure for an SPDE like Eq.~\eqref{Eq_SPDE_cubic}, the generation of an RDS that is either of white-noise type or weakly mixing in the sense of \eqref{weak_mix}; see \cite{flandoli2017synchronization} for examples in the latter case. 

We also introduce the following notion of  pullback parameterization defect. 
\bd
Given a parameterization $\Phi: \Omega\times H_\c \rightarrow H_\s$ that is measurable (in the measure-theoretic sense), and a mild solution $u$ of Eq.~\eqref{Eq_SPDE_cubic}, the pullback parameterization defect associated with $\Phi$ over the interval $[0,T]$ and  for a noise-path $\omega$, is defined as:
\be\label{Def_Q2}
 \mathcal{Q}_T(\omega,\Phi)=\frac{1}{T}\int_{0}^T \norm{u_\s(t,\theta_{-t}\omega) 
 -\Phi (\omega,u_\c(t,\theta_{-t}\omega))}^2 \d t.
\ee
\ed

Given the  statistical equilibrium $\rho_\omega$, 
we denote by $\m_\omega$  the push-forward of $\rho_\omega$ by the projector $\Pi_\c$ onto $H^\c$, namely 
 \be
    \m_\omega (B)=\rho_{\omega} (\Pi_\c^{-1}(B)), \;\;  B \in \mathcal{B}(H_\c),
    \ee
where $\mathcal{B}(H_{\c})$ denotes the family of Borel sets of $H_\c$; i.e. the family of sets that can be formed from open sets (for the topology on $H_\c$ induced by the norm $\|\cdot\|_{H_\c}$) through the operations of countable union, countable intersection, and relative complement.

The random measure $\m_\omega$ allows us to consider the following functional space $\mathcal{X}=L^2_{\m}(  \Omega \times H_\c; H_\s)$ of parameterizations, defined as
\bea\label{E_class}
&\mathcal{X}=\bigg\{\Phi:    \Omega \times H_\c\rightarrow H_\s\, \textrm{measurable} \; \Big| \\
 & \hspace{1cm} \; \int_{H_\c} \|\Phi(\omega,X) \|^2\d \m_{\omega}(X) <\infty \; a.s.\bigg\},
\eea
namely the Hilbert space constituted by $H_\s$-valued random functions of the resolved variables $X$ in $H_\c$, that are square-integrable with respect to $\mathfrak{m}_\omega$, almost surely.

We are now in position to formulate the main result of this section, namely the following generalization to the stochastic context of  Theorem  4 from \cite{CLM20_closure}.

\bt\label{Thm_variational-pb2}
Assume that one of the following properties hold:
\bi
\item[(i)] The RDS generated by Eq.~\eqref{Eq_SPDE_cubic} is of white-noise type and the Markov semigroup associated with Eq.~\eqref{Eq_SPDE} admits a (unique) ergodic invariant measure $\mu$. 
\item[(ii)] The RDS generated by Eq.~\eqref{Eq_SPDE_cubic} is weakly mixing in the sense of \eqref{weak_mix}. 
\ei
Let us denote by $\rho_{\omega}$ the weak$^\ast$-limit of $\mu$ defined in \eqref{def1} for case (i), and by $\rho_{\omega}$
the statistical equilibrium ensured by Lemma \ref{Lemma_random_from_mu}, in case (ii).  
We denote by $\rho_{\omega}^X$ the disintegration of $\rho_\omega$ on the small-scale subspace $H_\s$, conditioned on the coarse-scale variable $X$.

Assume that the small-scale variable $Y$ has a finite energy in the sense that  
\be\label{Eq_finite-energy2}
\int\int \norm{Y}^2 \d \rho_\omega \d \mathbb{P} <\infty.
\ee

Then 
the  minimization problem
\be\label{Eq_variational-pb2}
\underset{\Phi  \in \mathcal{X}}\min \int_{\Omega} \int_{(X,Y)\in  H_\c \times H_\s} \hspace{-1cm}\norm{Y -\Phi(\omega,X)}^2 \d \rho_{\omega} (X,Y)\d \mathbb{P} (\omega),
\ee
possesses a unique solution whose argmin  is given 
by
\be\label{Def_h2rand}
 \Phi^\ast(\omega,X)=\int_{H_\s} Y \d \rho_{\omega}^{X}(Y), \;\; X \in H_\c, \;\;\; \mathbb{P}\mbox{-}\rm{a.s.}
\ee

Furthermore, if the RDS possesses a pullback attractor $\mathcal{A}(\omega)$ in $H$, and $\rho_{\omega}$ is pullback mixing in the sense that for all $u$ in $H$,
 \bea\label{strong_ergo_pullback}
& \underset{T\rightarrow \infty}\lim \frac{1}{T} \int_0^T f(\omega,S(t,\theta_{-t}\omega)u) \d t= \\
 &\hspace{.5cm}\int_{\mathcal{A}(\omega)} f(\omega,v)\d \rho_{\omega}(v),  \; \; f(\omega,\cdot)\in C_b(H), \; \mathbb{P}\mbox{-a.s.},
\eea 
then $\mathbb{P}$-almost surely
\be\label{Var_propertyQ2}
\underset{T\rightarrow \infty}\lim\mathcal{Q}_T(\omega,\Phi^\ast) \leq \underset{T\rightarrow \infty}\lim\mathcal{Q}_T(\omega,\Phi), \;\; \Phi  \in \mathcal{X}.
\ee

\et

\bp
The proof follows the same lines of \cite[Theorem 4]{CLM20_closure}.
It consists of replacing:
\bi 
\item[(i)] The probability measure $\mu$ therein by the probability measure 
$\rho$ on $\Omega\times  H$ that is naturally associated with the family of random measures $\rho_{\omega}$ and whose marginal on $\Omega$ is given by $\mathbb{P}$; see  \cite[Prop.~3.6]{crauel2008measure}, 
\item[(ii)] The function space $L^2_{\mu}( H_\c\times H_\s; H_\s)$ therein by the function space $\mathcal{E}=L^2_{\rho}(\Omega\times H_\c\times H_\s; H_\s)$. 
\ei
By applying to the ambient Hilbert space  $\mathcal{E}$, the standard projection theorem onto closed convex sets \cite[Theorem 5.2]{brezis_book}, one defines (given $\Pi_\c$) the conditional expectation $\mathbb{E}_\rho[g\vert \Pi_\c]$ of $g$ as the unique function in $\mathcal{E}$ that satisfies the inequality 
\be\label{Eq_best_app}
\mathbb{E}_\rho[\|g-\mathbb{E}_\rho[g| \Pi_\c]\|^2] \leq \mathbb{E}_\rho[\| g-\Psi \|^2], \; \mbox{for all } \Psi \in\mathcal{E}.
\ee

Now by applying the general disintegration theorem of probability measures, applied to $\rho$ (see \cite[Eq.~(3.18)]{CLM20_closure}),
we obtain the following explicit representation of the random conditional expectation 
\be
\mathbb{E}_\rho[g| \Pi_\c](\omega,X)=\int_{H_\s} g(\omega,X,Y) \d\rho_{\omega,X}(Y),
\ee
with $\rho_{\omega,X}$ denoting the disintegrated measure of $\rho$ over $\Omega\times H_\c$. By noting that this disintegrated measure is the same 
as the disintegration $\rho_{\omega}^{X} $ (over $H_\c$) of the random measure $\rho_\omega$, we conclude by taking $g(\omega,X,Y)=Y$ as for the proof of \cite[Theorem 4]{CLM20_closure} that 
$ \Phi^\ast$ given by \eqref{Def_h2rand} solves the minimization problem \eqref{Eq_variational-pb2}. 

The inequality \eqref{Var_propertyQ2} results then from the pullback mixing property \eqref{strong_ergo_pullback} and the definition \eqref{Def_Q2} of the pullback parameterization defect. 
\ep

\subsection{Non-Markovian optimal reduced model, and conditional expectation}

The mathematical framework of Section \ref{Sec_OPM_IM} allows us to provide a useful interpretation of non-Markovian optimal reduced model for SPDE of type Eq.~\eqref{Eq_SPDE_cubic}. To simplify the presentation, we restrict ourselves here to the case $G_3=0$. 

We denote then by $G(v)$ the vector field in $H$, formed by gathering the linear and  nonlinear parts of Eq.~\eqref{Eq_SPDE_cubic} in this case, namely 
\be\label{Eq_vector_field}
G(v)=A v + G_2(v,v), \quad  \; v \in H.
\ee 

The theorem formulated below characterizes the relationships between the non-Markovian optimal reduced model  and the random conditional expectation associated with the projector $\Pi_\c$ onto the reduced state space $H_\c$; i.e.~the resolved modes.  Its proof is almost  identical to that of  \cite[Theorem 5]{CLM20_closure} 
with only slight amendment and the details are thus omitted.

 \bt\label{Thm_OPMclosure2}
 Consider the SPDE of type Eq.~\eqref{Eq_SPDE_cubic} with $G_3=0$. 
 Let $\rho_\omega$  be a (random) statistical equilibrium satisfying either \eqref{def1} (in case (i) of Theorem \ref{Thm_variational-pb2}) or ensured by  Lemma \ref{Lemma_random_from_mu} (in case (ii) of Theorem \ref{Thm_variational-pb2}).  

Then, under the conditions of Theorem \ref{Thm_variational-pb2},  the random  conditional expectation associated with $\Pi_\c$,
\bes
\mathbb{E}_{\rho_{\theta_t\omega}}[G| \Pi_\c](X)=\int_{Y\in H_\s} G(X+Y) \d\rho_{\theta_t\omega}^{X}(Y), \; X\in H_\c,
\ees
 satisfies 
\bea
& \mathbb{E}_{\rho_{\theta_t\omega}}[G| \Pi_\c](X)=A_\c X+\Pi_\c G_2(X,X) \\
& \hspace{.6cm}+2 \Pi_\c G_2(X,\Phi^\ast(\theta_t \omega,X))+\langle G_2(Y,Y) \rangle_{\rho_{\theta_t\omega}^{X}},
\eea
with 
\be
\langle G_2(Y,Y) \rangle_{\rho_{\theta_t\omega}^{X}}=\int_{Y \in H_\s} \hspace{-.2cm}\Pi_\c G_2(Y,Y)\d\rho_{\theta_t\omega}^{X}(Y),
\ee
for which  $\rho_{\omega}^{X} $ denotes the disintegration of  $\rho_\omega$ over $H_\c$.

Then, given the reduced state space $H_\c$ and the statistical equilibrium $\rho_{\omega}$,  
if  $\langle G_2(Y,Y) \rangle_{\rho_{\theta_t\omega}^{X}}=0$, the {\bf non-Markovian optimal reduced model} of Eq.~\eqref{Eq_SPDE_cubic}  
with noise acting into the ``orthogonal direction'' of $H_\c$ (assumption ({\bf H}) in Section \ref{Sec_LIA_cubic_case}), is given by:
\be\label{OPM_closure2}
\d X =  \mathbb{E}_{\rho_{\theta_t\omega}}[G| \Pi_\c](X) \d t.
\ee 
\et
Hereafter, we simply refer to  Eq.~\eqref{OPM_closure2} as the (random) conditional expectation.

\subsection{Practical implications}\label{Sec_classP}
 Thus, Theorem \ref{Thm_OPMclosure2} teaches us that an approximation of the (actual) non-Markovian optimal parameterization, $\Phi^\ast$ given by \eqref{Def_h2rand}, 
provides {\it in fine} an approximation of the conditional expectation involved in Eq.~\eqref{OPM_closure2}.

The non-Markovian optimal parameterization involves, from its definition, averaging with respect to the unknown probability measure $\rho_\omega$.
As such, designing a practical approximation scheme with rigorous error estimates remains challenging. Near the instability onset, the non-Markovian optimal parameterization simplifies to the stochastic invariant manifold. In this regime, probabilistic error estimates have been established for a wide range of SPDEs, including those relevant to fluid dynamics  \cite{Chekroun_al2023}. However, obtaining such guarantees becomes significantly more difficult for scenarios away from the instability onset, even for low-dimensional SDEs \cite{chekroun2019grisanov}.

This is where the data-informed optimization approach from Section \ref{Sec_Optim} offers a practical solution within a relevant class of parameterizations.  Variational inequality \eqref{Var_propertyQ2} provides a key insight: the parameterization defect serves as a good measure of a parameterization's quality. Lower defect values indicate a better parameterization, one that is expected to yield a reduced model closer to the (theoretical) non-Markovian optimal reduced model  defined by Eq.~\eqref{OPM_closure2}.

Even with a very good approximation of the actual non-Markovian optimal parameterization, $\Phi^\ast$, a key question remains: under what conditions does the theoretical conditional expectation (Eq.~\eqref{OPM_closure2}) provide a sufficient system's closure on its own? We address this question in Sections \ref{Sec_ACE} and \ref{Sec_Jump}, through concrete examples. Analyzing two universal models for non-equilibrium phase transitions, we demonstrate that the non-Markovian reduced models constructed using the formulas from Section \ref{Sec_OPMwhite}  become particularly relevant for understanding and predicting the underlying stochastic transitions when the cutoff scale exceeds the scales forced stochastically (Sections \ref{Sec_unravel},\ref{Sec_stat_reduced},\ref{Sec_weak_timescale} and \ref{Sec_opti_jumps}).

To deal with these non-equilibrium systems, we consider the class $\mathcal{P}$ of parameterizations for approximating the true non-Markovian optimal parameterization, $\Phi^\ast$. This class $\mathcal{P}$ consists of continuous deformations away from the instability onset point of parameterizations that are valid near this onset (as discussed in Section \ref{Sec_LIA_cubic_case}).

Our numerical results demonstrate that this class  $\mathcal{P}$  effectively models and approximates the true non-Markovian terms in the optimal reduced model (Eq.~\eqref{OPM_closure2}) across various scenarios. After learning an optimal parameterization within class $\mathcal{P}$ using a single noise realization, the resulting approximation of the optimal reduced model typically becomes an ODE system with coefficients dependent on the specific noise path. These path-dependent coefficients encode the interactions between noise and nonlinearities, as exemplified by the reduced models Eqns.~\eqref{Eq_closure} and \eqref{Eq_reduced_S_shaped} below, with their ability to predict noise-induced transitions within the reduced state spaces.

\section{Predicting Stochastic Transitions with Non-Markovian Reduced Models}\label{Sec_ACE}
\subsection{Noise-induced transitions in a stochastic Allen-Cahn model}
We consider now the following stochastic Allen-Cahn Equation (sACE)
\cite{allen1979microscopic,Chafee_al74,bray2002theory}
\be  \label{eq:Chafee}
\d u = (\partial_{x}^2 u   + u-u^3) \d t +  \mathrm{d} \W,  \quad 0<x<L, \;\, t>0,
\ee
with homogeneous Dirichlet boundary conditions. The ambient Hilbert space is $H=L^2(0,L)$ endowed with its natural  inner product denoted by $\langle \cdot, \cdot \rangle$.

This equation and variants have a long history. The deterministic part of the equation provides the gradient flow of the Ginzburg-Landau free energy functional \cite{ginzburg1950theory} 
\be\label{GL_engergy}
\mathcal{E}(u)=\int_0^L \left(\frac{|\partial_x u (x)|^2}{2} +V(u(x))\right) \d x, 
\ee
 with the potential $V$ given by the standard double-well function $V(u)=(u^2-1)^2/4$. Indeed, one can readily check that $- D\mathcal{E}(u) =  \partial_{x}^2 u   + u-u^3$, where $D\mathcal{E}$ denotes the Fr\'echet derivative of $\mathcal{E}$ at $u$.
 The sACE provides a general framework for studying pattern formation and interface dynamics in systems ranging from liquid crystals and ferromagnets to tumor growth and cell membranes \cite{hohenberg1977theory,tauber2014critical}. 

Its universality makes it relevant across diverse fields, contributing to a unified understanding of these complex phenomena. By incorporating the spatio-temporal noise  term, $\d \W$, into the deterministic Allen-Cahn equation, the sACE accounts for inherent fluctuation and uncertainty present in real systems \cite{bray2002theory}. This allows for a more realistic description of transitions between phases, providing valuable insights beyond what pure deterministic models can offer. Near critical points, where phases coexist, the sACE reveals the delicate balance between deterministic driving forces and stochastic fluctuations. Studying these transitions helps us understand critical phenomena in diverse systems, from superfluid transitions to magnetization reversal in magnets \cite{hohenberg1977theory,bray2002theory}. 

Phase transitions for the sACE and its variants have indeed been analyzed both theoretically utilizing the large deviation principle \cite{FW12} and also numerically through different rare event algorithms. For instance, in \cite{berglund2013sharp} the expectation of the transition time between two metastable states is rigorously derived for a class of 1D parabolic SPDEs with bistable potential that include sACE as a special case. The phase diagram and various transition paths have also been computed by either an adaptive multilevel splitting algorithm \cite{Rolland_al16} or a minimum-action method \cite{zakine2023minimum}, where the latter deals with a nongradient flow generalization of the sACE that includes a constant force and also a nonlocal term.

Despite its importance in numerous physical applications, the reduction problem for the accurate reproduction of transition paths in the sACE has received limited attention, with existing works focusing primarily on scenarios near instability onset and low noise intensity (e.g., \citep{Caraballo_al07,BM13,wang2013macroscopic,blomker2020impact,Gao24}).

\begin{figure}
\centering
\includegraphics[width=\linewidth]{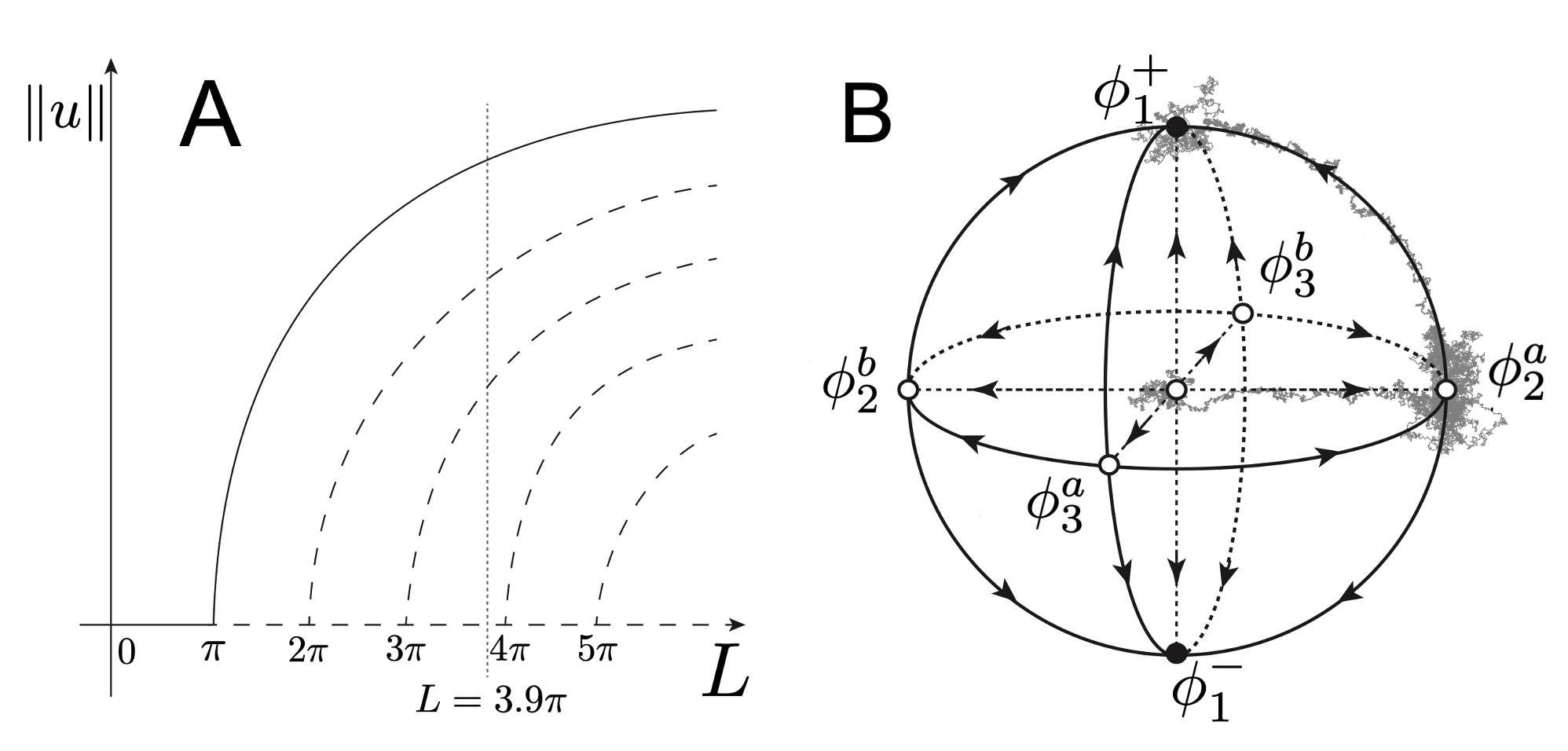}
\caption{{\bf Panel A:} Bifurcation diagram of from Eq.~\eqref{eq:Chafee} (without noise) as $L$ is varied; see \cite{bai1993numerical}.  The dashed vertical line marks the parameter regime analyzed here.  {\bf Panel B:} Schematic of the deterministic flow structure for this parameter regime. Are shown here, the deterministic global attractor on which is superimposed a solution path to the sACE.}
\label{fig:Chafee_bif_combo}
\end{figure}

In this study, we consider the sACE placed away from instability onset, after several bifurcations have taken place where multiple steady states coexist; see Fig.~\ref{fig:Chafee_bif_combo}. It corresponds to the  parameter regime
\be \label{Chafee_param}
L=3.9\pi, \; q=4,\; N=8, 
\ee 
in which the parameter $q$ and $N$ indicates the modes forced stochastically in Eq.~\eqref{eq:Chafee}, according to
\be  \label{eq:Wt_highOnly}
\W(x,\omega) = \sum_{j={q+1}}^{N} \sigma_j  W_t^j(\omega)\boldsymbol{e}_j(x), \; t \in \mathbb{R},
\ee
where we set $\sigma_j=\sigma= 0.2$ for $q+1\leq j \leq N$.
Note that as in Section \ref{Sec_LIA_cubic_case}, $\boldsymbol{e}_j$, denote the eigenmodes of the operator  $Au= \partial_{x}^2 u+u$ (with homogeneous Dirichlet boundary conditions), given here by
\be
 \boldsymbol{e}_j(x) = \sqrt{\frac{2}{L}} \sin\left(\frac{j\pi x}{L}\right), \quad j \in \mathbb{N},
\ee
with corresponding eigenvalues $\lambda_j = 1 - j^2\pi^2/L^2$.
We refer to Appendix \ref{Sec_numACE} for the numerical details regarding the simulation of Eq.~\eqref{eq:Chafee}.

More specifically, the domain size $L$ is chosen so that when $\sigma=0$ the trivial steady state has  experienced three successive bifurcations leading thus to a total of seven steady states with two stable nodes ($\phi_1^{+}$ and $\phi_1^{-}$) and five saddle points ($0$, $\phi_2^{a}$, $\phi_2^{b}$, $\phi_3^{a}$, and $\phi_3^{b}$). These steady states are characterized here by their zeros: the saddle states are of sign-change type while the stable ones are of constant sign; see \cite{Chafee_al74,bai1993numerical}. Note also that the same conclusion holds for the case with either Neumann boundary conditions or periodic boundary conditions; see e.g.~\cite[Figure 2]{berglund2013sharp}. 
In our notation, $a$ and $b$ denote the string made of sign changes over $(0,L)$, with leftmost symbol in $a$ to be $+$. For instance, for $\phi_3^{a}$ (resp.~$\phi_3^{b}$), $a=\{+,-,+\}$ (resp.~$b=\{-,+,-\}$), while $a=\{+,-\}$ (resp.~$b=\{-,+\}$) for $\phi_2^{a}$ (resp.~$\phi_2^{b}$).

By analyzing the energy distribution of these steady states carried out by the eigenmodes (see Table \ref{Table_Energy_frac}),  one observes that  $\phi_1^{\pm}$ are nearly collinear to $\bm{e}_1$, while $\phi_2^{a/b}$ are nearly collinear to  $\bm{e}_2$, and $\phi_3^{a/b}$ to  $\bm{e}_3$. 
\begin{table}[tbh!]
\caption{Energy distribution of the steady states  across the first few eigenmodes}
\label{Table_Energy_frac}
\centering
\begin{tabular}{ccccccccc}
\toprule\noalign{\smallskip}
 &  $\bm{e}_1$ & $\bm{e}_2$  &  $\bm{e}_3$ & $\bm{e}_4$ & $\bm{e}_5$ &  $\bm{e}_6$ & $\bm{e}_7$ &  $\bm{e}_8$   \\
\noalign{\smallskip}\hline\noalign{\smallskip}
$\phi_1^{\pm}$ & 94.69\% & 0 & 4.79\% & 0 &  0.46\% & 0 & 0.05\% & 0 \\
$\phi_2^{a/b}$ & 0 & 99.16\% & 0 & 0 &  0 & 0.83\% & 0 & 0 \\
$\phi_3^{a/b}$ & 0 & 0 & 99.93\% & 0 & 0  & 0 & 0 & 0  \\
\noalign{\smallskip} \bottomrule
\end{tabular}
\end{table}

A sketch of the bifurcation diagram as $L$ is varied (for $\sigma=0$) is shown in Fig.~\ref{fig:Chafee_bif_combo}A; the vertical dashed line marks the domain size $L=3.9\pi$ considered in our numerical simulations. For this domain size, the system exhibits three unstable modes whose interplay with the nonlinear terms yields the  deterministic flow structure depicted in Fig\ref{fig:Chafee_bif_combo}B. There, the deterministic global attractor is made of the seven steady states mentioned above (black filled/empty circles) connected by heteroclinic orbits shown by solid/dashed curves with arrows. To this attractor, we superimposed a solution path to the stochastic model \eqref{eq:Chafee} emanating from $u_0= 0$ (light grey ``rough'' curve).   We adopt the Arnol'd's convention for this representation: black filled circles indicate the stable steady states while the empty ones indicate the saddle/unstable steady states.

 These steady states organize the stochastic dynamics when the noise is turned on: a solution path to  Eq.~\eqref{eq:Chafee} transits typically  from one metastable states to another; the latter corresponding to profiles in the physical space that are   random perturbations of an underlying steady state's profile; see insets of Fig.~\ref{Fig_potential} below.  For the parameter setting \eqref{Chafee_param}, we observe though, over a large ensemble of noise paths,  that the SPDE's dynamics exhibit transition paths connecting $u_0=0$ to neighborhoods of the steady states $\phi_1^{\pm}$ and $\phi_2^{a/b}$, but does not meanders near  $\phi_3^{a/b}$ over the interval $[0,T]$ considered ($T=40$). 

\begin{figure}[htbp]
\centering
\includegraphics[width=1\linewidth, height=.5\textwidth]{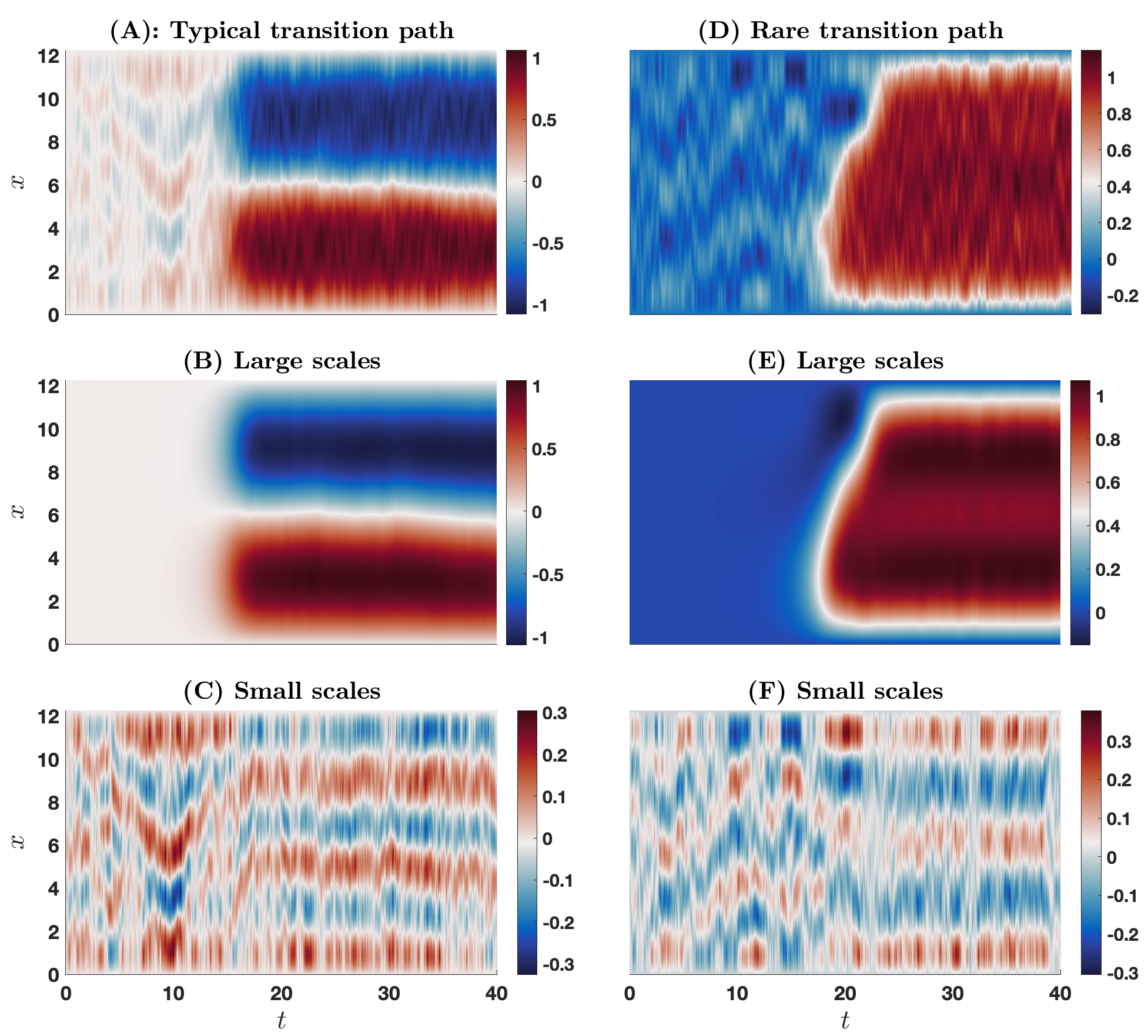}
\caption{{\bf Typical and rare transition paths for Eq.~\eqref{eq:Chafee}}. One path connecting $u_0=0$ to a sign-changing metastable state (typical, left panel), and one path connecting  $u_0=0$ to a constant-sign metastable state (rare event, right panel). Are shown also the large ($u_\c$) and small ($u_\s$) scale contributions.} \label{Fig_spacetime}
\end{figure}

Noteworthy is the visual rendering of the sojourn time of the solution path  shown in  Fig.~\ref{fig:Chafee_bif_combo}B. It is more substantial  in the neighborhood of $\phi_2^{a}$ than in that of $\phi_1^{+}$. It is a depiction  of a solution path to  Eq.~\eqref{eq:Chafee}  over  $[0,T]$: most of the trajectory spends time near a sign-change metastable state ($\phi_2^a$) rather than to the constant-sign one ($\phi_1^+$). Large ensemble simulations of the sACE for the parameter regime \eqref{Chafee_param} reveals that what is 
shown for a single path here is actually observed at the final time $t=T$, when the noise path is varied: the sign-change states $\phi_2^{a/b}$ are the rule whereas the  constant-sign states $\phi_1^{\pm}$ are the exceptions, corresponding here to rare events. 
Figure \ref{Fig_spacetime} shows  two solution paths in the space-time domain connecting  $u_0=0$ to either a sign-changing metastable state (typical) or a  constant-sign metastable state (rare event). 
In Section \ref{Sec_unravel} below we shall discuss the mechanisms behind these patterns and their transition phenomenology, and  the extent to which non-Markovian reduced models described below are able to reproduce it.

\subsection{Non-Markovian optimal reduced models}\label{Sec_above_forcingscaleACE}
Recall that our general goal is the derivation of reduced models able to reproduce the emergence of noise-driven spatiotemporal patterns and predict their transitions triggered by the subtle coupling between noise and the nonlinear dynamics.

For the sACE \eqref{eq:Chafee} at hand, this goal is declined into the following form: To derive reduced models  not only able to accurately estimate ensemble statistics (averages across many realizations) of the states reached at a given time $t=T$, but also to anticipate the system's typical metastable states as well as those reached through the rarest transition paths pointed above.  

To address this issue, we place ourselves in the challenging situation for which the cutoff scale is larger than the scales forced stochastically; i.e.~the case where noise acts only in the unresolved part of the system.  We show below that stochastic parameterization framework of Section \ref{Sec_OPMwhite} allows us to derive such an effective reduced model. This model is built from the stochastic parameterization \eqref{eq:h1tau_stoch} obtained from the BF system \eqref{Eq_BF_SPDEcubic} applied to Eqns.~\eqref{eq:Chafee} and \eqref{eq:Wt_highOnly}, that is optimized---on a {\it single training path}---following Section \ref{Sec_Optim}. We describe below the details of this derivation.

Thus, we take the reduced state space to be spanned by the  $q$ first unforced eigenmodes, i.e.
\be
H_{\c} = \mathrm{span}\{\boldsymbol{e}_1,\cdots,\boldsymbol{e}_q\}.
\ee 
In other words,  the  reduced state space is spanned by low modes that are not forced stochastically; i.e.~$q=m_c$ in \eqref{eq:Wt_highOnly}. 
The subspace $H_{\s}$ is taken to be the orthogonal complement of $H_{\c}$ in $L^2(0,L)$ and contains the stochastically forced modes; see \eqref{eq:Wt_highOnly}.

Equation \eqref{eq:Chafee} fits into the stochastic parameterization framework of Section \ref{Sec_OPMwhite} (see Eq.~\eqref{Eq_SPDE_cubic}) with $Au = \partial_x^2 u + u$ subject to homogenous Dirichlet boundary conditions, $G_2(u) = 0$  and 
with $G_3$ denoting the cubic {\it Nemytski'i operator} defined as $G_3(u,v,w)(x)=-u(x)v(x)w(x)$ for any $u,v$ and $w$ in $L^2(0,L)$, $x$ in $(0,L)$.

The resulting stochastic parameterization $\Phi_n$
given by \eqref{eq:h1tau_stoch} (with $Y=0$) takes then the following form: 
\bea\label{eq:h1_tau_Chafee}
\Phi_{n}& (\tau,  X, t; \omega)  \\
& =   \sigma_n \big(W^n_{t}(\omega) -  e^{\tau \lambda_n} W^n_{t-\tau}(\omega) \big)   + Z^{n}_{\tau}(t; \omega)  \\
&  \qquad + \sum_{i,j,k = 1}^q \Big( E_{ijk}^n(\tau)  C_{ijk}^n  X_{i} X_{j}X_{k} \Big),
\eea
for each $n$ in $(q+1,N)$. Here $E_{ijk}^n$ is given in \eqref{Eq_Ehat},  and $Z^{n}_{\tau}(t; \omega)$ is given by \eqref{Eq_Zn_decomp}.  

Since the linear operator $A$ is here self-adjoint, we have $\boldsymbol{e}^*_n = \boldsymbol{e}_n$, and the coefficients $C_{ijk}^n$ are given by
\bea \label{Eq_coef_C}
C_{ijk}^{n} &= -\langle \boldsymbol{e}_{i} \boldsymbol{e}_{j} \boldsymbol{e}_{k}, \boldsymbol{e}_n \rangle \\
& = -\int_0^L  \boldsymbol{e}_{i}(x) \boldsymbol{e}_{j}(x) \boldsymbol{e}_{k}(x) \boldsymbol{e}_{n}(x) \d x, 
\eea
where $1\leq i,j,k\leq q$, and $q+1 \leq n\leq N$.

For each $n$ in $(q+1,N)$, the free parameter $\tau$  in \eqref{eq:h1_tau_Chafee} is optimized over a {\it single, common, training path}, $u(t)$, solving the sACE, by  minimizing $\mathcal{Q}_n(\tau)$ in \eqref{Eq_minQnHn} with  $\Phi_n$ given by \eqref{eq:h1_tau_Chafee}; see Section \ref{Sec_Optim}. To do so, we follow the procedure described in Section \ref{Sec_Zn_simul} to simulate  efficiently 
the required random random coefficients (here the $Z^n$-terms), solving in particular the random ODEs \eqref{Eq_for_I}.

Once $\bftau^* = (\tau^*_{q+1},\ldots,\tau^*_N)$ is obtained, the corresponding $q$-dimensional optimal reduced system in the class $\mathcal{P}$ (see Section \ref{Sec_classP}),
is given component-wisely by  
\be \label{Eq_reduced_chafee_SI}
\dot{y}_i =   \lambda_i y_i - \Big\langle \Big( \sum_{j=1}^q y_j \boldsymbol{e}_{j} + \sum_{n=q+1}^N \Phi_{n}(\tau^*_{n}, \bm{y}, t; \omega) \boldsymbol{e}_{n} \Big)^3, \boldsymbol{e}_{i} \Big\rangle, 
\ee
where  $i = 1, \cdots, q$ and $\boldsymbol{y}=(y_1, \cdots, y_q)^{\textrm{T}}$.

The optimal reduced system \eqref{Eq_reduced_chafee_SI}, also referred to as the OPM closure hereafter, can be further expanded as
\begin{widetext}
\bea\label{Eq_closure}
\dot{y}_i & =  \lambda_i y_i + \sum_{j,k,\ell=1}^q C_{jk\ell}^i  y_{j} y_{k} y_{\ell} \\
& + 3  \sum_{j,k=1}^q  \sum_{n=q+1}^N C_{jkn}^i y_j y_k  \Phi_{n}(\tau^*_{n},  \bm{y},  t; \omega) \\
& +3  \sum_{j=1}^q  \sum_{n, n'=q+1}^N C_{jnn'}^i y_j \Phi_{n}(\tau^*_{n}, \bm{y}, t; \omega)\Phi_{n'}(\tau^*_{n'}, \bm{y}, t; \omega) \\
& +  \sum_{n, n', n''=q+1}^N C_{nn'n''}^i \Phi_{n}(\tau^*_{n}, \bm{y}, t; \omega) \Phi_{n'}(\tau^*_{n'}, \bm{y}, t; \omega) \Phi_{n''}(\tau^*_{n''}, \bm{y}, t; \omega).
\eea
\end{widetext}
Here also, the coefficients $C_{ijk}^{\ell}$ are those defined in \eqref{Eq_coef_C} for $1\leq i,j,k,\ell \leq N$.
Note that the required $Z^n$-terms to simulate Eq.~\eqref{Eq_closure}  are advanced in time by solving  the corresponding Eq.~\eqref{Eq_for_I}. A total of $N-q$ such equations are solved, each corresponding to a single random coefficient per mode to parameterize. 
Thus,  a total of, $N=q+(N-q)$, ODEs with random coefficients are solved to obtain a reduced model of the amplitudes of the $q$ first modes. We call Eq.~\eqref{Eq_closure} together with its auxiliary equations \eqref{Eq_for_I}, the optimal reduced model for the sACE.

A compact form of Eq.~\eqref{Eq_closure}  can be written as follows:
\bea\label{Eq_reduced_chafee}
 \dot{\bm{y}} =\Pi_\c A \bm{y} + \Pi_\c G_3\left(\bm{y}+\Phi_{\bftau^\ast}(\bm{y},t; \omega)\right),  \; \bm{y} \in H_{\c},\\
 \Phi_{\bftau^\ast}(\bm{y},t;\omega)=\sum_{n= q+1}^N  \Phi_{n}(\tau_{n}^\ast,  \bm{y}, t; \omega)   \boldsymbol{e}_n,
\eea
with $\Phi_n$ defined in Eq.~\eqref{eq:h1_tau_Chafee}. We refer to Appendix \ref{Sec_numACE} for the numerical details regarding the simulation of  Eq.~\eqref{Eq_reduced_chafee}.

Noteworthy is that an approximation of the sACE solution paths can be obtained from the solutions  $\boldsymbol{y}(t,\omega)$ to the optimal reduced model. Such an approximation consists of simply lifting the surrogate low-mode amplitude $\boldsymbol{y}(t,\omega)$  simulated from Eq.~\eqref{Eq_reduced_chafee}, to $H_\s$ via the optimal parameterization, $ \Phi_{\bftau^\ast}$, defined in Eq.~\eqref{Eq_reduced_chafee}. This approximation takes then the form
\bea \label{Eq_uPM}
u^\mathrm{opm}(t,x, \omega) & = \sum_{j=1}^q  y_j(t,\omega) \boldsymbol{e}_j(x)  \\
& +  \sum_{n=q+1}^N \Phi_{n}(\tau^*_{n}, \bm{y}(t,\omega), t; \omega) \boldsymbol{e}_n(x).  
\eea
The insets of Fig.~\ref{Fig_min_max}A below show the approximation skills achieved by this formula obtained thus only by solving the reduced equation  Eq.~\eqref{Eq_reduced_chafee}.

\begin{figure}[htbp]
\centering
\includegraphics[width=.5\textwidth,height=.4\textwidth]{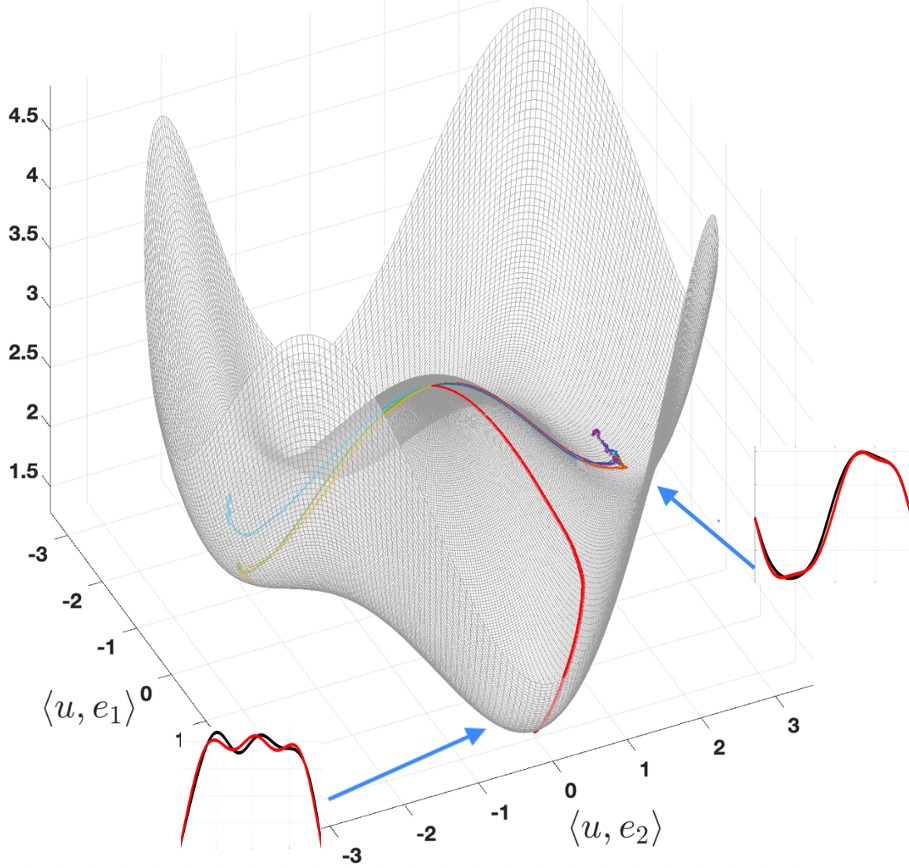} %
\caption{{\bf Ginzburg-Landau free energy functional}. Here, the Ginzburg-Landau free energy functional $\mathcal{E}(u)$ given by \eqref{GL_engergy} is shown over the reduced state made of the two dominant modes' amplitudes. A few solution paths navigating across this energy's landscape (over $[0,T]$) are shown. The ``shallower" local minima correspond to $\phi_2^{a/b}$ (sign-change profile), while the  deeper ones correspond to $\phi_1^{\pm}$ (constant-sign profile). The insets show examples of metastable states, reached at $t=T$ by the sACE, located near these steady states.  The rare paths are those trapped in the deeper wells due to the degenerate noise employed here; see Text.}
\label{Fig_potential}
\end{figure}

\begin{figure}[htbp]
\centering
\includegraphics[width=0.48\textwidth, height=0.4\textwidth]{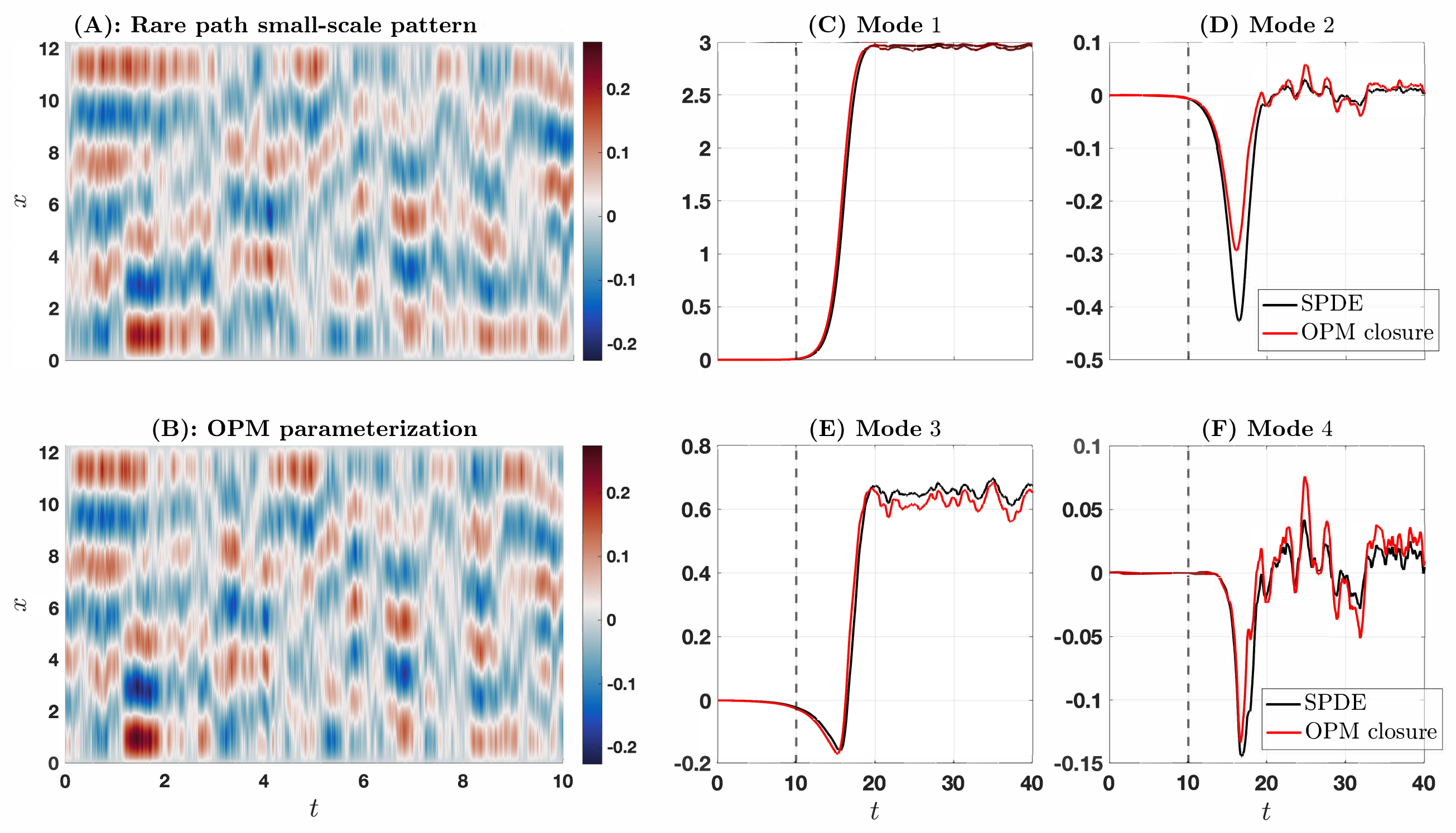}
\caption{{\bf Non-Markovian terms and rare path's small-scale component.} Panel A shows the small-scale component of  the solution to  the sACE,  experiencing a rare transition path. Panel B shows the OPM parameterization of this small-scale component over the interval (0,10), prior the transition takes place.
Panels C to F show  the ability of the optimal  reduced model, Eq.~\eqref{Eq_closure},  to predict the large fluctuations experienced by the large-scale modes' amplitudes, beyond this  interval.}\label{Fig_Importance_Non-Markov}
\end{figure}

\begin{figure}[htbp]
\centering
\includegraphics[width=0.48\textwidth, height=0.4\textwidth]{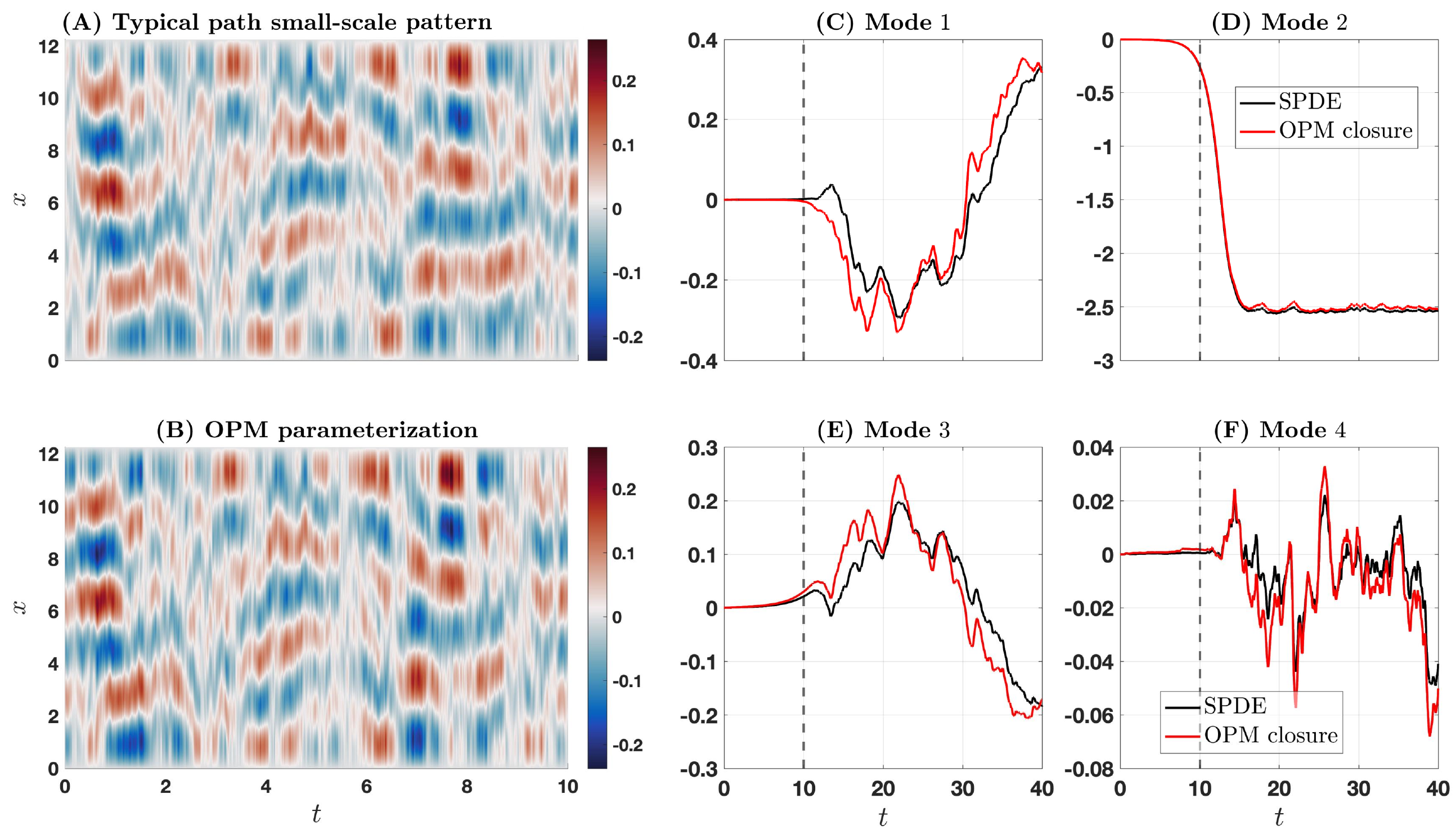}
\caption{{\bf Non-Markovian terms and typical path's small-scale component.} Same as Fig.~\ref{Fig_Importance_Non-Markov} but for a typical path.} 
\label{Fig_Importance_Non-Markov2}
\end{figure}


\subsection{Why energy flows in a certain direction: The role of non-Markovian effects}\label{Sec_unravel}
 
The sACE \eqref{eq:Chafee} is integrated over a large ensemble of noise paths ($10^6$ members), from $u_0=0$ at $t=0$ up to a time $t=T$. 
In our simulations with a fixed time horizon ($T = 40$), reaching a constant-sign state at the end ($t = T$) is a rare event, even though these states correspond to the deeper energy wells in the system's landscape (as shown in Fig.~\ref{Fig_potential}). Conversely, sign-change states, which reside in shallower wells, are more common. This might seem counterintuitive from an energy perspective. Nevertheless, this can be explained by the degenerate noise (Eq.~\eqref{eq:Wt_highOnly}) used to force the sACE: only a specific group of  unresolved (stable) modes is excited by noise and this may limit the system's ability to explore the deeper energy wells efficiently within the simulation timeframe.

Over time (as $t$ increases), the noise term in the sACE interacts with the nonlinear terms, gradually affecting the unforced resolved modes.
Our experiments show that the amplitude of the sACE's solution amplitude  carried by $\bm{e}_2$ grows the fastest for most noise paths.  Because $\phi_2^{a/b}$ is nearly aligned with $\bm{e}_2$ (as shown in Table \ref{Table_Energy_frac}), this means that for most noise paths, the solution is expected to reach a metastable state close to $\phi_2^{a/b}$ when the nonlinear effects become significant (nonlinear saturation kicks in).

Understanding why energy primarily accumulates in  the $\bm{e}_2$ direction is complex. It stems from how noise initially affecting unresolved modes (modes $5$ to $8$) propagates through nonlinear interactions, gradually influencing larger-scale resolved modes. As the initial state ($u_0=0$) lacks energy, this noise-induced energy transfer to larger scales takes some time. Our analysis reveals that during the interval $I = (0, 10)$, the solution is dominated by these small-scale features inherited from the noise forcing (Table \ref{Table_Energy_small_scales} and Fig.~\ref{Fig_spacetime}).

This finding has significant implications for closure methods: to accurately reproduce both rare and typical transitions, the underlying parameterization must effectively represent the system's small-scale response to stochastic forcing before the transition occurs.

Figures \ref{Fig_Importance_Non-Markov}  and \ref{Fig_Importance_Non-Markov2} illustrate how the optimal parameterization and reduced model (Eq.~\eqref{Eq_reduced_chafee}) capture the subtle differences in the small-scale component that ultimately lead to rare or typical events in the full sACE system. Notably, for both rare (Figure \ref{Fig_Importance_Non-Markov}) and typical paths (Figure \ref{Fig_Importance_Non-Markov2}), the low-mode amplitudes of the sACE solution are near zero during the interval $I = (0, 10)$ (Panels C-F and Table \ref{Table_Energy_small_scales}). Therefore, it is the spatial patterns of the small-scale component, rather than its energy content, that control the triggering of rare or typical events (Figure \ref{Fig_spacetime}A and \ref{Fig_spacetime}D). The successful parameterization of these small-scale features (Figures  \ref{Fig_Importance_Non-Markov}B and  \ref{Fig_Importance_Non-Markov2}B) is crucial for the optimal reduced model's ability to reproduce the key behaviors of the sACE's large-scale dynamics.

The non-Markovian terms are essential for accurately capturing these small-scale patterns before a rare or typical event occurs. When large-scale mode amplitudes are near zero during interval I (Table \ref{Table_Energy_small_scales}), the optimal parameterization (Eq.~\eqref{Eq_reduced_chafee})  is dominated by the (optimized) non-Markovian field $\xi_t$ given by:
\bea\label{Non-Markov_field}
&\xi_t(\omega,x)=\sum_{n=q+1}^N M_n(t,\tau_n^\ast,\omega) \bm{e}_n (x), \mbox{ with} \\
&M_n(t,\tau_n^\ast,\omega)=\sigma_n \big(W^n_{t}(\omega) -  e^{\tau_n^\ast \lambda_n} W^n_{t-\tau_n^\ast}(\omega) \big) \\
&\hspace{5cm} + Z^{n}_{\tau_n^\ast}(t; \omega).
\eea
After this initial transient period, nonlinear effects become significant, and the large-scale mode amplitudes undergo substantial fluctuations. Panels C-F in Figures  \ref{Fig_Importance_Non-Markov} and  \ref{Fig_Importance_Non-Markov2} demonstrate this behavior and the optimal reduced model's ability to predict these large-scale fluctuations beyond interval I. The OPM reduced system's ability to capture such transitions is rooted in its very structure. Its coefficients in Eq.~\eqref{Eq_closure}  are nonlinear functionals of the non-Markovian $M_n$-terns within $\xi_t$ allowing for an accurate representation of the genuine nonlinear interactions between noise and nonlinear terms in the original sACE, which drive fluctuations in the large-mode amplitudes (Modes 1 to 4).

After a transition occurs and energy is transferred to the low modes, the stochastic component $\xi_t$ and its non-Markovian terms become secondary while  the nonlinear terms in the stochastic parameterization \eqref{eq:h1_tau_Chafee}  become crucial for capturing the average behavior of the sACE dynamics, as highlighted in Section \ref{Sec_weak_timescale} below.

\begin{table}[tbh!]
\caption{Energy distribution over $(0,10)$ for the sACE's rare and typical paths}
\label{Table_Energy_small_scales}
\centering
\begin{tabular}{ccccccccc}
\toprule\noalign{\smallskip}
 &  $\bm{e}_1$ & $\bm{e}_2$  &  $\bm{e}_3$ & $\bm{e}_4$ & $\bm{e}_5$ &  $\bm{e}_6$ & $\bm{e}_7$ &  $\bm{e}_8$  \\
\noalign{\smallskip}\hline\noalign{\smallskip}
Rare path & 0 & 0.05\% & 0.01\% & 0 &  49.6\%& 26.1\%  & 15.1\% & 9.2\% \\
Typical path & 0 & 0.22\% & 0 & 0 &  34.8\%& 33\%  & 19.7\% & 12.4\% \\
\noalign{\smallskip} \bottomrule
\end{tabular}
\end{table}

\begin{figure*}[htbp]
\centering
\includegraphics[width=1\linewidth, height=.35\textwidth]{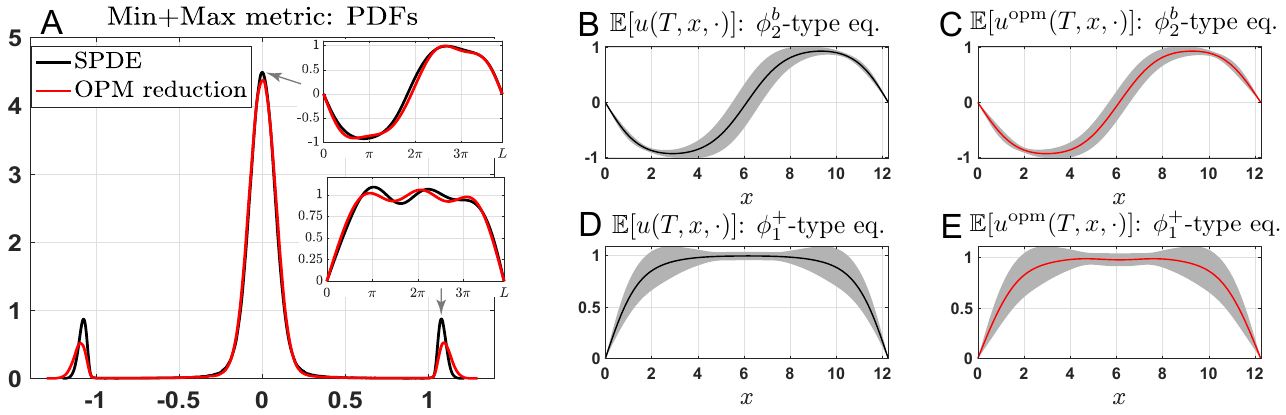}
\caption{{\bf Comparison of metastable state statistics between the full sACE model  (Eq.~\eqref{eq:Chafee}) and the optimal reduced model (Eq.~\eqref{Eq_reduced_chafee}) at time $t = T$}. Panels B and D show the distribution of $u(T, x, \cdot)$ from the sACE solution at $t = T$. Panels C and E show the corresponding approximation $u^{\mathrm{opm}}(T,x,\cdot)$ obtained from the optimal reduced model's solutions at $t = T$, leveraging Eq.~\eqref{Eq_uPM}. One million test paths are used to generate these statistics. See Text for more details.}\label{Fig_min_max}
\end{figure*}

\subsection{Rare and typical events statistics from non-Markovian reduced models}\label{Sec_stat_reduced}
Here, we evaluate the ability of our non-Markovian reduced model (Eq.~\eqref{Eq_reduced_chafee}) to statistically predict the rare and typical transitions discussed in the previous section. We test the reduced model's  performance in this task using a vast ensemble of one million noise paths.

\begin{figure}[htbp]
\centering
\includegraphics[width=1\linewidth,height=.25\textwidth]{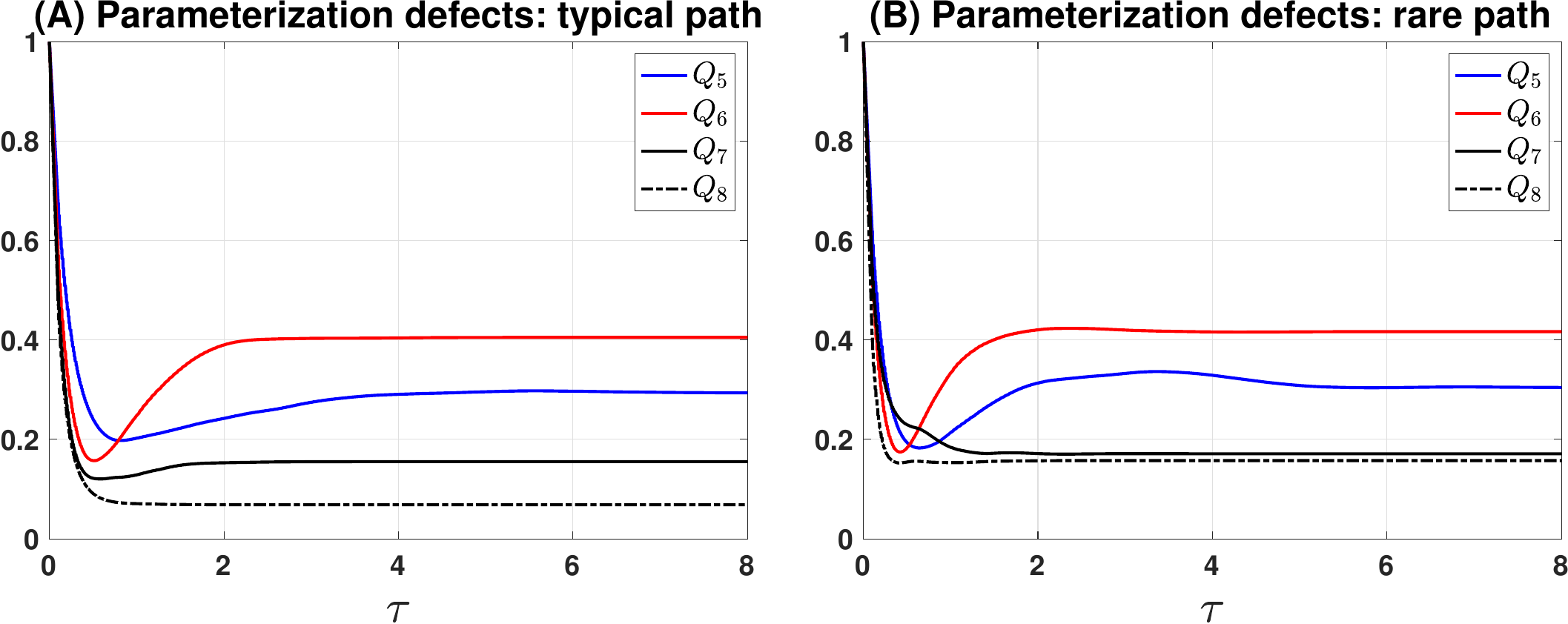}
\caption{{\bf The (normalized) parameterization defects, $Q_n(\tau)$.}  
{\bf Panel A:} The parameterization  $\Phi_n$ given by \eqref{eq:h1_tau_Chafee} is trained offline on a (single)  typical transition path taken over the interval $(t_s,T)$, with $t_s=10$ corresponding to the time at which nonlinear effects kick in (see Fig.~\ref{Fig_Importance_Non-Markov}). This optimized stochastic parameterization is used into the reduced system \eqref{Eq_reduced_chafee}  to produce the online ensemble results of Fig.~\ref{Fig_min_max}. {\bf Panel B:} The $\tau$-dependence of  the $Q_n(\tau)$ is shown on a rare transition path, for information purposes. See Section \ref{Sec_weak_timescale} for more details.}
\label{Fig_Qn_Chafee}
\end{figure}

Recall that the optimal parameterization ($\Phi_n$, defined in \eqref{eq:h1_tau_Chafee}) is trained on a single, ``typical" noise path (see Fig.~\ref{Fig_Qn_Chafee}A). In other words, we minimized the parameterization defect $\mathcal{Q}_n (\tau)$ in equation \eqref{Eq_minQnHn} for this specific path.  Our key objective here is to assess how well the reduced model (Eq.~\eqref{Eq_reduced_chafee}) can predict the distribution of final states (metastable states at time $T$) across a much larger set of out-of-sample noise paths.

To do so, we look at a natural ``min+max metric.''  This metric consists
of adding the minimum to the maximum values of the sACE's solution profile at time $T$. The shape of this probability distribution shows that a value close to zero corresponds to a typical sign-change metastable state, while a value close to $1$ or $-1$, corresponds to a rare constant-sign metastable state.

Figure \ref{Fig_min_max}A demonstrates the effectiveness of the non-Markovian optimal reduced model in capturing rare events. The  optimal reduced system (Eq.~\eqref{Eq_reduced_chafee}) can predict noise-driven trajectories connecting an initial state ($u_0=0$) to rare final states (sACE solution's profile with a constant sign at time $T$) with statistically significant accuracy (see Fig.~\ref{Fig_min_max}A). While these rare events occur less frequently in the reduced system compared to the full sACE (smaller red bumps compared to black ones near 1 and -1), the ability of the reduced equation Eq.~\eqref{Eq_reduced_chafee} to predict them highlights the optimal parameterization's success in capturing the system's subtle interactions between the noise and nonlinear effects responsible for these rare events.

Thus, the optimal reduced system reproduces faithfully the probability distribution of both common (frequently occurring sign-change profiles) and rare (constant-sign profiles) final states at time $T$ (Fig.~\ref{Fig_min_max}A). This is particularly impressive considering the model was trained on just a single path, yet generalizes well to predict the behavior across a large ensemble of out-of-sample paths.

Figure \ref{Fig_Qn_Chafee}A shows the training stage over a typical transition path. The optimization of the free parameter $\tau$ in the  stochastic parameterization $\Phi_n$ (see Eq.~\eqref{eq:h1_tau_Chafee}) is executed over a  typical transition path,  by minimizing the normalized parameterization defects $Q_n(\tau,T)= \mathcal{Q}_n(\tau,T)/\overline{\big|u_n(t) \big|^2}$ ($\mathcal{Q}_n(\tau,T)$ defined in  Eq.~\eqref{Eq_minQnHn}) for all the relevant $n$ (here $n=5,\cdots,8$).   
 We observe that $Q_5$ and $Q_6$ exhibit clear minima, whereas  the other modes have their parameterization defect not substantially diminished at  finite $\tau$-values compared to its asymptotic value (as $\tau\rightarrow \infty$). Same features are observed over a rare transition path; see Figure \ref{Fig_Qn_Chafee}B.
Recall that our reduced system is an ODE system whose stochasticity comes exclusively from its random coefficients as the resolved modes in the full system are not directly forced by noise.  It is this specificity that enables the reduced equations to effectively predict noise-induced transitions, including rare events. This remarkable capability relies on the optimized memory content embedded within the path-dependent, non-Markovian coefficients.

Thus, the optimized, non-Markovian and nonlinear reduced equation \eqref{Eq_closure} is able to track efficiently how the noise, acting on the ``unseen" part of the system, interacts with the resolved modes through these non-Markovian coefficients.  
While these coefficients effectively encode the fluctuations for a good reduction of the sACE system, they alone are insufficient to handle regimes with weak time-scale separation such as dealt with here.  Understanding the average behavior of the resolved dynamics is equally important. The next section explores how our parameterization, particularly its nonlinear terms, play a crucial role in capturing this average behavior.

\subsection{Weak time-scale separation, nonlinear terms and conditional expectation}\label{Sec_weak_timescale}
We assess the ability of our reduced system to capture this average dynamical behavior in the physical space. 
In that respect, Figures \ref{Fig_min_max}B and \ref{Fig_min_max}D show the expected profiles of the sACE's states obtained at $t=T$ (solid black curves) from the sACE \eqref{eq:Chafee}, along with their standard deviation (grey shaded areas), for the ensemble used for  Fig.~\ref{Fig_min_max}A. 
Similarly, Figures \ref{Fig_min_max}C and \ref{Fig_min_max}E show the expected profiles (solid red curves) with their   standard deviation  as obtained at $t=T$ from  $u^\mathrm{opm}(T,x, \cdot)$ defined in Eq.~\eqref{Eq_uPM}, once the low-mode amplitude is simulated by the optimal reduced model \eqref{Eq_reduced_chafee} and lifted to $H_\s$ by using the stochastic parameterization $ \Phi_{\bftau^\ast}$.
\begin{figure}[htbp]
\centering
\includegraphics[width=1\linewidth,height=.55\textwidth]{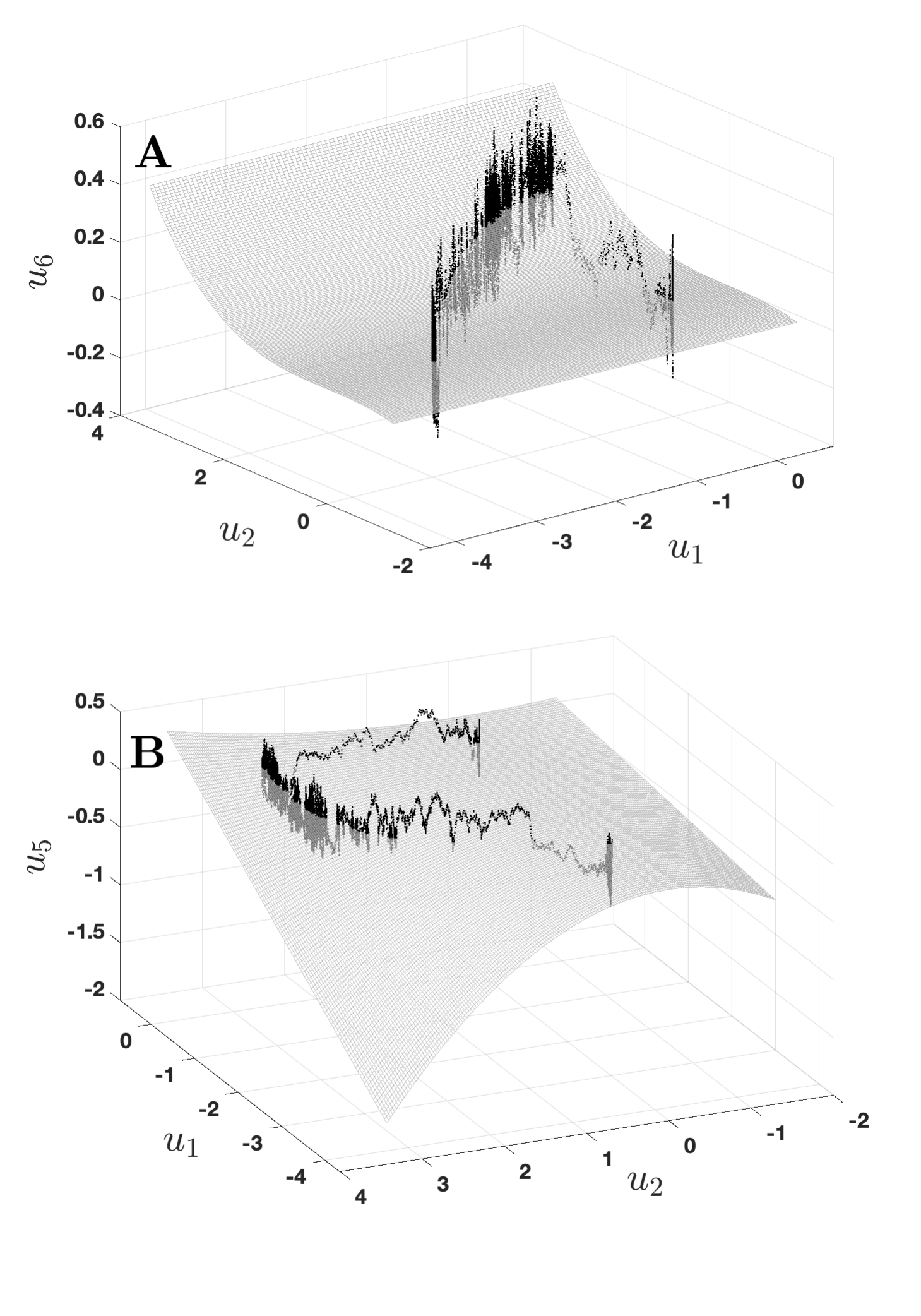} %
\caption{{\bf Markovian optimal parameterization for the 5th and 6th modes' amplitudes}. {\bf Panel A}: 
 The black curve shows the time evolution of the 6th mode's amplitude, $u_6(t)$, obtained by solving the governing equation, Eq.~\eqref{eq:Chafee}, over the time interval 0 to 400. This 6th mode's amplitude is shown against the corresponding first and second modes' amplitudes, $u_1(t)$ and $u_2(t)$. For visualization purposes, the other solution components in the four-dimensional subspace $H_{\c}$, $u_3(t)$ and $u_4(t)$, are fixed at their most probable values, $p_3$ and $p_4$, respectively. The manifold shown is the underlying Markovian OPM for the 6th mode's amplitude (Eq.~\eqref{Eq_Markovian_OPM} with $n=6$). 
{\bf Panel B}: Same for the 5th mode's amplitude $u_5(t)=\langle u(t),\bm{e}_5\rangle$ and its corresponding Markovian OPM (Eq.~\eqref{Eq_Markovian_OPM} with $n=5$).}
\label{Fig_toto_visual}
\end{figure}

By conducting large-ensemble (online) simulations of our reduced model (Eq.~\eqref{Eq_reduced_chafee}), we accurately reproduce the mean motion and second-order statistics of the original sACE. This success stems from our parameterization's ability to effectively represent the average behavior of the unresolved variables during the offline training phase, as defined by the probability measure  $\mu$ in \eqref{cond1}. Figure \ref{Fig_toto_visual} illustrates this point for the fifth and sixth modes' amplitudes of the sACE system.
Notably, the nonlinear nature of the Markovian optimal parameterization, evident after taking the expectation (see Eq.~\eqref{Eq_Markovian_OPM} below), highlights the crucial role of nonlinear terms in our stochastic parameterization (Eq.~\eqref{eq:h1_tau_Chafee}) in achieving this accurate representation.

 The manifold shown in  Fig.~\ref{Fig_toto_visual}A (resp.~Fig.~\ref{Fig_toto_visual}B) is obtained as the graph of the following Markovian OPM for the 6th (resp.~5th) mode's amplitude, given by the expectation of Eq.~\eqref{eq:h1_tau_Chafee} with $n=6$ (resp.~$n=5$) after optimization (Fig.~\ref{Fig_Qn_Chafee}), 
\bea\label{Eq_Markovian_OPM}
(X_1,X_2)\mapsto   & \mathbb{E} [\Phi_{n}(\tau_n^\ast,\bm{X})] \\
 & = \sum_{i,j,k = 1}^q \Big( E_{ijk}^n(\tau_n^*)  C_{ijk}^n  X_{i} X_{j}X_{k} \Big),
\eea
for $\bm{X}=(X_1,X_2,X_3,X_4)$ with $X_3=p_3$ and $X_4=p_4$, where  
$p_3$ and $p_4$ denote the  most probable values of the 3rd and 4th solution components, $u_3(t)$ and $u_4(t)$, to enable visualization as a 2D surface while favoring a certain "typicalness" of the visualized manifold. Note that the Markovian OPM corresponds simply to the deterministic cubic term in \eqref{eq:h1_tau_Chafee} after replacing $\tau$ by the optimal $\tau_n^*$, since the remaining stochastic terms equal to  $\sigma_n \int_{t-\tau_n^*}^t  e^{\lambda_n(t-s)}  \d W^n_{t}(\omega)$ (see Remark~\ref{Rmk_OU}), whose expectation is zero as stochastic integral in the sense of It\^o.

Without the optimization stage, the parameterizations of Section \ref{Sec_Invariance_Eq} experience deficient accuracy when the memory content is taken to be $\tau=\infty$, due to small values that the spectral gaps $\delta_{ijk}^n $ take in denominators. As a reminder, small spectral gaps are indeed a known limitation of traditional invariant manifold techniques, often leading to inaccurate results in this case. Optimizing the parameter $\tau$ helps address this issue. It introduces corrective terms like $(1-\exp(-\delta_{ijk}^n)\tau)$ in the cubic coefficients ($E_{ijk}^n$, Eq.~\eqref{Eq_Ehat}) that effectively compensate for these small gaps. This is especially important for fifth and sixth (unresolved) modes whose parameterization defects show a clear minimum during optimization (Figure \ref{Fig_Qn_Chafee}).

\begin{figure}[htbp]
\centering
\includegraphics[width=1\linewidth,height=.4\textwidth]{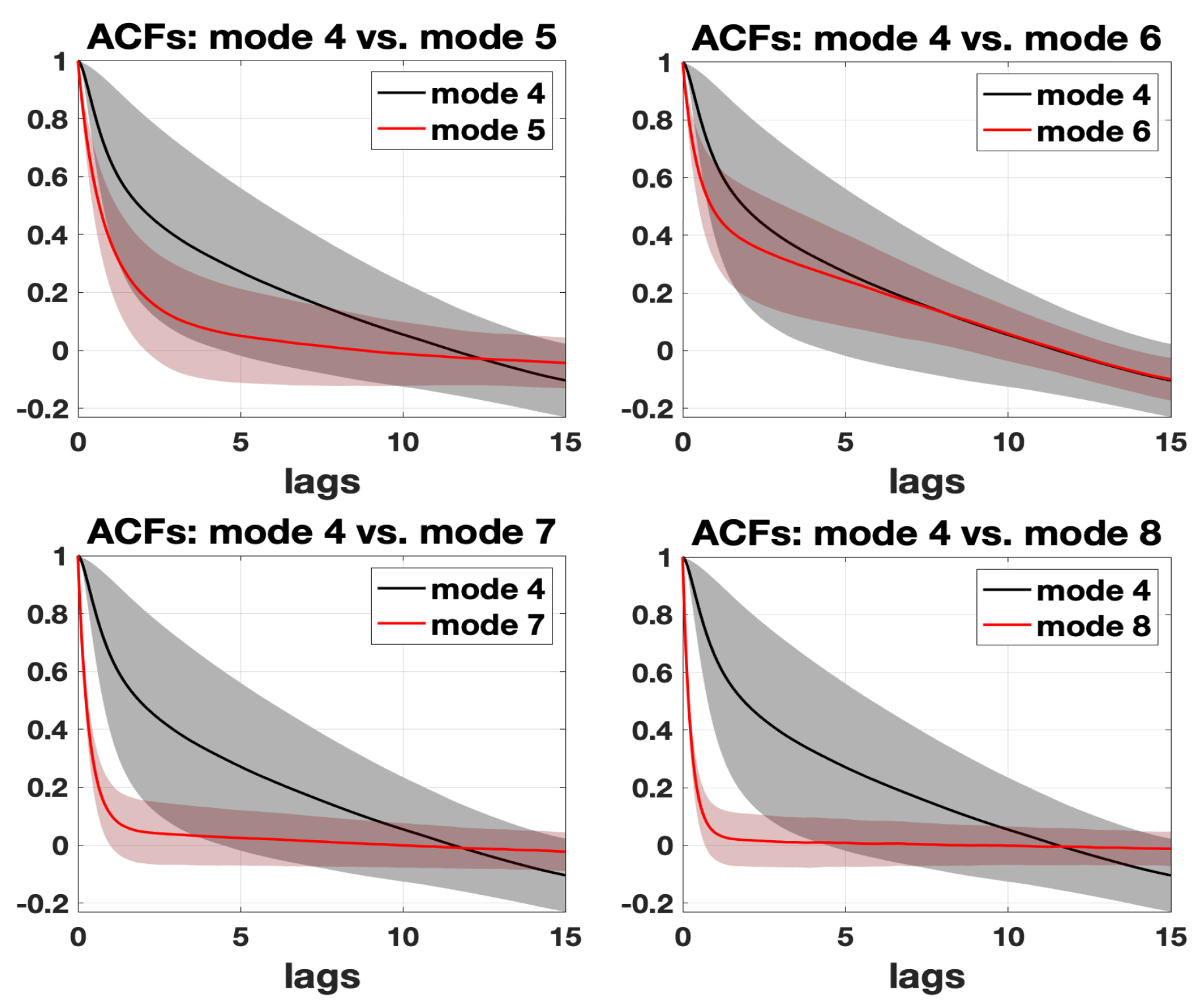} %
\caption{{\bf Autocorrelation functions (ACFs).}  Ensemble averages of the ACFs for the sACE's solution amplitude carried by the resolved mode $\bm{e}_4$ (black curve) and those carried by the unresolved modes $\bm{e}_5, \ldots, \bm{e}_8$ (red curves). The gray (resp.~pinky) shaded area around the black (resp.~red) curve shows the ACF's standard deviation for $\bm{e}_4$ (resp.~the unresolved modes). These statistics are computed from the same ensemble of sACE solution paths as in Fig.~\ref{Fig_min_max}.}  
\label{Fig_ACFs_Chafee}
\end{figure}

Small spectral gaps often translate into a more observable dynamical effect: weak time-scale separation between the resolved and unresolved variables near the cutoff scale. This is where the optimization process plays a particularly important role.  In our sACE system, we observe indeed that the timescales of the unresolved modes near cutoff have a greater impact on the optimization of the backward integration time $\tau$ (Fig.~\ref{Fig_Qn_Chafee}) from our backward-forward  system (Eq.~\eqref{Eq_BF_SPDEcubic}) which underpins the parameterization (Eq.~\eqref{eq:h1_tau_Chafee}). Modes like the fifth and sixth (unresolved) have correlation decay rates closer to the fourth resolved mode, compared to the much faster decay rates of the seventh and eighth modes (Fig.~\ref{Fig_ACFs_Chafee}). Notably, the fourth and sixth modes exhibit very similar decay rates (except for short time lags). This similarity makes finding the optimal value for $\tau$ crucial.
 Interestingly, the minimum in the optimization is more pronounced for the sixth mode compared to the fifth (Fig.~\ref{Fig_Qn_Chafee}).  This difference could be caused by non-linear interactions between these modes or their relative energy content. The seventh and eighth modes exhibit much faster correlation decay (Fig.~\ref{Fig_ACFs_Chafee}) and consequently their parameterization defects exhibit marginal optimization at finite $\tau$-values compared to their asymptotic values.

 Thus, the closure results shown in Figure \ref{Fig_min_max} demonstrate our reduction approach's effectiveness in handling weak time-scale separation scenarios, while also suggesting potential for further improvement (see Section \ref{Sec_discussion}).
Our data-driven optimization of analytical parameterization formulas offers a compelling alternative to traditional multiscale methods or those based on invariant manifolds \citep{Caraballo_al07,BM13,wang2013macroscopic,blomker2020impact,Gao24}. These methods, while providing rigorous error estimates, often require strong time-scale separation for the derivation of their amplitude equations. Additionally, they can achieve higher-order nonlinearities in their reduced equations too but through a different route:  the It\^o formula is applied after exploiting a slow timescale of the system (e.g., \cite[Eq.~(12)]{BM13} and \cite[Eq.~(15)]{Mohammed_al14}).  Integrating optimization concepts into these alternative frameworks could be a promising avenue to explore how it performs in weak time-scale separation regimes.

\section{Predicting Jump-Induced Transitions}\label{Sec_Jump}

\subsection{Spatially extended system with double-fold bifurcation}\label{Sec_jumpA}
We consider the following spatially extended system,
\bea\label{Eq_Sbif}
&\partial_{x}^2 u   + \lambda(1+u^2-\epsilon u^3)=0, \;\mbox{on } (0,L),\\
&u(0)=u(L)=0.
\eea
This system serves as a basic example of a spatially extended system exhibiting a double-fold bifurcation (also known as an S-shaped bifurcation). In this bifurcation, two equilibrium points collide and vanish, leading to tipping phenomena within a hysteresis loop.
 Similar elliptic problems to Eq.~\eqref{Eq_Sbif} with these properties emerge in various applications such as gas combustion and chemical kinetics \cite{bebernes2013mathematical,frank2015diffusion}, plasma physics and Grad-Shafranov equilibria in Tokamaks \cite{chandrasekhar1957introduction,grad1958hydromagnetic,Shafranov1958}, and 
gravitational equilibrium of polytropic stars \cite{chandrasekhar1957introduction,fowler1931further,hopf1931emden}. 
Beyond polynomial nonlinearities, other spatially extended systems are known to exhibit S-shaped bifurcation diagrams. Examples include: Earth's energy balance models \cite{ghil1976climate,north1981energy,hetzer1997s}, and oceanic models based on hydrostatic primitive or Boussinesq equations employed in the study of the Atlantic meridional overturning circulation (AMOC) \cite{thual1992catastrophe,dijkstra1997symmetry,weijer2019stability}.
We also note that the OPM approach has already been successfully applied to predict tipping phenomena in a prototype Stommel--Cessi model of the AMOC driven by Gaussian noise \cite{chekroun2023optimal}.

\begin{figure}
\centering
\includegraphics[width=.45\textwidth, height=0.32\textwidth]{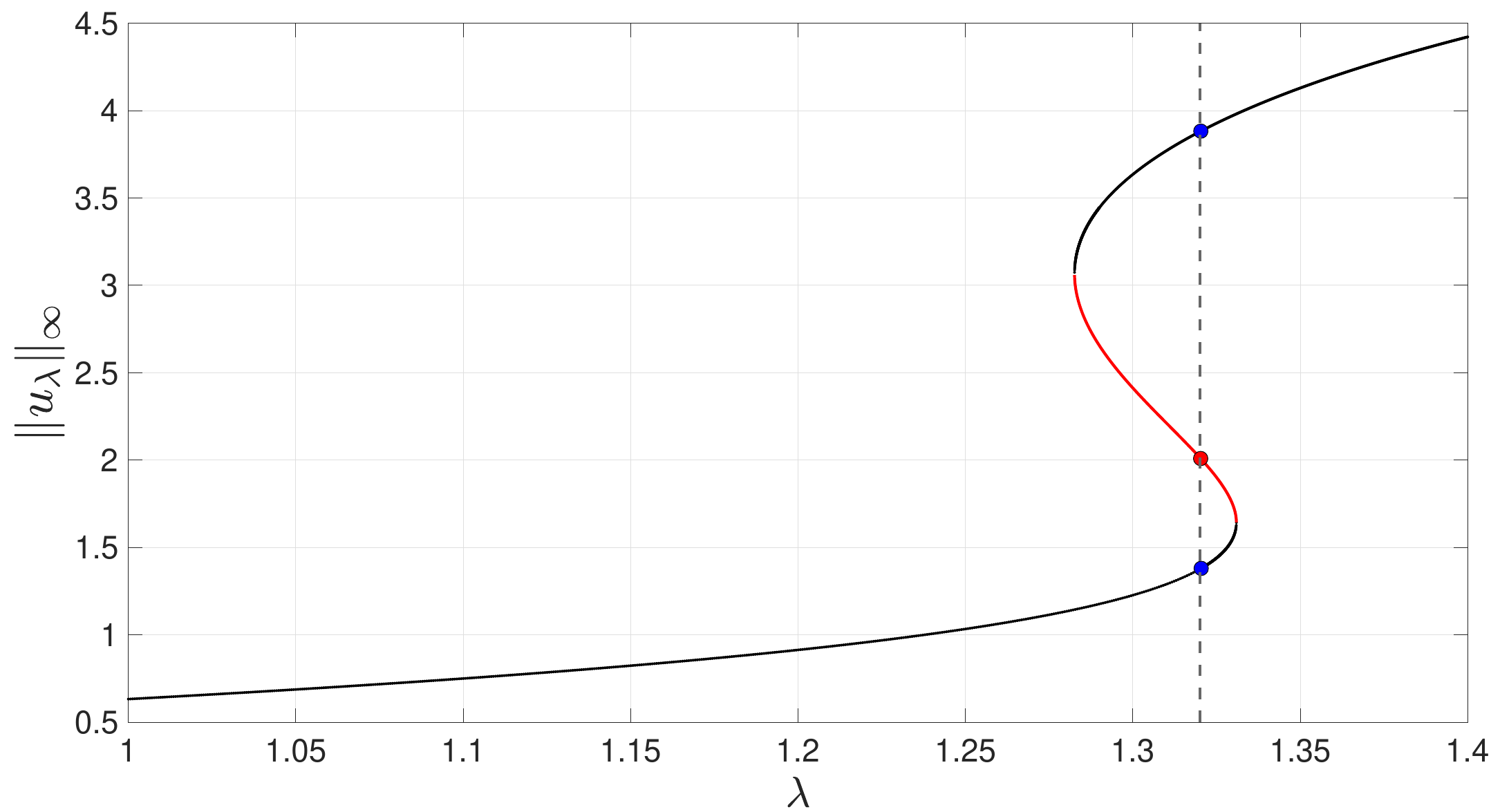}
\caption{{\bf  Double-fold bifurcation in Eq.~\eqref{Eq_Sbif}}. 
The red curve shows unstable steady states for $\epsilon=\epsilon^\ast/2$.
The black curves represent stable steady states.
The double-fold bifurcation occurs at the two turning points where the unstable states (red curve) collide with the stable states (black curves).
The vertical dashed line indicates the value $\lambda=1.32$ used in Table \ref{Table_Bimodal_regime} for stochastic  simulation in Section \ref{Sec_jump-driven}, and prediction experiments in Section \ref{Sec_OPM_Jumpcase}.} \label{Fig_Sbif}
\end{figure}

For equation \eqref{Eq_Sbif}, the bifurcation parameter is $\lambda$. 
For $L = 2$, it is known that the solution curve  $\lambda\mapsto u_\lambda$ exhibits an S-shaped bifurcation when $\epsilon <\epsilon^* \approx 0.3103$ \cite[Theorem 3.2]{korman1999exactness}.  This critical constant $\epsilon^\ast$ is derived from a general analytic formula applicable to a broader class of equations where  $u^3$ is replaced by $u^p$ with a fractional exponent $p > 2$ \cite{korman1999exactness}.
More precisely,  for $\epsilon$ satisfying 
\be
\epsilon^{\frac{p}{p-2}}<\frac{1}{3}\bigg[ \frac{2}{p(p-1)}\bigg]^{\frac{p}{p-2}} -\frac{p-1}{p+1}\bigg[ \frac{2}{p(p-1)}\bigg]^{\frac{p}{p-2}},
\ee
Theorem 3.2 of \cite{korman1999exactness} ensures that the solution curve to Eq.~\eqref{Eq_Sbif} (with $u^3$ replaced by $u^p$) is exactly S-shaped. 
We focus here on the case $p=3$.

To locate the interval of $\lambda$ over which an exact multiplicity of three steady states occurs (as predicted by the $S$-shape) we compute the bifurcation diagram of Eq.~\eqref{Eq_Sbif} for $\epsilon=\epsilon^\ast/2$;  see Appendix~\ref{Sec_Num_bif_diag} for details.
The result is shown in Figure \ref{Fig_Sbif}. The vertical dashed line at $\lambda=1.32$ marks the $\lambda$-value at which stochastic perturbations are included below to trigger jumps  across these  steady states.

\subsection{System's sensitivity to jumps, and jump-driven dynamics}\label{Sec_jump-driven}
In this section, we consider stochastic perturbations of jump type carried by some of the eigenmodes of  the linearized problem at the unstable steady state $U^\ast_\lambda$. The latter problem is given by
\bea\label{Eq_eigenpb}
& \partial_{x}^2 \psi +  \lambda (2{U^\ast_\lambda}  - 3  \epsilon {U^\ast_\lambda} ^2) \psi  = \beta \psi, \;\mbox{on } (0,L),\\
& \psi (0) = \psi (L) = 0,
\eea
and is solved below using a Chebyshev spectral method \cite[Chapters 6 \& 7]{trefethen2000spectral} for the numerical experiments.  Hereafter we denote by $(\beta_n,\bm{e}_n)$ the eigenpairs solving Eq.~\eqref{Eq_eigenpb}. Note that the eigenvalues are real and simple for this problem.

Thus, the jump-driven SPDE at the core of our reduction study is given by:
\bea\label{Eq_SPDEjump}
&\partial_t u -\partial_{x}^2 u   - \lambda(1+u^2-\epsilon u^3)=\lambda \eta_t,\\
&u(0)=u_0 \in H,
\eea
supplemented with Dirichlet boundary conditions. The ambient Hilbert space is $H=L^2(0,L)$ endowed with its natural  inner product denoted by $\langle \cdot, \cdot \rangle$.

Here, $\eta_t$ is a jump noise term that takes the form
\be\label{forcing_eta}
\eta_t(x)= \sigma \zeta_t f(t) (\bm{e}_3(x) +\bm{e}_5(x)), \;x\in (0,L),
\ee 
where $\sigma \geq 0$, $\bm{e}_3$  and $\bm{e}_5$ are the third and fifth eigenmodes for Eq.~\eqref{Eq_eigenpb}, $\zeta_t$ is a random variable uniformly distributed in $(-1,1)$, and $f(t)$ is a square signal, whose activation is randomly distributed in the course of time. 
More precisely, given a firing rate $f_r$  in $(0,1)$, and duration $\Delta t > 0$, we define the following real-valued jump process:
\begin{equation}\label{Eq_f}
f(t) = \mathds{1}_{\{\xi_n\leq f_r\}}, \quad n \Delta t \leq t < (n+1)\Delta t, 
\end{equation}
where $\xi_n$ is a uniformly distributed random variable taking values in $[0,1]$ and $\mathds{1}_{\{\xi \leq f_r\}}=1$ if and only if $0 \leq \xi\leq f_r$.  We refer to \cite[Sec.~3.4]{debussche2013dynamics} for existence results  in Hilbert spaces  of mild solutions to SPDEs driven by jump processes, that apply to Eq.~\eqref{Eq_SPDEjump}. This type of jump process is also known as a dichotomous Markov noise \cite{bena2006dichotomous}, or a two-state Markov jump process in certain fields \cite{stechmann2011stochastic}. It is encountered in many applications \cite{horsthemke1984noise}.

The specific values for the noise parameters (frequency rate, $f_r$, time scale,  $\Delta t$, and intensity, $\sigma$) are provided in Table \ref{Table_Bimodal_regime}. The random forcing in our simulations is chosen with two key considerations:
\bi
\item {\bf Time scale ($\Delta t$)}: The time scale $\Delta t$ is large enough to allow the system to gradually relax towards its stable steady states over time.
  \item {\bf Frequency rate ($f_r$)}: The frequency rate, $f_r$, is adjusted to ensure this relaxation process can occur effectively.
\ei
This configuration ensures that if a perturbation is applied for a duration of $\Delta t$ and then removed for several subsequent $\Delta t$ intervals, the system naturally return to its closest stable steady state.

\begin{table}[h] 
\caption{Parameter values for Eq.~\eqref{Eq_SPDEjump}}
\label{Table_Bimodal_regime}
\centering
	\setlength{\tabcolsep}{7pt}
\begin{tabular}{cccccc}
\toprule\noalign{\smallskip}
    $\epsilon$ & $f_r$   & $\Delta t$  & $\sigma$ & $\lambda$  & $L$   \\ 
\noalign{\smallskip}\hline\noalign{\smallskip}
$\frac{\epsilon^\ast}{2}$ & 0.35 & $1$ & 300 & 1.32 & 2\\
\noalign{\smallskip} \bottomrule 
\end{tabular}
\end{table}

\begin{figure*}[htbp]
\centering
\includegraphics[width=1\textwidth, height=0.4\textwidth]{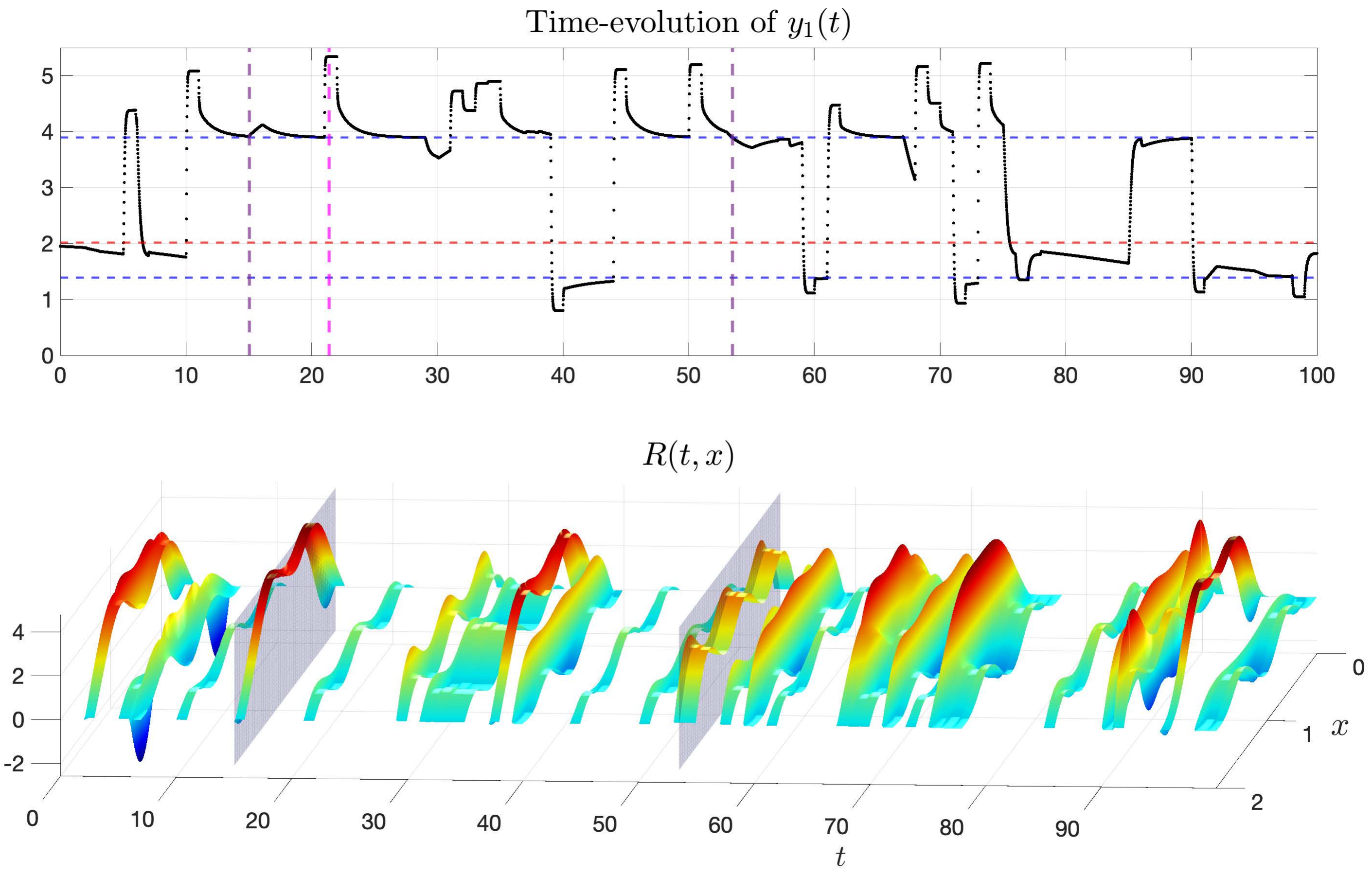}
\caption{{\bf Jump-driven dynamics.} {\bf Top panel:} The top (resp.~bottom) horizontal blue dashed line shown corresponds to $\langle U^M_\lambda,\bm{e_1}\rangle $ (resp.~$\langle U^m_\lambda,\bm{e_1}\rangle $) with $U^M_\lambda$ (resp.~$U^m_\lambda$) denoting the stable steady state to Eq.~\eqref{Eq_Sbif} with maximum (resp.~minimum) energy. The horizontal dashed red line corresponds to $\langle U^\ast_\lambda,\bm{e_1}\rangle $ with  $U^\ast_\lambda$ denoting the unstable steady state. The black dotted curve going across these horizontal dashed lines (and the associated potential barriers) is showing $y_1(t)=\langle u(t) , \bm{e}_1\rangle$ with $u$ solving Eq.~\eqref{Eq_SPDEjump} for the parameter values of Table \ref{Table_Bimodal_regime}. {\bf Bottom panel:} Here, is shown $R(t,x)=u(t,x)/\texttt{std}(u(t)) -\eta_t(x)/\texttt{std}(\eta_t)$ when $\eta_t\neq0$. The shaded planes mark two  different patterns exhibited by the system's response, for two close $y_1$-values (vertical purple dashed lines), illustrating the difficulty of closing accurately the $y_1$-variable. See Text for more details.}  \label{Fig_Jump_Dynamics}
\end{figure*}

 This relaxation behavior after a perturbation is clearly visible in the first mode's amplitude, $y_1(t)=\langle u(t) , \bm{e}_1\rangle$ (where $u$ solves Eq.~\eqref{Eq_SPDEjump}). The top panel of Figure \ref{Fig_Jump_Dynamics} illustrates this concept.  Shortly after the vertical dashed magenta line at $t = 21.4$, we observe $y_1(t)$ relaxing towards the constant value $\langle U^M_\lambda,\bm{e_1}\rangle $. Here, $U^M_\lambda$ denotes the stable steady state with the largest energy for $\lambda = 1.32$ (represented by the top blue dot in Figure \ref{Fig_Sbif}). Similar relaxation episodes occur throughout the simulation, with $y_1(t)$ sometimes relaxing towards $\langle U^m_\lambda,\bm{e_1}\rangle $, involving the stable steady state $U^m_\lambda$ with the smallest energy (refer again to Figure \ref{Fig_Sbif}).

While the noise intensity ($\sigma=300$) in Table Table \ref{Table_Bimodal_regime} may seem large, this value is necessary because we specifically target modes $\bm{e}_3$ and $\bm{e}_5$ used  in the stochastic forcing defined in Eq.~\eqref{forcing_eta}. These modes contribute less than 1\% of the total energy in the unforced parabolic equation.  When the jump noise is activated with this intensity, the amplitudes of $y_3(t)$ and $y_5(t)$ fluctuate significantly, becoming comparable in magnitude to $y_1(t)$. This eliminates the previous order-of-magnitude difference.

As a result, the random perturbations introduced by modes $\bm{e}_3$ and $\bm{e}_5$, through the system's nonlinearities, can now drive the SPDE solution across the potential barriers associated with the steady states. This explains the multimodal behavior observed in the time evolution of $y_1(t)$ as shown in Figure \ref{Fig_Jump_Dynamics}.

A crucial question remains: with such a strong random force ($\sigma=300$), are the observed dynamics in the system solely a reflection of this forcing, or is there more to it? To address this, we calculate a quantity called $R(t, x)$. This quantity compares the normalized solution, $u(t,x)/\texttt{std}(u(t))$, with the normalized forcing term, $\eta_t(x)/\texttt{std}(\eta_t)$, at times when forcing is present ($\eta_t\neq 0$); see caption of Figure \ref{Fig_Jump_Dynamics}. Here,  $\texttt{std}$ denotes standard deviation. 

The key point is that the forcing term, $\eta_t(x)$, consists of a time-dependent coefficient multiplied by a specific spatial pattern defined by $p(x)=\bm{e}_3(x) +\bm{e}_5(x)$.  If the system's response were strictly proportional to the forcing, then $R(t, x)$ would simply match this spatial pattern ($p(x)$ with three peaks and two valleys).

However, the bottom panel of Figure \ref{Fig_Jump_Dynamics} tells a different story. We see that $R(t, x)$ deviates significantly and in a complex way from $p(x)$ as time progresses. This implies that the system's response is not simply a mirror image of the forcing term. The nonlinear terms in the equation play a crucial role in shaping the response dynamics.

Developing reduced models for this system faces another significant hurdle: the system's high sensitivity to small changes, particularly in the first mode's amplitude ($y_1$). Here is the key point: even tiny variations in $y_1$, especially near the potential barrier of a steady state, can lead to vastly different spatial patterns in the full SPDE solution. Figure  \ref{Fig_Jump_Dynamics} illustrates this sensitivity.
For example, consider points $t_1 = 15.04$ and $t_2 = 53.49$ marked by the two vertical purple dashed lines in the top panel of Figure \ref{Fig_Jump_Dynamics}. The corresponding values of $y_1$ are very close ($y_1(t_1) = 3.92$ and $y_1(t_2) = 3.89$). However, the SPDE solutions at these times (shown in the vertical shaded planes) exhibit dramatically different spatial patterns.
 
The significant difference between the response of the first mode's amplitude ($y_1$) and the full SPDE solution ($u$) in equation \eqref{Eq_SPDEjump}  underscores the critical need for an accurate parameterization of the neglected variables.  Capturing this kind of sensitive behavior with a reduced model that only tracks $y_1$ requires a comprehensive representation of the influence from these excluded variables. We demonstrate in the next section that our optimal parameterization framework effectively tackles this challenging parameterization issue.

\subsection{Non-Markovian optimal reduced model}\label{Sec_OPM_Jumpcase}
We describe below how our Backward-Forward framework  introduced in this work to derive stochastic parameterizations of SPDEs driven by white noise can be readily adapted  to SPDEs driven by jump processes.  
Equation \eqref{Eq_SPDEjump}, driven by the jump process defined in Eq.~\eqref{forcing_eta}, serves as a concrete example. Our objective is to develop a reduced model for the first mode's amplitude in Eq.~\eqref{Eq_SPDEjump}. In other words, we aim to create a simplified, one-dimensional stochastic model that can accurately reproduce the complex dynamics discussed in Section \ref{Sec_jump-driven}, including the stochastic transitions.

To do so, we rewrite Eq.~\eqref{Eq_SPDEjump} for the fluctuation variable 
\be
v = u - {U^\ast_\lambda} ,
\ee
where ${U^\ast_\lambda} $ denotes the unstable steady state of Eq.~\eqref{Eq_Sbif} (for $\lambda=1.32$) corresponding to the red dot shown in Fig.~\ref{Fig_Sbif}. The fluctuation equation in the $v$-variable reads then as follows: 
\bea \label{Eq_fluc}
\partial_t v - \partial_{x}^2 v & -  \lambda (2{U^\ast_\lambda}  - 3 \epsilon {U^\ast_\lambda} ^2)v  - \lambda (1 - 3 \epsilon {U^\ast_\lambda} )v^2 \\
& + \lambda \epsilon v^3 = \lambda \eta_t, \;\mbox{on } (0,L)\times (0,T).
\eea

Each $n$th eigenmode, $\bm{e}_n$, solving Eq.~\eqref{Eq_eigenpb}  turns out to be very close to the Fourier sine mode $\sqrt{2/L}\sin(n\pi x/L)$, and the eigenvalues decay quickly.
The first few eigenvalues are $\beta_1 = 0.1815$, $\beta_2 = -7.4966$, $\beta_3 = -19.9665$, $\beta_4 = -37.2840$, and $\beta_5 = -59.5108$.  Here only the first mode $\bm{e}_1$ is unstable while the others are all stable. 

Our reduced state space is thus taken to be 
\be
H_\c=\mbox{span}\{\bm{e}_1\}.
\ee 
 For the parameter setting of Table \ref{Table_Bimodal_regime}, most of the energy ($99.9\%$) is distributed among  the unstable mode $\bm{e}_1$, on one hand,  and the forced modes, $\bm{e}_3$ and $\bm{e}_5$, on the other,  with energy carried by the third mode's  amplitude representing up to $60\%$ depending on the noise path. An accurate parameterization of these forced modes along with their nonlinear interactions with the unstable mode is thus key to achieve for a reliable reduced model of the first mode's  amplitude.

To do so, we consider the BF parameterization framework of  Section \ref{Sec_LIA_cubic_case} for cubic SPDEs driven by white noise that we adapt  to the jump noise case. If one replaces the Brownian motion in Eq.~\eqref{Eq_BF_SPDEcubic} by the jump process term, $\sigma \zeta_t f(t)$ (from Eq.~\eqref{forcing_eta}), one gets the following parameterization of the $n$th mode's amplitude:
\bea\label{eq:h1tau_jumpnoise}
\Phi_{n} (\tau_{n},&  X, t; \omega) \\
& = e^{\beta_n \tau} Y+J^{n}_{\tau}(t; \omega) + D^{n}_{11}(\tau) B_{11}^n X^2 \\
& \;\;  + E_{111}^n(\tau)  C_{111}^n  X^3,  \; n\in\{3,5\},
\eea
where the  coefficients $B_{11}^n$, $C_{111}^n$, $D^{n}_{11}$, and $E_{111}^n$ are those defined in Section \ref{Sec_LIA_cubic_case}, while the stochastic term is now given by the non-Markovian coefficient
\be \label{Eq_J_inital}
J^{n}_{\tau}(t; \omega) = \sigma  \lambda \int_{t-\tau}^t e^{\beta_n (t-s)}  \zeta_s f(s) \d s.
\ee
This stochastic term is the counterpart of the  $Z^{n}$-term defined in  \eqref{Eq_Zn_decomp} for the white noise case. 
As for the $Z^{n}$-term,  the $J^{n}$-term is efficiently generated by solving the following  ODE with random coefficients: 
\be\label{Eq_for_J}
\frac{\d J}{\d t} =  \beta_n J +  \sigma \lambda \left( \zeta_t f(t) - e^{\beta_n \tau} \zeta_{t-\tau} f(t-\tau)\right),
\ee
for $n=3$ and $n=5$; cf.~Eq.~\eqref{Eq_for_I}. 
Note that the parameterization Eq.~\eqref{eq:h1tau_jumpnoise} takes the same structural form as Eq.~\eqref{eq:h1_tau_Chafee}  in the case of white noise forcing (see Remark \ref{Rmk_OU}), with the bounded jump noise replacing the latter here.  
This is due to the underlying backward-forward systems at the core of the derivation of these formulas where only the nature of the forcing changes; see also Sec.~\ref{Sec_chart} below.

Then, after optimization over a single training path of the  parameterization defects $Q_n(\tau)$ (see Section \ref{Sec_opti_jumps} below),  an analogue to the optimal reduced model (Eq.~\eqref{Eq_reduced_chafee}) becomes in our case the following 1-D  ODE with path-dependent coefficients:  
\bea \label{Eq_reduced_S_shaped}
&\dot{y}  = \beta_1 y  + \Big\langle  (\lambda - 3\lambda \epsilon {U^\ast_\lambda} )   \big(y (t) \bm{e}_{1} +\P(t)\big)^2, \boldsymbol{e}_1\Big\rangle\\
& \qquad  \qquad \qquad  - \lambda \epsilon  \Big\langle (y (t) \bm{e}_{1} +\P(t))^3, \boldsymbol{e}_1\Big\rangle,\\
&\P(t) = \Phi_{3}(\tau^*_{3}, y(t), t) \bm{e}_{3}+ \Phi_{5}(\tau^*_{5},  y(t), t) \bm{e}_{5}.
\eea
Here, $\P(t)$ is the optimized stochastic parameterization in the class $\mathcal{P}$ (see Section \ref{Sec_classP})
for the SPDE \eqref{Eq_fluc} driven by the jump process $\eta_t$ defined in Eq.~\eqref{forcing_eta}. Note that the required non-Markovian $J^n$-terms to simulate Eq.~\eqref{Eq_reduced_S_shaped}  are obtained by solving  Eq.~\eqref{Eq_for_J} in the course of time. So a total of three ODEs with random coefficients are solved to obtain a reduced model of the first mode's amplitude. Thus, Eq.~\eqref{Eq_reduced_S_shaped} along with its auxiliary $J$-equations, form the optimal reduced system.

 The inner product in  Eq.~\eqref{Eq_reduced_S_shaped}  leverages the analytical expression of $\P(t)$ provided in the same equation. This inner product is computed once offline, resulting in a degree-9 polynomial with random coefficients after expanding the  $\Phi_3$ and $\Phi_5$ terms as described in Eq.~\eqref{eq:h1tau_jumpnoise}. However, this expression can be further simplified. By analyzing the order of magnitude of the terms in  Eq.~\eqref{eq:h1tau_jumpnoise} for the parameter regime of Table  \ref{Table_Bimodal_regime} ($\sigma$ large), we observe that $\P(t)$ simplifies to
\be\label{OU_approx}
\P(t) \approx J^{3}_{\tau^*_3} (t; \omega)\bm{e}_3 +J^{5}_{\tau^*_5}(t; \omega)\bm{e}_5,
\ee
setting $Y=0$ in $\Phi_n$ ($n=3,5$) defined in Eq.~\eqref{eq:h1tau_jumpnoise}\footnote{Note that the other unresolved (and unforced) modes contain less than $0.1\%$ of the total energy and are thus not taken here in account into the parameterization $\P(t)$.}.
This approximation expresses essentially a prominence of  the stochastic terms over the nonlinear ones in the parameterization  \eqref{eq:h1tau_jumpnoise}, unlike for the reduction of the sACE analyzed in Section \ref{Sec_stat_reduced}. The parameterization becomes then purely stochastic without nonlinear terms. Yet,  the optimal reduced system exhibits nonlinear products of stochastic coefficients.

To make this point, we introduce a few more notations.
We denote by $a_\lambda^\epsilon(x)$ the spatial coefficient, $\lambda - 3\lambda \epsilon {U^\ast_\lambda}(x)$, involved in the quadratic term of Eq.~\eqref{Eq_reduced_S_shaped} (as inherited from Eq.~\eqref{Eq_fluc}).
To simplify the notations we denote also  by $J_{3}^\ast$ (resp.~$J_{5}^\ast$)  the term $J^{3}_{\tau^*_3}(t; \omega)$ (resp.~$J^5_{\tau_5^*}(t; \omega)$), and introduce the bilinear operator $G_2(u,v)(x)=u(x)v(x),$ for $x$ in $(0,L)$ and any functions $u,v$ in $L^2(0,L)$. We have then, using the approximation formula \eqref{OU_approx}, that the quadratic term in Eq.~\eqref{Eq_reduced_S_shaped}
expands as follows 
\begin{widetext}
\bea\label{Expansion_quadterm}
 &\Big\langle  a_\lambda^\epsilon  \big(y (t) \bm{e}_{1} +\P(t)\big)^2, \boldsymbol{e}_1\Big\rangle=y^2(t) \langle a_\lambda^\epsilon  G_2(\bm{e}_1,\bm{e}_1),\bm{e}_1\rangle+ (J_{3}^\ast(t))^2\langle a_\lambda^\epsilon G_2(\bm{e}_3,\bm{e}_3),\bm{e}_1\rangle+(J_{5}^\ast(t))^2 \langle a_\lambda^\epsilon G_2(\bm{e}_5,\bm{e}_5),\bm{e}_1\rangle\\
 &\hspace{2cm}+2 \langle a_\lambda^\epsilon G_2(\bm{e}_3,\bm{e}_5),\bm{e}_1\rangle  J_{3}^\ast(t) J_{5}^\ast(t) +2 \langle a_\lambda^\epsilon G_2(\bm{e}_1,\bm{e}_3),\bm{e}_1\rangle y(t) J_{3}^\ast (t) +2 \langle a_\lambda^\epsilon G_2(\bm{e}_1,\bm{e}_5),\bm{e}_1\rangle  y(t) J_{5}^\ast (t),
\eea
\end{widetext}
 with a similar expression obtained from the cubic terms  in Eq.~\eqref{Eq_reduced_S_shaped}. 
Already from \eqref{Expansion_quadterm} one observes the production of stochastic constant terms driven by $(J_{3}^\ast(t))^2$,   $(J_{5}^\ast(t))^2 $ and $J_{3}^\ast(t) J_{5}^\ast(t)$ as well as linear terms driven by  $J_{3}^\ast(t)$  and $J_{5}^\ast(t)$ contributing to the dynamics of the surrogate first mode's amplitude $y$. These stochastic terms turn out to be essential for a good emulation of the first mode's amplitude. 
This point is illustrated by the following experiment.

\begin{figure}
\centering
\includegraphics[width=.5\textwidth, height=0.2\textwidth]{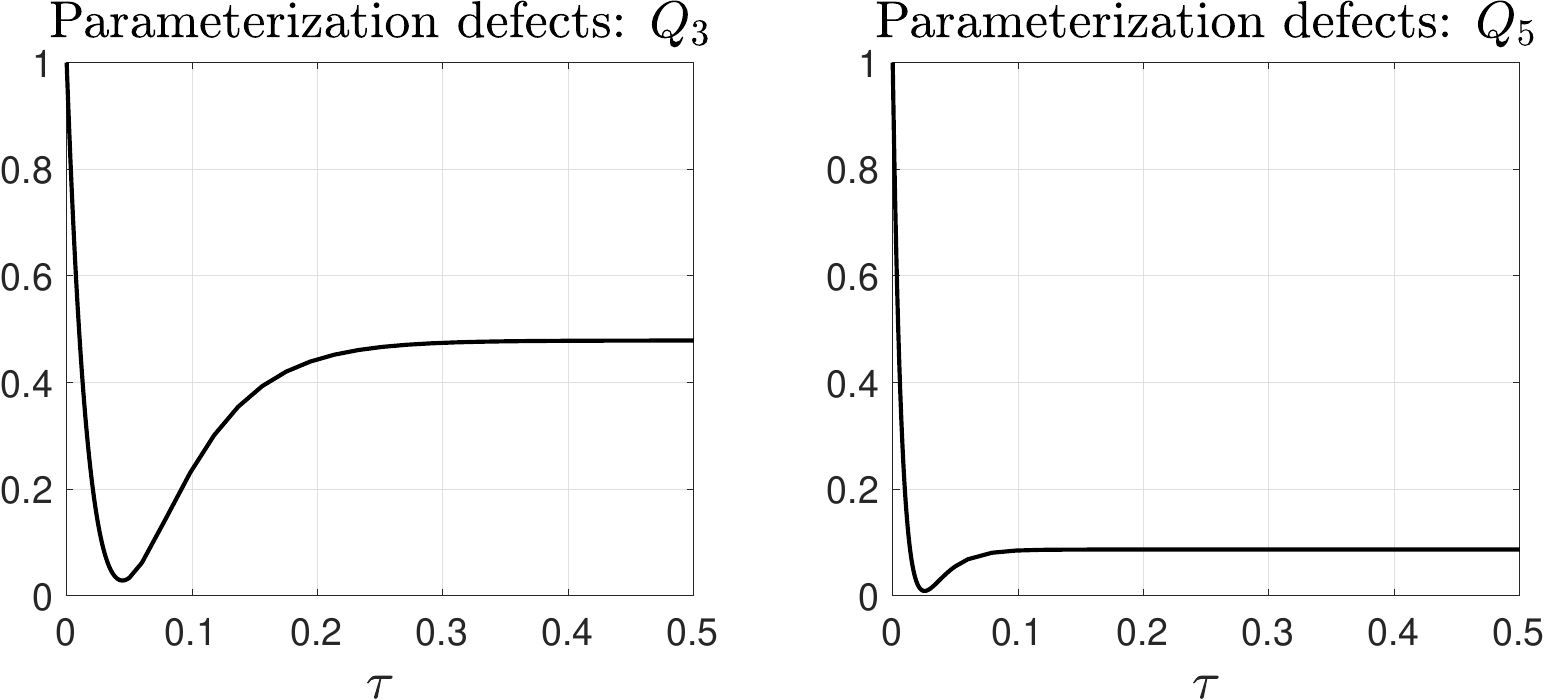}
\caption{{\bf The normalized parameterization defects, $Q_n$ for $n=3$ and $n=5$.} The optimization is performed on a single path. See Text.} 
\label{Fig_S_eqn_Qtau}
\end{figure}

\subsection{Training data, parameterization defects, and reduced model's skills}\label{Sec_opti_jumps}
The training data for the optimal reduced model is comprised of $N = 40,000$ snapshots. These sACE snapshots are generated over a training interval from $0$ to $T$, with $T$ = 400, for a single, noise path driving the sACE. The corresponding normalized parameterization defects (quantified by $Q_n$) are visualized in Figure \ref{Fig_S_eqn_Qtau}.

Here, it's important to note that the distinct minima observed in $Q_3$ and $Q_5$ are not caused by weak time-scale separations (as discussed in Section \ref{Sec_weak_timescale} for the sACE). Instead, they are primarily a consequence of the high noise intensity ($\sigma$) used in this case, which leads to solutions with large amplitudes. This is further supported by the observation that these minima become less pronounced as $\sigma$ is reduced (verified for arbitrary  noise paths).

\begin{figure}
\centering
\includegraphics[width=.5\textwidth, height=0.3\textwidth]{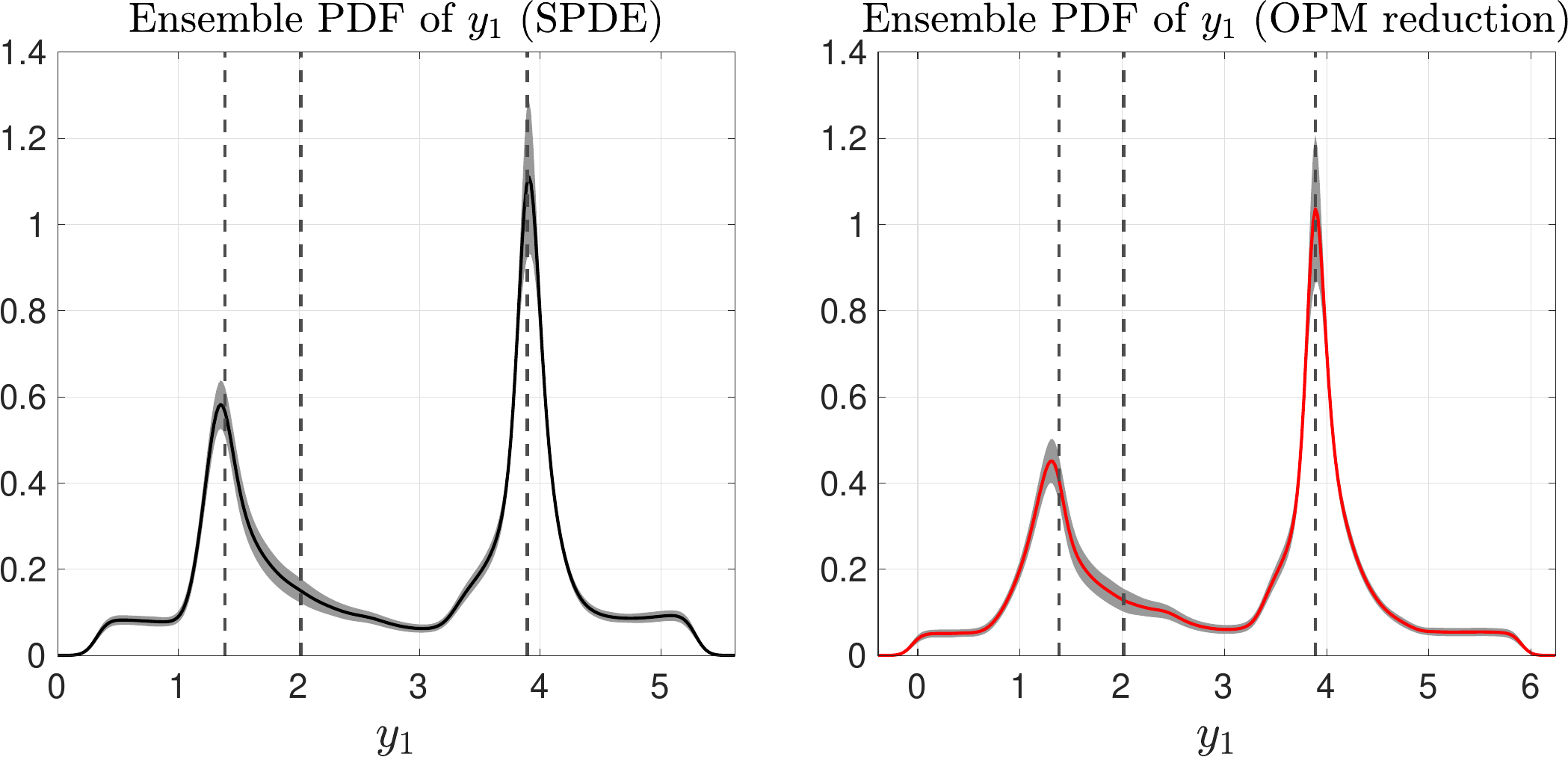}
\caption{{\bf Ensemble PDF}.  Here, $10^5$ out-of-sample test paths are used to estimate these ensemble PDFs. Each underlying solution path is made of  $N = 2\times 10^5$ iterations for a time-step of $\delta t= 10^{-2}$.} \label{Fig_EnsPDF}
\end{figure}

\begin{figure}
\centering
\includegraphics[width=.5\textwidth, height=0.3\textwidth]{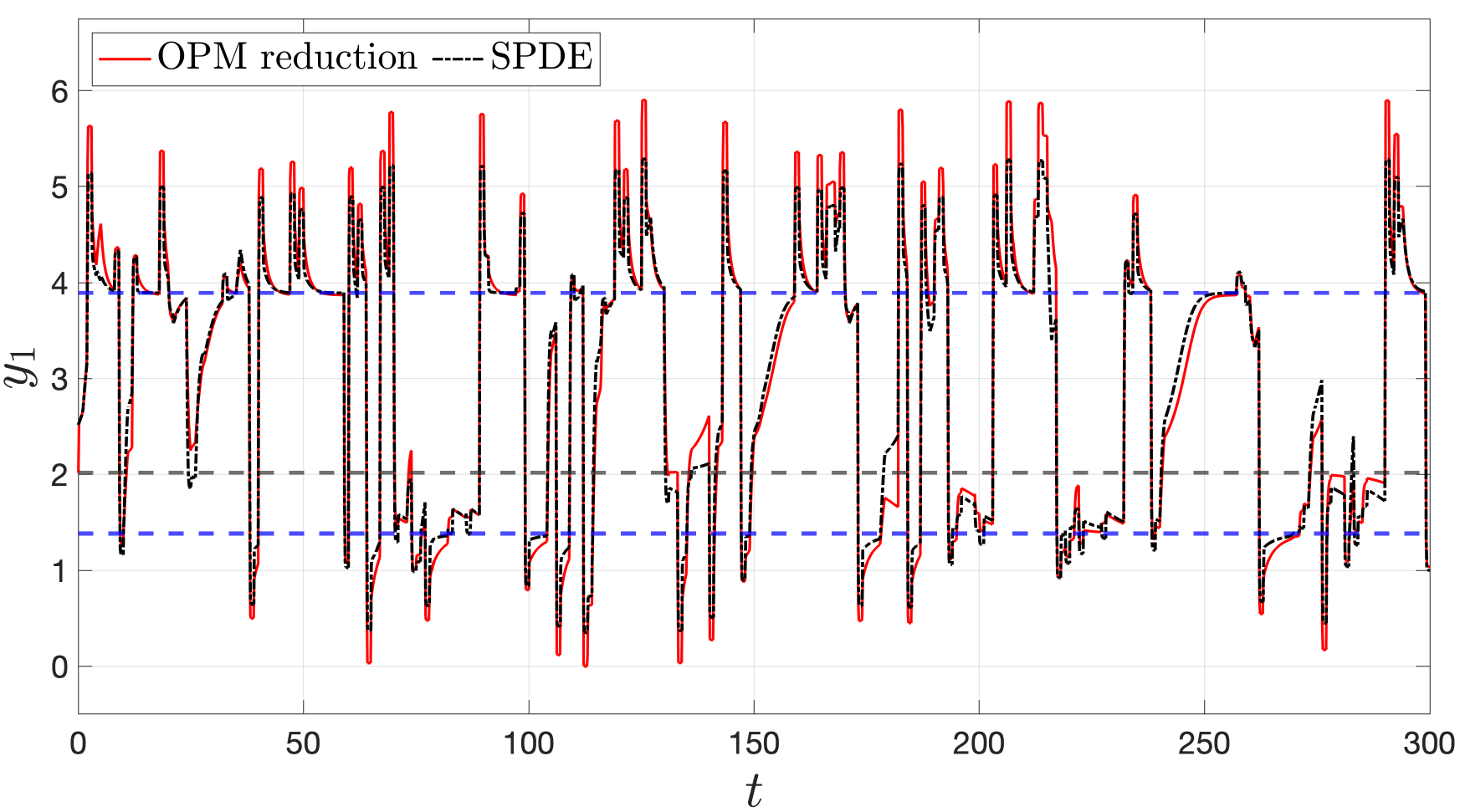}
\caption{{\bf First's mode amplitude dynamics: Full system vs optimal reduction}.} \label{Fig_y1timeseries}
\end{figure}

Figure  \ref{Fig_EnsPDF} shows the optimal parameterization's ability to reproduce first mode's complex behavior.
This figure showcases the effectiveness of the optimal reduced model (Eq.~\eqref{Eq_reduced_S_shaped} with $\P(t)$ from Eq.~\eqref{OU_approx}) in capturing the bimodal nature of the first mode's amplitude. The model's performance is evaluated using a vast set of unseen test paths. When compared to the actual SPDE solution's amplitude for the first mode (across these test paths), the optimal reduced model closely follows the various transitions between metastable states, crossing as needed the potential barriers. However, it is worth noting that this reduced model tends to slightly overestimate the intensity of these larger excursions (as shown in Figure \ref{Fig_y1timeseries}).

\section{Towards the Reduction of More General L\'evy-driven Dynamics} \label{Sec_Levy}

We prolong the discussion and results of Section \ref{Sec_Jump} towards more general L\'evy noise. 
For the sake of simplicity we focus on finite-dimensional SDEs driven by L\'evy processes.

There is a vast literature on L\'evy processes which roughly speaking are processes given as the sum of a Brownian motion and a jump process \cite{tankov2003financial,peszat2007stochastic,applebaum2009levy,kuhn2017levy,oksendal2019stochastic}. Unlike the case of diffusion processes (SDEs driven by Brownian motions), the representation of Kolmogorov operators is non-unique in the case of L\'evy-driven SDEs, and may take the form of operators involving fractional Laplacians or singular integrals among other representations \cite{Kwasnicki:2017aa}. We adopt here the definition commonly used in probability theory \cite[Theorem 3.5.14]{Applebaum_2019} which presents the interest of being particularly intuitive as recalled below. 
We consider SDEs of the form
\be\label{Eq_Levy_additive}
\d X_t = {\bm F} (X_t) \d t +{\bm \Sigma}(X_t) \d \W + \d  \bm{L}_t,
\ee  
where $\W$ is a Brownian motion in $\mathbb{R}^d$, and $\bm{L}_t$ is a L\'evy process on $\mathbb{R}^d$ independent from $\W$. Here, we take ${\bm F}$ to be a smooth vector field on $\mathbb{R}^d$, and $ {\bm \Sigma}(\x)$ to be, for each $\x$, a $d\times d$ matrix with smooth entries, such that $(a_{ij}(\x))= {\bm \Sigma}(\x){\bm \Sigma}^T(\x)$ is a positive definite matrix for each $\x$ in $\mathbb{R}^d$.

Roughly speaking, a L\'evy process $\bm{L}_t$ on $\mathbb{R}^d$ is a non-Gaussian stochastic process  with independent and stationary increments that experience sudden, large jumps in any direction. 
The probability distribution of the these jumps is characterized by a non-negative Borel measure $\nu$ defined on $\mathbb{R}^d$ and concentrated on  
$\mathbb{R}^d\backslash \{0\}$ that satisfies the property  $\int_{\mathbb{R}^d\backslash\{0\}} \min(1,\y^2)\nu (\d \y) <\infty$.  This measure $\nu$ is called the jump measure of the L\'evy process $\bm{L}_t$. 
Sometimes $X_t$ itself is referred to as a L\'evy process with triplet $({\bm F},{\bm \Sigma},\nu)$. Within this convention, we reserve ourselves the terminology of a L\'evy process to a process with  triplet $(0,0,\nu)$.
We refer to \cite{applebaum2009levy} and \cite{peszat2007stochastic} for the mathematical background on L\'evy processes.

Under suitable assumptions on ${\bm F}$, $ {\bm \Sigma}$, and the L\'evy measure $\nu$, the solution $X_t$ to Eq.~\eqref{Eq_Levy_additive}
is a Markov process (e.g.~\cite{protter2005stochastic}) and even a Feller process \cite{kuhn2018solutions}, with associated Kolmogorov operator 
taking the following integro-differential form for (e.g.)~$\psi$ in $C^{\infty}_c(\mathbb{R}^d)$ (Courr\`ege theorem \cite{courrege1965forme,kuhn2017levy}):
\begin{widetext}
\bea\label{Eq_Kolmo2_general}
\mathcal{L}_K \psi (\x) &={\bm F} (\x) \cdot \nabla \psi+ \sum_{i,j=1}^d a_{ij}(\x) \partial_{ij} \psi + \Gamma \psi(\x), \mbox{with}\\
&\Gamma \psi (\x)=\int_{\mathbb{R}^n\backslash\{0\}} \Big[\psi(\x+\y)-\psi(\x)- \y\cdot \nabla \psi (\x) \mathds{1}_{\{\norm{\y}<1\}} \Big] \nu(\d \y),
\eea
\end{widetext}
where $ \mathds{1}_{\{\norm{\y}<1\}}$ denotes the indicator function of the (open) unit ball in $\mathbb{R}^n$; see also \cite[Theorem 3.3.3]{applebaum2009levy}.

The first-order term in Eq.~\eqref{Eq_Kolmo2_general} is the drift term caused by the deterministic, nonlinear dynamics. 
The second-order differential operator  represents the diffusion part of the process $X_t$. It is responsible for the continuous component of the process. 

The $\Gamma$-term, involving the integral, represents the jump part of the process. 
 It captures the discontinuous jumps that the process experience due to the sudden changes caused by the L\'evy process $\bm{L}_t$.
Its intuitive interpretation breaks down as follows. 
The term, $\psi(\x+\y)-\psi(\x)$, calculates the difference in the test function value before and after the potential jump, capturing the change in the test function due to the jump.

The term $-\y\cdot \nabla \psi (\x) \mathds{1}_{\{\norm{y}<1\}}$ represents a first-order correction for small jumps. It aims to account for the fact that a small jump might not land exactly on the grid point ($\x + \y$),
but somewhere in its vicinity. This term is often referred to as the Girsanov correction. 
Thus, the integral term $\Gamma$ accounts for all possible jump sizes ($\y$) within a unit ball, as 
weighted by the L\'evy measure $\nu(\d \y)$. The notion of Kolmogorov operators can be extended to SPDEs driven by cylindrical Wiener processes, with $\x$ lying in a Hilbert space; see \cite{da2004kolmogorov}. Building up on such a functional framework, it is possible to provide a rigorous meaning of the Kolmogorov operator $\mathcal{L}_K $ in the case of an SPDE driven by a general L\'evy process, forcing possibly infinitely many modes.

In the case of SPDE  driven by a scalar jump process (forcing only two modes) studied in Section \ref{Sec_Jump} (see Eq.~\eqref{forcing_eta}), the jump measure is simple to describe.  It corresponds to  the measure  $\nu(\d s)$ ($s$ real) associated with the real-valued jump process $f(t)$ given by Eq.~\eqref{Eq_f}.
This measure  is the Dirac delta $\lambda \delta_{s=1}$, where $\lambda$ is the intensity of the Poisson process associated with the random comb signal. Recall that this  intensity is the limit of the probability of a single event occurring in a small interval divided by the length of that interval as the interval becomes infinitesimally small. Therefore, $\lambda$ is equal to the firing rate $f_r$ of $f(t)$.

In the case of SDEs driven by $\alpha$-stable  L\'evy processes, 
\be
\nu(\d \y)=c_{d,\alpha} \|\y\|^{-d-\alpha} \d \y,
\ee
where  the non-Gaussianity index   $\alpha$ lies  in $(0,2)$, and $c_{d,\alpha}$ is a constant that depends on the dimension and involves the gamma function of Euler \cite{Kwasnicki:2017aa}. A $2$-stable ($\alpha= 2$) process is simply a Brownian motion. 

\begin{figure*}[htbp]
\centering
\includegraphics[width=.9\textwidth, height=0.6\textwidth]{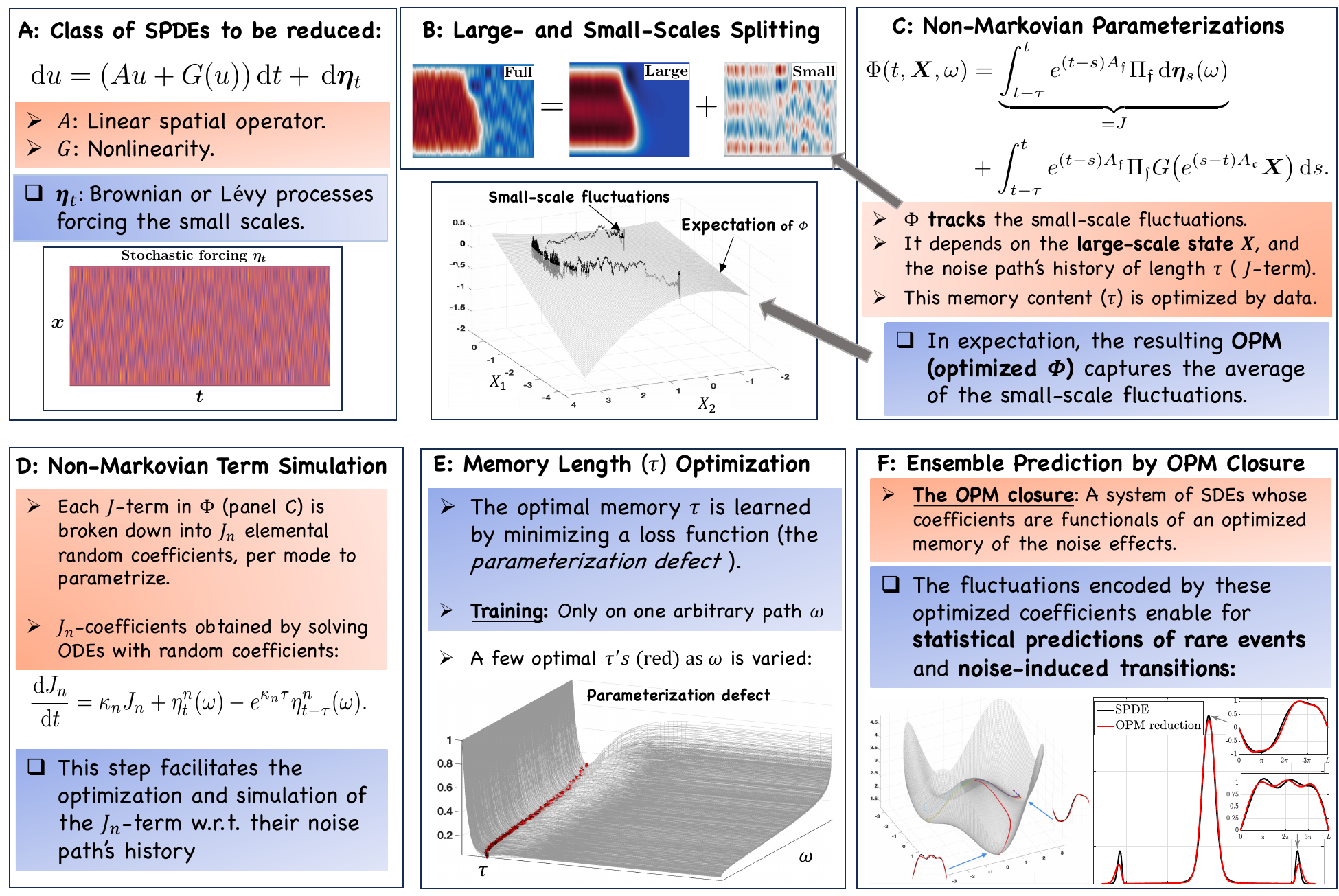}
\caption{{\bf Chart presenting the main  steps, objectives, and characteristics of the OPM reduction approach.}} 
\label{Fig_chart}
\end{figure*}

The rigorous theory of slow and center manifolds are in their infancy for SDEs driven by   $\alpha$-stable  L\'evy processes \cite{yuan2019slow,yuan2021stochastic}, but the main formulas share the same ingredients  as recalled in Section \ref{sec_inv_eq}, with analogous backward-forward interpretations of those discussed in Section \ref{Sec_fully_coupled_BF_gen} for  the Lyapunov-Perron integrals involved therein (\cite[Sec.~4]{yuan2019slow}). 
The non-Markovian parameterizations of Section \ref{Sec_OPMwhite} and their data-driven optimization, can thus be extended, at least formally, for the case of SDEs driven by   $\alpha$-stable  L\'evy processes and will be pursued elsewhere on concrete applications. 


\section{A Chart Summarizing the Approach}\label{Sec_chart}
 To streamline the application of the OPM reduction approach for stochastically forced systems presented in this paper, we provide a concise overview of its key steps, objectives, and characteristics in the chart shown in Figure Fig.~\ref{Fig_chart}. This summary is tailored to the class of SPDEs discussed in Section \ref{Sec_SPDEs}, with the primary goal of performing a reduction with a cutoff scale just below the forcing scale. This reduction involves a solution splitting where only the disregarded "small" scales are stochastically forced (Figures \ref{Fig_chart}A and \ref{Fig_chart}B).

Given this splitting, where the disregarded scales are those directly forced by stochasticity, our approach reveals the importance of non-Markovian terms in the parameterization of these scales. These terms depend on the noise's path history, enabling the tracking of small-scale fluctuations. The accuracy of this tracking is primarily governed by the memory content  ($\tau$) of the parameterization, which is optimized from data by minimizing a natural loss function (Figure \ref{Fig_chart}E) over a single noise realization $\omega$.

Once optimized, the resulting non-Markovian parameterization, the Optimal Parameterizing Manifold (OPM), allows us to capture the average behavior of small-scale fluctuations as a functional of large scales by averaging over noise paths (see inset below Figure \ref{Fig_chart}B). To efficiently optimize the memory content ($\tau$) and avoid the computation of cumbersome quadratures,   we decompose the $J$-term in Figure \ref{Fig_chart}, $J(t,\x)=\sum  J_n(t) {\bm e}_n (\x)$, into random coefficients $J_n$ that solve simple auxiliary ODEs dependent on $\tau$ and the noise path $\omega$. The resulting OPM closure, exemplified by Equation \eqref{Eq_closure}, takes the form of an SDE system with coefficients that are nonlinear functionals of these $J_n$ terms (Figure \ref{Fig_chart}F) reflecting the nonlinear interactions between the noise and nonlinear terms in the original system.

As demonstrated with SPDEs driven by Gaussian noise (Section \ref{Sec_ACE}) or jump processes (Section \ref{Sec_Jump}), the OPM closure effectively predicts noise-induced transitions and rare event statistics. This success is attributed to our hybrid approach, which combines analytical insights for constructing the parameterizing manifold $\Phi$ (Figure \ref{Fig_chart}C and Section \ref{Sec_LIA_cubic_case}) with data-driven learning of the optimal memory content $\tau$ (Section III C).

\section{Discussion and Perspectives}\label{Sec_discussion}

Thus, this work demonstrates that Optimal Parameterizing Manifolds (OPMs) constructed from Backward-Forward (BF) systems and optimized using a single solution path can accurately reproduce noise-induced transitions across a broad range of  out-of-sample, test paths. This approach is effective for non-Gaussian  noise with jumps. A remarkable feature of the approach is that, even though these two different types of noise typically lead to different dynamical responses, the OPMs in these two settings take exactly the same structural form with the corresponding non-Markovian terms depending on the past history of the corresponding noise forcing the original SPDE;  cf.~Eqns.~\eqref{eq:h1tau_stoch}--\eqref{Eq_Zn_decomp} and cf Eqns.~\eqref{eq:h1tau_jumpnoise}--\eqref{Eq_J_inital}.

 Interestingly, the optimal parameterization formulas align with the approximation of stochastic invariant manifolds when they exist. However, when these manifolds are no longer applicable, the optimal parameterizations transcends the invariance constraint. In such cases, it provides the optimal stochastic parameterization with respect to the random invariant measure $\rho_\omega$ (defined in Section \ref{Sec_OPM_IM}), conditioned on the resolved variable $X$ (as proven in Theorem \ref{Thm_variational-pb2}).

The formulas derived from the BF systems in Section \ref{Sec_OPMwhite} enable practical approximation of the theoretical optimal parameterization, particularly its non-Markovian random coefficients. These coefficients become crucial, especially when noise affects unresolved modes. They essentially encode how noise propagates to the resolved modes through nonlinear interactions.

 We have shown, through various examples, that training on a single noise path 
equips the coefficients with an optimized historical memory (of the noise) that plays a key role in the accurate prediction of noise-induced excursion statistics, even when applied to a diverse set of test paths.  This demonstrates the strength of our hybrid framework, which combines analytical understanding with data-driven insights. This approach allows us to overcome the ``extrapolation problem," a significant challenge for purely data-driven machine-learned parameterizations.

 Section  \ref{Sec_fully_coupled_BF_gen} discusses the possibility of constructing more intricate parameterizations using an iterative scheme based on more general BF systems. However, these elaborations require more computational effort due to the repeated stochastic convolutions. Implementing them efficiently necessitates deriving auxiliary ODEs (analogues to Eq.~\eqref{Eq_for_I}) for the coefficients ($a_n$ and $b_{in}$) in  Eqns.~\eqref{Eq_an} and \eqref{Eq_bn}. Reference \cite{CLW15_vol2}, Section 7.3, provides examples of such auxiliary ODEs for handling repeated stochastic convolutions.  

 These more general stochastic parameterizations have the potential to significantly enhance the parameterization capabilities of $\Phi$ (Eq.~\eqref{Eq_us1_for_LIA}). This is because they encode higher-order interaction terms with respect to the noise amplitude, as shown in Eq.~\eqref{Eq_q2_explicit}. Resolving these interactions is potentially  important for improving the prediction of rare event statistics by non-Markovian reduced systems, especially when the noise intensity is low. This opens doors to potential applications using large deviation or rare event algorithms (e.g., \cite{E_al04,cerou_2007,Vanden_West08,grafke2015instanton,Rolland_al16,Ragone2018,galfi2021applications,simonnet2021multistability}). Additionally, iterating the BF system (Eq.~\eqref{Eq_BF_fully_coupled_additive2}) beyond $\ell= 2$ can also introduce higher-order terms involving the low-mode amplitude $X$.

The inclusion of higher-term in $X$ is non-unique though. For instance, inspired by \cite[Sec.~7.2.1]{CLW15_vol2}, stochastic parameterizations  accounting for higher-order terms in $X$ are produced by solving analytically the following BF system:
\bea\label{Eq_h2}
& \mathrm{d} p^{(1)} =  A_\c p^{(1)} \d s,  \\
&\mathrm{d} p^{(2)} =  \Big( A_\c p^{(2)} + \Pi_c G_2\big(p^{(1)}\big)  \Big)\d s, \\%
& \mathrm{d} q_n  = \Big( \lambda_n q_n  +  \Pi_{n} G_2\big(p^{(2)}\big)  \Big) \d s + \sigma_n \d W_{s}^n, \\
&  p^{(1)}(t)= p^{(2)}(t) = X \in H_{\c}, \; q_{n}(t-\tau)=Y,
\eea
for which the two first equations for $p^{(1)}$ and $p^{(2)}$ are simultaneously integrated backward over $[t-\tau,t]$, followed by the forward integration of the $q_n$-equation over the same interval.
Compared to Eq.~\eqref{Eq_BF_SPDEcubic},  this BF system has two-layer in the low-mode variable $p$ resulting into parameterizations of order 4 in $X$; see \cite[Theorem 2]{CL15} for an example. These higher-order terms come from the nonlinear self-interactions between the low
modes as brought by the  term $\Pi_c G_2\big(p^{(1)}\big)$. Accounting for additional interactions either through \eqref{Eq_h2} or the iterative procedure of Section \ref{Sec_fully_coupled_BF_gen} can help improve the parameterization skills (see also \cite{BM13,Gao24}),  and should be considered if accuracy is key to resolve.  

An alternative route to enhance accuracy would involve to parameterize stochastically the parameterization residual from the loss functions in \eqref{Eq_minQnHn}, post-minimization. 
This technique has proven effective in modeling the fast motion transversal to OPM manifolds by means of networks of nonlinear stochastic oscillators \cite{Chekroun2021c,LC2023}, as diagnosed by Ruelle-Pollicott resonances \cite{Chekroun_al_RP2,Tantet_al_Hopf}. Similarly, rapid dynamical fluctuations, transversal to (Markovian) OPM manifolds (Figure  \ref{Fig_toto_visual}), could benefit from refined parameterizations through the incorporation of additional stochastic components following \cite{Chekroun2021c}, potentially leading to more accurate representation of rare event statistics by the corresponding closures (Figure \ref{Fig_min_max}). 

Extending beyond SPDEs driven by Markov processes, our framework can naturally accommodate SPDEs driven by non-Markovian stochastic processes. A prime example is when the stochastic forcing is a fractional Brownian motion (fBm) \cite{mandelbrot1968fractional}, for which the existence of unstable manifolds has been demonstrated \cite{garrido2010unstable}. By leveraging Lyapunov-Perron techniques, as employed in the construction of these manifolds (Section \ref{sec_inv_eq}), our extension program to parameterization can proceed analogously to Section \ref{Sec_OPMwhite}.
This involves here as well, finite-time integration of relevant backward-forward systems to derive optimizable parameterizations from the governing equation. The random coefficients in these parameterizations would then become again functionals of the past history of the forcing noise, albeit this time incorporating a higher degree of memory compared to the stochastic process $Z^n_\tau$ used in the parameterizations Section \ref{Sec_OPMwhite}, due to the original memory content of the original forcing.  Such generalized parameterizations should provide valuable insights into the interplay between nonlinearity and noise in physical systems exhibiting long-term memory effects \cite{fraedrich2009continuum,rypdal2014long,pereira20191}.

Regardless of the Markovian or non-Markovian nature of the stochastic forcing, our reduction approach and its potential extensions offer valuable tools for addressing critical questions. By elucidating the intricate interplay between noise and nonlinearity, we can gain deeper dynamical insights. For example, we anticipate applying this reduction program to investigate early warning indicators of tipping phenomena, from a dynamical perspective, by tracing the energy flow from forced scales to physically relevant modes \cite{van2024physics,Lucarini_Chekroun_PRL24}. Climate tipping elements provide a compelling domain for exploring these concepts and applying our methods  \cite{lenton2008tipping,lenton2019climate,mckay2022exceeding}.  

The proposed stochastic parameterization method offers a new framework  for understanding noise-nonlinearity interactions in stochastic multiscale systems. It is important to note that alternative techniques based on stochastic time scale separation, such as the stochastic path perturbation approach, have been employed to study related phenomena, including the random escape of stochastic fields from criticality in explosive systems \cite{caceres1999stochastic} and first passage time problems in nonlocal Fisher equations \cite{fuentes2003nonlocal,caceres2015first}. 
Future research directions could explore the potential of our method to obtain analytical approximations for first passage times in nonlocal SPDEs near criticality. Such an extension could provide valuable insights into pattern formations for a wide range of physical and biological problems such as in infectious disease modeling \cite{ruan2007spatial}, nonlinear optics \cite{fernandez2013strong}, or vegetation models \cite{gilad2004ecosystem}.

\begin{acknowledgements}
Several insights discussed in this work have substantially benefited from the reviewers’ constructive comments, and we express our deep gratitude to them.
This work has been supported by the Office of Naval Research (ONR) Multidisciplinary University Research Initiative (MURI) grant N00014-20-1-2023, and by the National Science Foundation grants DMS-2108856, DMS-2407483, and DMS-2407484. 
This work has been also supported by the European Research Council (ERC) under the European Union’s Horizon 2020 research and innovation program (grant agreement No. 810370) and by a Ben May Center grant for theoretical and/or computational research. 
\end{acknowledgements}

\appendix

{\small 
\section{Proof of Theorem \ref{Lemma_BF}}\label{Sec_Proof}

\bp
{\bf Step 1}.
Let us remark that for any $n$
\beas\label{Eq_pin}
&\Pi_{n}  G_k(e^{(s-t)A_{\c}}X + z_{\c})=\\
&\sum_{(j_1,\cdots, j_k)\in I^k} \prod_{\ell=1}^{k} \bigg(\exp((s-t) \lambda_{j_{\ell}}) X_{j_{\ell}}+z_{j_{\ell}}\bigg) G_{j_1 \cdots j_k}^n,
\eeas
where $I=\{1,\cdots,m_c\}$ and $G_{j_1 \cdots j_k}^n= \langle G_k(\boldsymbol{e}_{j_1},\cdots,\boldsymbol{e}_{j_{k}}),\boldsymbol{e}_n^{\ast}\rangle$.

Note that the product term above can be re-written as
\bea\label{Eq_pin2}
\prod_{\ell=1}^{k} & \bigg( \exp((s-t) \lambda_{j_{\ell}}) X_{j_{\ell}}+z_{j_{\ell}}\bigg)=\\
&\sum_{K_z} \bigg(\prod_{p\in K\setminus K_z} \exp((s-t) \lambda_{j_p}) X_{j_{p}}\bigg)\bigg( \underset{q \in K_z}{\prod} z_{j_{q}}\bigg),
\eea 
where $K=\{1,\cdots,k\}$ and $K_z$ runs through the subsets (possibly empty) of disjoint elements of $K$.

Thus,
\begin{widetext}
\bea
\mathfrak{J}(t,X,\omega) & =\sum_{n\geq m_c+1} \bigg( \int_{-\infty}^t e^{(t-s)\lambda_n}\Pi_{n} G_k(e^{(s-t)A_{\c}}X + z_{\c}(s,\omega))\d s\bigg) \boldsymbol{e}_n\\
 &=\sum_{n\geq m_c+1} \sum_{(j_1,\cdots, j_k)\in I^k}  \bigg(\int_{-\infty}^t  \prod_{\ell=1}^{k} \bigg( \exp((s-t) \lambda_{j_{\ell}}) X_{j_{\ell}}+z_{j_{\ell}} \bigg) \exp((t-s) \lambda_n) \d s\bigg) G_{j_1 \cdots j_k}^n  \boldsymbol{e}_n,
\eea
\end{widetext}
and by using \eqref{Eq_pin2},  $\mathfrak{J}(t,X,\omega)$ exists if and only if for each $(j_1,\cdots, j_k)\in I^k$ with $G_{j_1 \cdots j_k}^n \neq 0$, the integrals
\be \label{Eq_I_Kz}
I_{K_z}(t)=\int_{-\infty}^t \hspace{-.1cm}\bigg( \underset{q \in K_z}{\prod} z_{q}(s,\omega) \bigg) \exp\bigg((t-s)\Big[\lambda_n-\hspace{-.25cm}\sum_{p \in K \setminus K_z} \hspace{-0.2cm}\lambda_{j_{p}}\Big]\bigg) \d s,
\ee
exist as $K_z$ runs through the subsets (possibly empty) of disjoint elements of $K=\{1,\cdots,k\}$. 

Since the OU process $z_{q}(s,\omega)$ has sublinear growth (see e.g.~\cite[Lemma 3.1]{CLW15_vol1}),  that is $\lim_{s \rightarrow \pm \infty} \frac{z_{q}(s,\omega)}{s} = 0$ for all $\omega$, one can readily check that the integral $I_{K_z}$ defined in \eqref{Eq_I_Kz} is finite for all $t$ provided that 
\be 
\Re\left(\lambda_n-\sum_{p \in K \setminus K_z} \lambda_{j_p}\right) < 0.
\ee
Thus, the Lyapunov-Perron integral $\mathfrak{J}$ defined in Eq.~\eqref{Eq_phi_stationary} is well-defined if the non-resonance condition \eqref{Eq_NR} holds. 

\medskip
{\bf Step 2.}
Thus $\mathfrak{J}$ is well defined under the condition \eqref{Eq_NR}. By direct calculation, we get 
\begin{widetext}
\bea
\frac{\partial \mathfrak{J}(t,X,\omega)}{\partial t} & = \Pi_{\s} G_k(X + z_{\c}(\theta_t\omega)) - A_{\s} \mathfrak{J}(t,X,\omega) \\
& \quad - \int_{-\infty}^t e^{(t-s)A_{\s}} \Pi_{\s} D G_k \Big(e^{(s-t)A_{\c}}X + z_{\c}(\theta_s\omega)\Big) e^{(s-t)A_{\c}} A_{\c} X \d s,
\eea
\end{widetext}
and that
\beas
D_X & \mathfrak{J}(t,X,\omega)  A_{\c} X =\\
& \int_{-\infty}^t e^{(t-s)A_{\s}} \Pi_{\s} D G_k\Big( e^{(s-t)A_{\c}}X + z_{\c}(\theta_s\omega) \Big)  e^{(s-t)A_{\c}} A_{\c} X \d s,
\eeas
where $DG_k$ denote the Fr\'echet derivative of $G_k$. As a result, $\mathfrak{J}(t,X,\omega)$  solves \eqref{Eq_homoligical}.

\medskip
{\bf Step 3.}
The solution to \eqref{Eq_BF_RPDE} is given explicitly 
\bea
p(s) &= e^{(s-t)A_{\c}} X, \\
q(s,\omega) &= \int_{t-\tau}^s e^{(s-s')A_{\s}} \Pi_{\s} G_k\Big(e^{(s'-s)A_{\c}}X + z_{\c}(s',\omega)\Big) \d s'.
\eea
In particular, the value of $q(s,\omega)$ at $s=t$, denoted by $q_{\tau}(t,X, \omega)$, is given by  
\bes
q_{\tau}(t,X, \omega)= \int_{t-\tau}^t e^{(t-s')A_{\s}} \Pi_{\s} G_k\Big(e^{(s'-t)A_{\c}}X + z_{\c}(s', \omega)\Big) \d s'.
\ees
We get then
\bea \label{Eq_diff_J}
\mathfrak{J}&(t,X,\omega) - q_{\tau}(t,X, \omega)= \\
& \int_{-\infty}^{t-\tau} e^{(t-s)A_{\s}} \Pi_{\s} G_k(e^{(s-t)A_{\c}}X + z_{\c}(s,\omega)) \d s.
 \eea
The boundedness of $\mathfrak{J}$ implies that the $H_\s$-norm of the integral on the RHS of \eqref{Eq_diff_J} converges to zero as $\tau$ goes to infinity,  and \eqref{Eq_PB_limit} follows.
\ep

\section{Solving the sACE and its optimal reduced system}\label{Sec_numACE}
 For the numerical integration of the stochastic Allen-Cahn equation \eqref{eq:Chafee}, we adopt a semi-implicit Euler method, where the nonlinearity $-u^3$ and the noise term $\d \W$ are treated explicitly, and the linear term $\partial_x^2 u + u$ is treated implicitly; the spatial discretization is handled with a pseudo-spectral method; see \cite[Section 6.1]{CLW15_vol2} for a detailed description of this scheme applied to a stochastic Burgers-type equation. We have set the time-step size to be $\delta t=10^{-2}$ and the spatial resolution to be $\delta x = L/201$. 
 
 The optimal reduced system  (Eq.~\eqref{Eq_closure})
is also integrated with a semi-implicit Euler method where the linear term is treated implicitly and the nonlinear terms explicitly, still with $\delta t=10^{-2}$. The ODEs with random coefficients such as \eqref{Eq_for_I}, used to simulate the non-Markovian coefficients $Z^{n}_{\tau^*_n}$ involved in the computation of 
$\Phi_{n}(\tau^*_{n},\bm{X},t)$, are  simply integrated by a forward Euler scheme, with again $\delta t=10^{-2}$. The initial condition \eqref{Eq_init} for integrating the RDE \eqref{Eq_for_I} is computed using the trapezoidal rule.

\section{Computations of the double-fold bifurcation diagram, and jump-induced transitions} \label{Sec_numS_shaped}

\subsection{Computation of the double-fold bifurcation diagram} \label{Sec_Num_bif_diag}
To compute the double-fold  bifurcation diagram shown in Figure \ref{Fig_Sbif} for Eq.~\eqref{Eq_Sbif}, we adopt a Fourier sine series expansion of the solutions, in which $u$ solving Eq.~\eqref{Eq_Sbif}  is approximated by $u_N(x) = \sum_{n=1}^N a_n \sqrt{2/L}\sin(n\pi x/L)$. By doing so, the elliptic boundary value problem  \eqref{Eq_Sbif} is reduced to a nonlinear algebraic system for the coefficients $\bm{a}_N = (a_1, \ldots, a_N)^T$, which is solved by the Matlab built-in solver \texttt{fsolve}. The results shown in Figure \ref{Fig_Sbif} are obtained from a six-dimensional algebraic system (i.e., with six sine modes, $N=6$), which turns out to provide high-precision approximations.

We use a simple continuation approach to compute the bifurcation diagram based on this algebraic system, benefiting from the knowledge that the solution curve is S-shaped. Starting from $\lambda = 1$ with a random initial guess, an increment of $\Delta \lambda = 10^{-3}$ is used to compute the lower branch $\lambda\mapsto U^m_\lambda$ of solutions to Eq.~\eqref{Eq_Sbif} (lower black curve) with minimal energy. In our continuation procedure, the initial guess for  $U^m_{\lambda+\Delta \lambda}$  is taken to be  $U^m_\lambda$ computed at the previous step. As $\lambda$ is further increased, the procedure is stopped when the associated Jacobian becomes nearly singular which indicates that $\lambda$ is approaching the lower turning point of the S-shaped curve.

  The location of this turning point is estimated to be $\lambda^\ast\approx 1.3309$ as obtained by using finer mesh as one gets close to singularity. We then select a $\lambda$-value to be just a few $\Delta \lambda$-increments below  $\lambda^\ast$ (to ensure multiplicity of solutions) and pick as many different random initial guesses as needed to get the two other solutions: the unstable one, $U_\lambda^\ast$, and the stable one, $U_\lambda^M$, with maximal energy. From each of these newly computed solutions, we then adopt the same continuation method to trace out the upper branch $\lambda\mapsto U^M_\lambda$  made of solutions with maximal energy (upper black curve)  and the ``in between" branch $\lambda\mapsto U_\lambda^\ast$ made of unstable solutions (red curve).

\subsection{Solving the jump-driven dynamics in Eq.~\eqref{Eq_fluc}, and its optimal reduced model} \label{Sec_Num_S-shaped}

The numerical integration of the stochastic equation \eqref{Eq_fluc} is performed in the same way as for the sACE outlined in Appendix~\ref{Sec_numACE} using a semi-implicit Euler method. Here also the nonlinearity and the noise term are treated explicitly, while the linear term is treated implicitly. The spatial discretization is handled with a pseudo-spectral method; see again \cite[Section 6.1]{CLW15_vol2}. We have set the time-step size to be $\delta t=10^{-2}$ and the spatial resolution to be $\delta x = L/257$. The initial data is set to be $u_0 = 0.5 \bm{e}_1$ for the results shown in Figs.~\ref{Fig_Jump_Dynamics} and \ref{Fig_y1timeseries}.  
 
 The optimal reduced system (Eq.~\eqref{Eq_reduced_S_shaped}) is integrated with a forward Euler method still with $\delta t=10^{-2}$. The initial condition $J^{n}_{\tau^*_n}(t_0; \omega)$ (cf.~\eqref{Eq_J_inital}) for integrating \eqref{Eq_for_J} is computed using the trapezoidal rule, and $t_0=\max\{\tau_3^*, \tau_5^*\}$. The initial data for Eq.~\eqref{Eq_reduced_S_shaped} is set to be the projection of the SPDE solution onto $\bm{e}_1$ at time $t = t_0$.

\bibliography{reference}
\end{document}